\documentclass[12pt]{article}
\usepackage[english]{babel}
\usepackage{textcomp}
\usepackage{amscd}
\usepackage{amssymb}
\usepackage{amsthm}
\usepackage{amsmath}
\usepackage{tikz}
\usepackage{wrapfig}
\usepackage{epsfig}

\newtheorem{theorem}{Theorem}

\newtheorem{proposition}{Proposition}
\newtheorem{lemma}{Lemma}

\newcommand{\Z}{\mathbb{Z}}

\newcommand{\C}{\mathbb{C}}

\newcommand{\n}{m}
\newcommand{\DD}{\mathbb{D}_2}
\newcommand{\midg}{mid grid point}
\newcommand{\tha}{\theta}
\newcommand{\rK}{\widetilde{K}}
\newcommand{\Sb}{\overline{S}}

\newcommand{\mm}{\rho}
\newcommand{\mirror}{\pi}

\newcommand{\me}{\mathrm{e}}
\newcommand{\eight}{{\footnotesize \boxminus}}
\begin{document}

\title{\bf Paperfolding morphisms, planefilling curves, and fractal tiles}

\author{Michel Dekking}

\maketitle

\begin{center}
{\small Michel Dekking, Delft Institute of Applied Mathematics,\\
 Technical University of Delft,  The Netherlands}
{ \tiny F.M.Dekking{@}tudelft.nl}
\end{center}

\begin{abstract}
An interesting class of automatic sequences emerges from iterated paperfolding. The sequences generate curves in the plane with an almost periodic structure. We generalize the results obtained by Davis and Knuth on the self-avoiding and planefilling properties of these curves, giving simple geometric criteria for a complete classification. Finally, we show how the automatic structure of the sequences leads to self-similarity of the curves, which turns the planefilling curves in a scaling limit into fractal tiles. For some of these tiles we give a particularly simple formula for the Hausdorff dimension of their boundary.
\end{abstract}

\begin{center}
\end{center}

\section{Introduction}

Folding starts with folding in two. In the early sixties John E.~Heighway discovered that the pattern of creases that results from repeatedly folding a piece of paper in half generates an intriguing curve in the plane when the creases are unfolded to $90^o$ angles. This curve is known as Heighway's dragon curve.

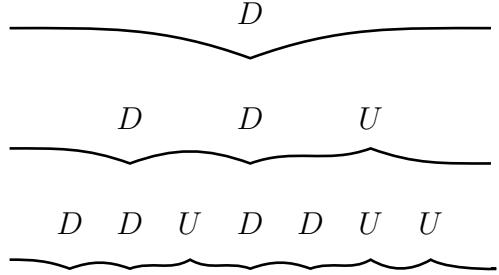
\begin{figure}[!t]
\centering
\begin{tikzpicture}[scale=.4]
% draw grid and axes
%%    \draw[black!20] (-2,-2) grid (18,10);
    \draw[line width=1pt] (0,8) to [out=0, in=160] (8,7) to [out=20, in=180] (16,8);\draw[] (8,8.5) node[] {$D$};
    \draw[line width=1pt] (0,4) to [out=0, in=160] (4,3.5) to [out=20, in=160] (8,3.5) to [out=20, in=200] (12,4) to [out=340, in=180] (16,3.5);
    \draw[] (4,5) node[] {$D$};\draw[] (8,5) node[] {$D$};\draw[] (12,5) node[] {$U$};
    \draw[line width=1pt] (0,.3) to [out=0, in=160] (2,0) to [out=20, in=160] (4,0) to [out=20, in=210] (6,0.3) to [out=330, in=160] (8,0) to [out=20, in=160] (10,0) to [out=20, in=210] (12,0.3) to [out=330, in=210] (14,0.3) to [out=330, in=0] (16,0);
    \draw[] (2,1.5) node[] {$D$};\draw[] (4,1.5) node[] {$D$};\draw[] (6,1.5) node[] {$U$};\draw[] (8,1.5) node[] {$D$};\draw[] (10,1.5) node[] {$D$}; \draw[] (12,1.5) node[] {$U$};\draw[] (14,1.5) node[] {$U$};
\end{tikzpicture}
\caption{\small Folding a strip of paper 1, 2, and 3 times, coding the folds.}\label{fold2}
\end{figure}

The abstract sequence of the symbols $D$ and $U$, coding the `down' folds and `up' folds (cf.~Figure \ref{fold2}) is known as \emph{the} paperfolding sequence\footnote{We follow the convention of \cite{Allou}: when $A$ is an alphabet, elements from $A^*$ will be called words or strings, elements from $A^{\infty}$ will be called infinite words or sequences.}, and is extensively treated in Allouche and Shallit's monograph~\cite{Allou}.

 On the monoid $\{D,U\}^*$  we have the morphism $\mirror$ which permutes letters:
$$\mirror(D)=U, \quad \mirror(U)=D,$$
and the anti-morphism $\tau$ which reverses  words
$$\tau(w_1w_2\dots w_\n)=w_\n\dots w_2w_1 \quad {\rm for\; all}\;w_1w_2\dots w_\n\in \{D,U\}^*.$$
The process of iterated folding is allied to the composition of these two morphisms, the morphism $\mirror \tau$. We therefore shortly write
$$\Sb:=\mirror \tau(S)\quad \text{ for any word } S \text{ in } \{D,U\}^*.$$
We next define the \emph{folding convolution} $S*T$ of two words $S=a_1\ldots a_{s-1}$ and $T=b_1\ldots b_{t-1}$ by
\begin{eqnarray}\label{eq:conv}
S*T = \left\{
\begin{array}{ll}
  Sb_1\Sb b_2S\ldots Sb_{t-1}\Sb& \textrm{ if } t  \textrm{ is even} \\
  Sb_1\Sb b_2S\ldots \Sb b_{t-1}S& \textrm{ if } t \textrm{ is odd.}
  \end{array} \right.
\end{eqnarray}
It is quickly checked that this operation is associative:
$$ (R \ast S)*T=R*(S*T) \quad\textrm{for all } R,S,T \in \{D,U\}^*.$$
Thus $n$-fold (sic) convolutions $S^{*n}$ are unambiguously defined, and since
$S^{*(n+1)}$ starts with $S^{*n}$ for all $n$ we can also talk about the infinite
folding sequence $S^{*\infty}$.

Note that the  three strings in Figure~\ref{fold2} are respectively $D$, $D*D$, and $D*D*D$.
It can be seen in general that folding first according to $T$ and then according to $S$ gives the pattern coded by $S*T$.

We  end this section with a few historic remarks. The  Highway dragon was brought to a larger audience by Martin Gardner in 1967 (\cite{Gardner}). Mathematical proofs of its properties were given by Davis and Knuth three years later in \cite{DK} (reprinted in 2010 in \cite{Fun}). Around 1975 the present author generalized  the notion of folding by introducing the folding convolution  in Equation (\ref{eq:conv}). He announced this in two letters send to  Martin Gardner. These were found by Donald Knuth last year in the Stanford archive. Knuth incited the author to write down the results with proofs.  Some of these results have been announced in the report \cite{Nijm} and the paper \cite{FoldsIII} jointly written with Mend\`es France and van der Poorten.

\section{Folding  morphisms}

Since the main interest in folding sequences is in the curves they generate in the (complex) plane, it
is more convenient to describe these curves directly. We do this by passing from an alphabet of
folding instructions $\{U,D\}$ to an alphabet $A=\{a,b,c,d\}$ of four letters, representing line segments in the positive real, positive imaginary, negative real, and negative imaginary direction.

It does, indeed, appear that the process of iterated folding and unfolding to $90^o$ angles can be described by iterating a uniform morphism on the alphabet $\{a,b,c,d\}$. For 2-folding this morphism $\tha$ is given by
\begin{equation*}\label{2-fold}
\tha(a)=ab,\quad
\tha(b)=cb,\quad
\tha(c)=cd,\quad
\tha(d)=ad.
\end{equation*}
How can this be seen?

Define the morphism $\sigma$ on the monoid $A^*$ which cyclically permutes the letters:
\begin{equation*}\label{sigma}
\sigma(a)=b,\quad
\sigma(b)=c,\quad
\sigma(c)=d,\quad
\sigma(d)=a.
\end{equation*}
Also let the anti-morphism $\tau$, which reverses the words from $A^*$\!, be as before.

The crux of folding morphisms is that they commute with the morphism $\sigma\tau$. For the 2-folding morphism this holds for the letters:
\begin{equation*}\label{2-fold2}
\tha\sigma\tau(a)=\tha(b)=cb=\sigma\tau(ab)=\sigma\tau\tha(a),\quad
\end{equation*}
similarly for the other three letters, and then for any word.

Suppose now that $\tha^n(a)$ corresponds to the curve obtained by folding in two $n$ times.
Then if we unfold a strip folded $n+1$ times, the first $2^n$ letters are equal to those we obtained before,
and the last $2^n$ letters are equal to the first  $2^n$, but in reverse order and cyclically permuted
(because of the $90^o$ angle that we make half way). In other words, this word should be
$\tha^n(a)\sigma\tau\tha^n(a).$

But $\tha^n$ also commutes with $\sigma\tau$, hence $$\tha^n(a)\sigma\tau\tha^n(a)=\tha^n(a)\tha^n\sigma\tau(a)=\tha^n(a)\tha^n(b)=\tha^n(ab)=\theta^{n+1}(a),$$
so that indeed it follows by induction that the morphism $\tha$ describes the iterated 2-folding.

We define in general a \emph{folding morphism} as any morphism $\tha$ on  $A^*$ that commutes with $\sigma\tau$.
Such morphisms are of course determined by their value at the letter $a$, and we will require that
$\tha(a)$ starts with the word $ab$.

 We then have a one-to-one correspondence between strings $S=D\ldots$ of length $m-1$ and folding $m$-morphisms. To see this, define the homomorphism $w:\{D,U\}^*\rightarrow\Z$ by
 \begin{equation}\label{eq:defw}
 w(D)=+1, \quad w(U)=-1.
 \end{equation}
 Then note that if $S=a_1\ldots a_{m-1}$, and $\tha(a)=e_0\ldots e_{m-1}$, then $e_0=a$ and
 \begin{equation}\label{eq:homs}
 e_i=\sigma^{w(a_1\ldots a_i)}(a),\quad \text{for } i=1,\dots,m-1,
 \end{equation}
 where $\sigma$ is the cyclic permutation morphism on $\{a,b,c,d\}^*$.

 Conversely, define the 2-block map $\mm$ by
$$\mm(e\,\sigma(e))=D, \quad  \mm(\sigma(e)\,e)=U \quad \text{for} \; e\in A.$$
Then $\mm(e_0\ldots e_{m-1})=\mm(e_0e_1)\mm(e_1e_2)\ldots\mm(e_{m-2}e_{m-1})$ gives the string
back.

\medskip

Our one-to-one correspondence is an operation-preserving isomorphism if we consider the convolution operation $\ast$ on strings and composition of morphisms:

\begin{lemma}\label{lem:connect}
Let $S=D\ldots$ and $T=D\ldots$ be two strings with associated folding morphisms $\tha_S$ and $\tha_T$ . Then  $$\tha_{S*T}=\tha_S\tha_T.$$
\end{lemma}

\noindent {\bf Proof:} We give a purely geometrical
%\fbox{\tiny How would the algebraic argument with $w$ and $\rho$ go?}
argument: $S*T$ describes the segments in the curve associated to first folding according to $T$, and then to $S$. But the same curve is described by replacing the segments in direction $i^k$ ($k=0,1,2$ or 3) of the curve described by $\tha_T(a)$ by the (scaled) curves  associated to $(\sigma\tau)^k\tha_S(a)$.

The algebraically inclined reader may provide an alternative proof, using equations (\ref{eq:defw}) and (\ref{eq:homs}).\hfill\qed

\medskip

In the sequel \emph{all} strings $S$ from $\{D,U\}^*$ will start with the letter $D$, and $\tha$ will always be the morphism on $A^*=\{a,b,c,d\}$ associated to $S$--- and so $\tha(a)$ will always start with $ab$.

\section{Representations of morphisms}\label{sec:rec}

We now do a more elaborate analysis of the situation based on (\cite{Rec}).

As usual we will write for a set $B\subset \C$, and $r,t\in \C$:
$$rB=\{rx: x\in B\}\quad {\rm and} \quad B+t=\{x+t: x\in B\},$$
and we will write \emph{segments} in the complex plane as
$$[y,y']=\{x\in \C: x=\alpha y + (1-\alpha)y',\; 0\le \alpha\le 1\}.$$
We introduce a homomorphism $f:A^*\rightarrow \C$ by
\begin{equation}\label{eq:f}
f(a)=1,\; f(b)=i,\; f(c)=-1,\; f(d)=-i,
\end{equation}
where the operation on $A^*$ is concatenation, and that on $\C$ addition.
To associate curves to words we define a map $K[\cdot]:A^*\rightarrow 2^\C$, by first putting
$$K[a]=[0,1],\; K[b]=[0,i],\; K[c]=[0,-1],\; K[d]=[0,-i],
$$
and then extending this by putting for any two words $v,w\in A^*$
\begin{equation}\label{eq:Ktha}
K[vw]=K[v]\cup (K[w]+f(v)).
\end{equation}
So for example $K[abcd]=K[a]\cup(K[b]+1)\cup(K[c]+1\!\!+\!i)\cup(K[d]+i)$ is nothing else but the unit square in $\C$.

For some morphisms $\tha$ we will have a commuting diagram
%\begin{equation*}\label{commute}
%\begin{array}{clllllllllllll}
%\quad     & f   & \quad       & \\[.3cm]
%A^*       &:    & \rightarrow &\C\\[.3cm]
%\downarrow&\tha & \downarrow  &L_\tha\\
%A^*       &:    & \rightarrow &\C
%\end{array}
%\end{equation*}
$$\begin{CD}
A^* @>\tha>>A^*\\
@VVf V @VV fV\\
\C @>L_\tha >> \C
\end{CD}$$
where $L_\tha$ is multiplication by a complex number.
It is not hard to see that $L_\tha$ exists for any folding morphism and is multiplication by $z=f(\tha(a))$. For example for the 2-folding morphism determined by $\tha(a)=ab$, $L_\tha$ is multiplication by $1+i$. As another example: $L_\sigma$ is multiplication by $i$ for the cyclic permutation morphism $\sigma$.

\medskip

Folding convolutions behave nicely with respect to the embedding in $\C$:

\begin{lemma}\label{endstar}
Let $S$ and $T$ be two folding strings, whose curves end in $z_S$ and $z_T$. Then
the curve of $S*T$ ends in $z_Sz_T$.
\end{lemma}

\noindent {\bf Proof:} We have that $L_{\tha_S}$ is multiplication by
$z_S=f(\tha_S(a))$,  and $L_{\tha_T}$ is multiplication by $z_T=f(\tha_T(a))$. Lemma~\ref{lem:connect} implies that the curve of  $S*T$ ends in the point $f(\tha_S\tha_T(a))$.
The result thus follows from a look at the following commutative diagram, recalling that $f(a)=1$.
$$\begin{CD}
A^* @>\tha_T>>A^*@>\tha_S>>A^*\\
@VVfV @VVfV @VVfV\\
\C @>L_{\tha_T} >> \C @>L_{\tha_S}>>\C
\end{CD}$$\hfill \qed

In Sections \ref{sec:planefill} and  \ref{sec:classify} we will need a map from unions of Gaussian integer segments to itself that describes the action of $\tha$.
We consider neigbouring Gaussian integers $z$ and $z'$, defining the $A$-code $e_{z,z'}$ of the segment $[z,z']$ by
$$e_{z,z'}=e\quad{\rm if}\quad z'-z=f(e), \quad  \text{where}\;e=a,b,c,\,\text{or}\;d.$$
If $[z,z']$ is such a segment then we define the \emph{$K$-$\tha$-image} of $[z,z']$ by
\begin{equation}\label{eq:Kthimage}
K[\tha]([z,z'])=K[\tha(e_{z,z'})]+L_\tha(z).
\end{equation}
For example, if $\tha$ is the 2-folding morphism, then
\begin{eqnarray*}
K[\tha]([1+i,i])&=&K[\tha(c)]+L_\tha(1+i)=K[cd]+2i\\
    &=&(K[c]\,\cup\,(K[d]-1))+2i=[2i,2i-1]\,\cup\,[2i-1,i-1].
\end{eqnarray*}
The map $K[\tha](\cdot)$ is extended to
 unions of $k$ segments, in the natural way:
\begin{equation}\label{eq:Kthimageunion}
K[\tha]([z_1,z'_1]\,\cup\ldots\cup\,[z_k,z'_k])=K[\tha]([z_1,z'_1])\,\cup \ldots\cup\,K[\tha]([z_k,z'_k]).
\end{equation}

Note that if the union $Q$ of the segments can be written as $Q=K[w]$ for some word $w$ from $A^k$, then by (\ref{eq:Ktha}) simply $K[\tha](Q)=K[\tha(w)]$.

\section{Symmetric folds}

We call a string of $U$'s and $D$'s \emph{symmetric} if $S=\Sb$. For example the string $DDUU$ coding  a 5-fold is symmetric. Obviously only strings of folds of odd length can be symmetric. We call a morphism $\tha$ \emph{palindromic} if it commutes with the reversal map, i.e., if $\tau\tha=\tha\tau$.

\begin{lemma}\label{lem:symm}
A  string is symmetric if and only if the associated folding morphism is palindromic.
\end{lemma}

\noindent {\bf Proof:} Suppose $S=a_1\ldots a_{m-1}$ is symmetric. It is easy to see that we only have to check that the associated $\tha(a)=e_0\ldots e_{m-1}$ is a palindrome word. First note that the symmetry of $S$ implies that the number of $U$'s in $S$ equals the number of $D$'s, and hence that $w(S)=0$, where $w$ is the homomorphism from (\ref{eq:defw}). But then we have for $j=1,\dots,m-1$
$$w(a_1\ldots a_{m-j-1})=-w(a_{m-j}\ldots a_{m-1})=w(a_1\ldots a_j),$$
since symmetry is equivalent to $w(a_{m-i})=-w(a_i)$ for $i=1,\dots m-1$.
Then we have by Equation (\ref{eq:homs}) that
$$e_{m-j-1}=\sigma^{w(a_1\ldots a_{m-j-1})}(a)=\sigma^{w(a_1\ldots a_j)}(a)=e_j\quad \text{for } j=1,\dots,m-1,$$
which exactly means that $\tha(a)=e_0\ldots e_{m-1}$ is a palindrome.

Conversely, if $\tha(a)$ is a palindrome, then for $i=1,\dots,m-1$
$$a_i=\rho(e_{i-1}e_{i})=\rho(e_{m-i}e_{m-i-1})=\overline{\rho(e_{m-i-1}e_{m-i})}=\overline{a_{m-i}},$$
and thus $S$ is symmetric. \hfill \qed

\section{Self-avoiding folding sequences}

We call a folding string $S$ \emph{self-avoiding} if  the curve associated to the sequence $S^{*\infty}$ does not traverse a segment twice. It will appear that all information that we need to decide whether this infinite sequence satisfies this property is in the finite word $\tha(abcd)$, if $\tha$ is the morphism associated to $S$.

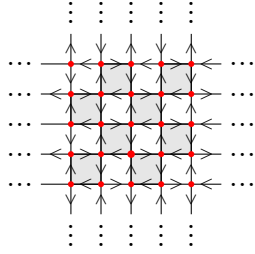
\begin{figure}[!h]
\centering
%\vspace*{-0.5cm}
\begin{tikzpicture}[scale=.4]
% Filled blocks with a light gray color
    \begin{scope}
        \filldraw[fill=black!10](-1,2)--(0,2)--(0,3)--(-1,3)--cycle;
        \filldraw[fill=black!10](1,2)--(2,2)--(2,3)--(1,3)--cycle;
        \filldraw[fill=black!10](-2,1)--(-1,1)--(-1,2)--(-2,2)--cycle;
        \filldraw[fill=black!10](0,1)--(1,1)--(1,2)--(0,2)--cycle;
        \filldraw[fill=black!10](-1,0)--(0,0)--(0,1)--(-1,1)--cycle;
        \filldraw[fill=black!10](1,0)--(2,0)--(2,1)--(1,1)--cycle;
        \filldraw[fill=black!10](-2,-1)--(-1,-1)--(-1,0)--(-2,0)--cycle;
        \filldraw[fill=black!10](0,0)--(1,0)--(1,-1)--(0,-1)--cycle;
    \end{scope}
% draw matrix
    % draw horizontal lines
    \foreach \x in {-2,-1,0,1,2} {\draw[] (\x,-2) --  (\x,4);}
    % draw vertical lines
    \foreach \y in {-1,0,1,2,3} {\draw[] (-3,\y) -- (3,\y);}
    % draw arrows "<" and ">"
    \foreach \y in {-1,1,3} {
        \draw[] (-1.5,\y) node[] {\tiny$<$};\draw[] (.5,\y) node[] {\tiny$<$};\draw[] (2.5,\y) node[] {\tiny$<$};
        \draw[] (-.5,\y) node[] {\tiny$>$};\draw[] (1.5,\y) node[] {\tiny$>$};}
    \foreach \y in {0,2} {
        \draw[] (-2.5,\y) node[] {\tiny$<$};\draw[] (-.5,\y) node[] {\tiny$<$};\draw[] (1.5,\y) node[] {\tiny$<$};
        \draw[] (-1.5,\y) node[] {\tiny$>$};\draw[] (.5,\y) node[] {\tiny$>$};\draw[] (2.5,\y) node[] {\tiny$>$};}
    % draw arrows "\vee" and "\wedge"
    \foreach \x in {-2,0,2} {
        \draw[] (\x,-1.5) node[] {\tiny$\vee$};\draw[] (\x,.5) node[] {\tiny$\vee$};\draw[] (\x,2.5) node[] {\tiny$\vee$};
        \draw[] (\x,-.5) node[] {\tiny$\wedge$};\draw[] (\x,1.5) node[] {\tiny$\wedge$};\draw[] (\x,3.5) node[] {\tiny$\wedge$};}
    \foreach \x in {-1,1} {
        \draw[] (\x,-.5) node[] {\tiny$\vee$};\draw[] (\x,1.5) node[] {\tiny$\vee$};\draw[] (\x,3.5) node[] {\tiny$\vee$};
        \draw[] (\x,-1.5) node[] {\tiny$\wedge$};\draw[] (\x,.5) node[] {\tiny$\wedge$};\draw[] (\x,2.5) node[] {\tiny$\wedge$};}
    % draw 3 dots vertical lines
    \foreach \x in {-2,-1,0,1,2} {\node at (\x,-3) {.}; \node at (\x,-2.7) {.}; \node at (\x,-2.4) {.};
        \node at (\x,4.4) {.}; \node at (\x,4.7) {.}; \node at (\x,5) {.};     }
    % draw 3 dots horizontal lines
    \foreach \y in {-1,0,1,2,3} {
        \node at (-4,\y) {.};   \node at (-3.7,\y) {.};
        \node at (-3.4,\y) {.}; \node at (3.4,\y) {.};
        \node at (3.7,\y) {.};  \node at (4,\y) {.};    }
% Drawing the red points
    \foreach \x in {-2,-1,0,1,2}{
        \foreach \y in {-1,0,1,2,3}{\node[fill=red,circle, inner sep=.8pt] at (\x,\y)  {}; }  }
% Null-point
    \node[fill=red,circle, inner sep=1pt] at (0,0)  {};
\end{tikzpicture}
\caption{\small The directed graph $\DD$.}\label{fig:D2}
\end{figure}

It is important to observe that all words generated by a folding morphism have the property that letters from $\{a,c\}$ alternate with letters from $\{b,d\}$. We will call such words \emph{$\DD$-words}, for reasons that will become clear soon (see also Figure~\ref{fig:D2}).

\begin{figure}[!h]
\centering
\begin{tikzpicture}[scale=.35]
\draw[blue]   (0., 0.)-- (1., 0.)-- (1., 1.)--(0., 1.)-- (0., 2.)--(-1., 2.)--(-1., 1.)--(-2., 1.)--
(-2., 2.)--  (-3., 2.)--(-3.,1.)--(-2., 1.)--(-2., 0.)--  (-3., 0.)--(-3., -1.)--  (-4., -1.)--
(-4., 0.)--  (-5., 0.)--(-5., -1.)--  (-4., -1.)--(-4., -2.)--  (-3., -2.)--(-3., -1.)--  (-2., -1.)--
(-2., -2.)--  (-3., -2.)--(-3., -3.)--  (-2., -3.)--(-2., -4.)--  (-3., -4.)--(-3., -5.)--  (-4., -5.)--
(-4., -4.)--  (-5., -4.)--(-5., -5.)--  (-4., -5.)--(-4., -6.)--  (-3., -6.)--(-3., -5.)--  (-2., -5.)--
(-2., -6.)--  (-1., -6.)--(-1., -5.)--  (-2., -5.)--(-2., -4.)--  (-1., -4.)--(-1., -3.)--  (0., -3.)--
(0., -4.)--  (-1., -4.)--(-1., -5.)--  (0., -5.)--(0., -6.)--  (1., -6.)--(1., -5.)--  (2., -5.)--
(2., -6.)--  (1., -6.)--(1., -7.)--  (2., -7.)--(2., -8.)--  (1., -8.)--(1., -9.)--(0., -9.)--(0,-8);
\node[fill=red,circle, inner sep=1pt] at (0,0)  {};
\end{tikzpicture}\hspace*{2cm}
\begin{tikzpicture}[scale=.35,rounded corners=2pt]
\draw[blue]  (0., 0.)--(0.750, 0.)--(1., 0.250)--(1., 0.750)--(0.750, 1.)--(0.250, 1.)--
 (0., 1.25)--  (0., 1.75)--(-0.250, 2.)--(-0.750, 2.)--(-1., 1.75)-- (-1., 1.25)--(-1.25, 1.)-- (-1.75, 1.)--
 (-2., 1.25)-- (-2., 1.75)--( -2.25, 2.)-- (-2.75, 2.)--(-3., 1.75)-- (-3., 1.25)--(-2.75, 1.)-- (-2.25, 1.)--
 (-2., 0.750)--(-2., 0.250)--(-2.25, 0.)-- (-2.75, 0.)--(-3.,-0.250)--(-3.,-0.750)--(-3.25, -1.)--(-3.75, -1.)--
(-4., -0.750)--(-4., -0.250)--(-4.25, 0.)-- (-4.75, 0.)--(-5.,-0.250)--(-5.,-0.750)--(-4.75, -1.)--(-4.25, -1.)--
 (-4., -1.25)-- (-4., -1.75)--(-3.75, -2.)--(-3.25, -2.)--(-3.,-1.75)--(-3., -1.25)--(-2.75, -1.)--(-2.25, -1.)--
 (-2., -1.25)-- (-2., -1.75)--(-2.25, -2.)--(-2.75, -2.)--(-3.,-2.25)--(-3., -2.75)--(-2.75, -3.)--(-2.25, -3.)--
 (-2., -3.25)-- (-2., -3.75)--(-2.25, -4.)--(-2.75, -4.)--(-3.,-4.25)--(-3., -4.75)--(-3.25, -5.)--(-3.75, -5.)--
 (-4., -4.75)-- (-4., -4.25)--(-4.25, -4.)--(-4.75, -4.)--(-5.,-4.25)--(-5., -4.75)--(-4.75, -5.)--(-4.25, -5.)--
 (-4., -5.25)-- (-4., -5.75)--(-3.75, -6.)--(-3.25, -6.)--(-3.,-5.75)--(-3., -5.25)--(-2.75, -5.)--(-2.25, -5.)--
 (-2., -5.25)-- (-2., -5.75)--(-1.75, -6.)--(-1.25, -6.)--(-1.,-5.75)--(-1., -5.25)--(-1.25, -5.)--(-1.75, -5.)--
 (-2., -4.75)-- (-2.,-4.25)--(-1.75, -4.)--(-1.25,-4.)--(-1.,-3.75)--(-1., -3.25)--(-0.750, -3.)--(-0.250, -3.)--
 (0., -3.25)-- (0.,-3.75)--(-0.250, -4.)--(-0.750,-4.)--(-1., -4.25)--(-1.,-4.75)--(-0.750, -5.)--(-0.250,-5.)--
 (0., -5.25)-- (0., -5.75)-- (0.250, -6.)-- (0.750, -6.)-- (1., -5.75)--(1., -5.25)--(1.25, -5.)-- (1.75, -5.)--
 (2., -5.25)-- (2., -5.75)-- (1.75, -6.)-- (1.25, -6.)--  (1., -6.25)--(1., -6.75)--(1.25, -7.)-- (1.75, -7.)--
 (2., -7.25)-- (2., -7.75)-- (1.75, -8.)-- (1.25, -8.)-- (1., -8.25)-- (1., -8.75)--(0.750, -9.)--
 (0.250, -9.)--(0., -8.75)--  (0., -8.);
% \node  at (0,0) {{\Huge $\cdot$}};
\node[fill=red,circle, inner sep=1pt] at (0,0)  {};
\end{tikzpicture}
\caption{\small The curves $K[w]$ and $\rK[w]$ with $w=\tha^6(a)$, where $\tha(a)=ab$ .}\label{fig:round}
\end{figure}
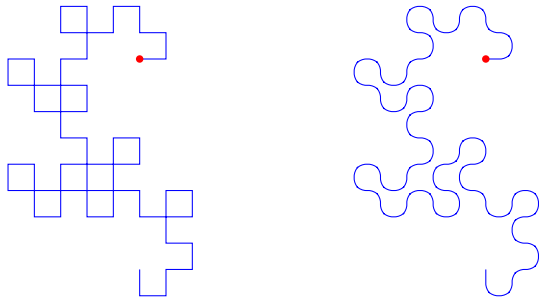

For the following development it is convenient  to consider curves $\rK[w]$, obtained from $K[w]$ by rounding of the corners (see Figure~\ref{fig:round}).

The point of this is that the curves $\rK[w]$
do not pass through the Gaussian integers $u+iv$, with $u,v\in \Z$. In this way, if $w$ is a $\DD$-word, and $K[w]$ does not traverse a segment twice, then $\rK[w]$ is a simple curve in the topological sense (i.e., the image of an injective map from an interval to the plane).

 \begin{figure}[b]
\centering
 \begin{tikzpicture}[scale=.5,rounded corners=2pt]
\draw[blue]        (0,0)--(0.50, 0.)--(0.75, 0.25) --(0.75, 0.75) -- (1., 1.) --(1.5, 1.) --
              (1.75, 1.25) --(1.75, 1.75) --(2., 2.) -- (2.5, 2.) --(2.75, 1.75) --(2.75, 1.25) --
          (3., 1.) --(3.5, 1.) --(3.75, 1.25) --(3.75, 1.75)--(3.5, 2.) --(3., 2.) -- (2.75, 2.25) --
           (2.75, 2.75) --(2.5, 3.) --(2., 3.) --(1.75, 3.25) --(1.75, 3.75) --(1.5, 4.)--(1., 4.) --
          (0.75, 3.75) --(0.75, 3.25) --(0.50, 3.) --(0., 3.) --(-0.25, 2.75) --(-0.25, 2.25) --
           (-0.50, 2.)--(-1., 2.) --(-1.25, 2.25) --(-1.25, 2.75)--(-1.5, 3.) --(-2., 3.)--(-2.25, 2.75) --
 (-2.25, 2.25) --(-2., 2.)--(-1.5, 2.) --(-1.25, 1.75)--(-1.25, 1.25)--(-1., 1.)--(-0.50, 1.)--(-0.25, 0.75) --
(-0.25, 0.25) -- (0., 0.);
% \node  at (0,0) {{\Huge $\cdot$}};
\node[fill=black,circle, inner sep=.8pt] at (2.75,1.0) {};
\node[fill=black,circle, inner sep=.8pt] at (1.75,4.0) {};
\node[fill=black,circle, inner sep=.8pt] at (-1.25,3.0) {};
\node[fill=red,circle, inner sep=.8pt] at (-0.25,0)  {};
\end{tikzpicture}
\caption{\small The $\tha$-loop of $\tha(a)=ababad$.}\label{fig:loop}
\end{figure}
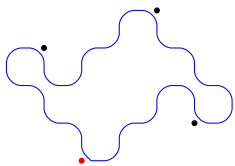

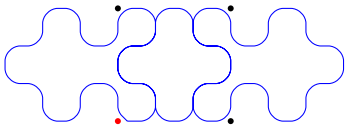
\begin{figure}[!h]
\centering
\begin{tikzpicture}[scale=.5,rounded corners=2pt]
\draw[blue]
  ( 0., 0.)-- ( 0.500, 0.)--(0.750,0.250)--(0.750, 0.750)--( 0.500, 1.)--   ( 0., 1.)--(-0.250,1.25)--(-0.250, 1.75)--
  ( 0., 2.)-- ( 0.500, 2.)--( 0.750, 2.25)--(0.750, 2.75)--  ( 1., 3.)--  ( 1.50, 3.)-- ( 1.75, 2.75)--( 1.75, 2.25)--
  ( 2., 2.)--  ( 2.50, 2.)-- ( 2.75, 1.75)--( 2.75, 1.25)-- ( 2.50, 1.)--   ( 2., 1.)--( 1.75, 0.750)--( 1.75,0.250)--
  ( 2., 0.)--  ( 2.50, 0.)--( 2.75, 0.250)--(2.75, 0.750)--  ( 3., 1.)--  ( 3.50, 1.)--( 3.75, 0.750)--(3.75, 0.250)--
  ( 4., 0.)--  ( 4.50, 0.)--( 4.75, 0.250)--(4.75, 0.750)--  ( 5., 1.)--  ( 5.50, 1.)-- ( 5.75, 1.25)--( 5.75, 1.75)--
 ( 5.50, 2.)--   ( 5., 2.)-- ( 4.75, 2.25)--( 4.75, 2.75)-- ( 4.50, 3.)--   ( 4., 3.)-- ( 3.75, 2.75)--( 3.75, 2.25)--
 ( 3.50, 2.)--   ( 3., 2.)-- ( 2.75, 2.25)--( 2.75, 2.75)-- ( 2.50, 3.)--   ( 2., 3.)-- ( 1.75, 2.75)--( 1.75, 2.25)--
  ( 2., 2.)--  ( 2.50, 2.)-- ( 2.75, 1.75)--( 2.75, 1.25)-- ( 2.50, 1.)--   ( 2., 1.)--( 1.75, 0.750)--(1.75, 0.250)--
 ( 1.50, 0.)--   ( 1., 0.)--(0.750, 0.250)--(0.750,0.750)--( 0.500, 1.)--   ( 0., 1.)--(-0.250, 1.25)--(-0.250,1.75)--
  ( 0., 2.)-- ( 0.500, 2.)--( 0.750, 2.25)--(0.750, 2.75)--( 0.500, 3.)--   ( 0., 3.)--(-0.250,2.75)--(-0.250, 2.25)--
( -0.500, 2.)-- ( -1., 2.)--( -1.25,2.25)--( -1.25, 2.75)--( -1.50, 3.)--  ( -2., 3.)--( -2.25, 2.75)--(-2.25, 2.25)--
( -2.50, 2.)--  ( -3., 2.)--( -3.25, 1.75)--(-3.25, 1.25)-- ( -3., 1.)-- ( -2.50, 1.)--(-2.25, 0.750)--(-2.25,0.250)--
(-2.,0.)--( -1.50,0.)--(-1.25, 0.250)--(-1.25,0.750)--( -1., 1.)--(-0.500,1.)--( -0.250, 0.750)--(-0.250, 0.250)--(0,0);
\node[fill=red,circle, inner sep=.8pt] at (-0.25,0)  {};
\node[fill=black,circle, inner sep=.8pt] at (-0.25,3)  {};
\node[fill=black,circle, inner sep=.8pt] at (2.75,3)  {};
\node[fill=black,circle, inner sep=.8pt] at (2.75,0)  {};
\end{tikzpicture}
\caption{\small The $\tha$-loop of $\tha(a)=abcbabadadcda$.}\label{fig:13-fold-loop}
\end{figure}

The central notion in the sequel will be that of a \emph{$\tha$-loop} defined for any morphism $\tha$ as the set $\rK[\tha(abcd)]$. See Figure~\ref{fig:loop} and  Figure~\ref{fig:13-fold-loop} for  examples.

\begin{lemma}\label{closed}
Let $\tha$ be a folding morphism. Then its $\tha$-loop is a closed curve.
\end{lemma}

\noindent {\bf Proof:} The endpoint of $K[\tha(abcd)]$ is given by
$$f(\tha(abcd))=L_\tha(f(abcd))=L_\tha (\,f(a)+\!f(b)+\!f(c)+\!f(d)\,)=0.\hfill\qed$$

\begin{theorem}\label{th:avoid}
Let $S$ be a folding string. Then $S$ is a self-avoiding string if and only if the $\theta$-loop of $S$ is a simple curve.
\end{theorem}

\noindent {\bf Proof:} \emph{Part 1:} Suppose that the $\tha$-loop is not simple, i.e., there is a segment that is traced twice. We have to show that there exist $n\ge1$ such that $\rK[\tha^n(a)]$ is not simple.

Let $z=f(\tha(a))$. The $\tha$-loop (ignoring the rounding off) consists of  the four curves
$K_a:=K[\tha(a)], K_b:=K[\tha(b)]+z,K_c:=K[\tha(c)]+(1+i)z,$ and $K_d:=K[\tha(d)]+i\,z,$
which can be mapped to each other by a rotation and a translation (see Lemma~\ref{lem:dir}).
Therefore, if one $K_e$ has a segment that it is traced by another segment of the \emph{same} $K_e$, then all four curves
have such a segment. So the curve $\rK[\tha(a)]$ is already not simple.

The next possibility is that a segment of a $K_e$ is traced by $K_{\sigma(e)}$, and again, if this happens for one $e$, then it happens for all $e\in A$. Thus, since $\tha(a)$ starts with $ab=a\sigma(a)$, the curve $\rK[\tha^2(a)]$ will not be simple.

The last (note that $K_e$ has a segment in common with $K_{\sigma^3(e)}$ iff $K_{\sigma(e)}$ has a segment in common with $K_e$) possibility is that a segment of a $K_e$ is traced by $K_{\sigma^2(e)}$, and now the situation is different, since this might happen for $K_a$ and $K_c$, but \emph{not} for $K_b$ and $K_d$ (cf. Figure~\ref{fig:13-fold-loop}), or vice versa. In the first case this overlap will occur in $K[\tha(abc)]$, in the second case in $K[\tha(bad)]$. So our goal is to show that both
$abc$ and $bad$ will occur in some $\tha^n(a)$.

 Now, since $\tha(a)$ starts with $ab$, there exists a $k\ge 1$ such that \emph{either} $\tha(a)=(ab)^k$ or $(ab)^ka$, \emph{or} $\tha(a)$ starts with $(ab)^kc$ or with $(ab)^kad$. The former case can not occur, since clearly the $\tha$-loops of $(ab)^k$ or $(ab)^ka$ are simple for any $k\ge 1$. In the second case $\tha(a)$ contains the word $abc$, or the word $bad$. But since $\tau\sigma(bad)=abc$ and $(\tau\sigma)^3(abc)=bad$, $\tha^3(a)$ will contain \emph{both} $abc$ and $bad$.
It follows that $\rK[\tha^3(a)]$ will not be simple.

\emph{Part 2:} Suppose the $\theta$-loop of $S$ is a simple curve.
The crucial observation here is that in whatever contorted way a curve $K[\tha(w)]$ (with $w$ a $\DD$-word, and $\tha$ a folding morphism) lies in the plane, it will always be a subset of a very regular checkerboard pattern of $\tha$-loops and anti-$\tha$-loops (ignoring rounding off). Here we define the \emph{anti-$\tha$-loop} as the curve $\rK[\tha(dcba)]$.

The regular pattern is created by  the structure of the directed graph $\DD$, which as vertices has the Gaussian integers, and directed edges between neighbouring vertices as in Figure~\ref{fig:D2}. To avoid cumbersome notation we shall often identify the vertices with the corresponding Gausssian integers (which we will also call \emph{grid points}), and the directed edges  with the corresponding segments in the plane.

We call a Gaussian integer $z=u+iv$ \emph{even} if $u+v$ is even, and \emph{odd} if  if $u+v$ is odd.

 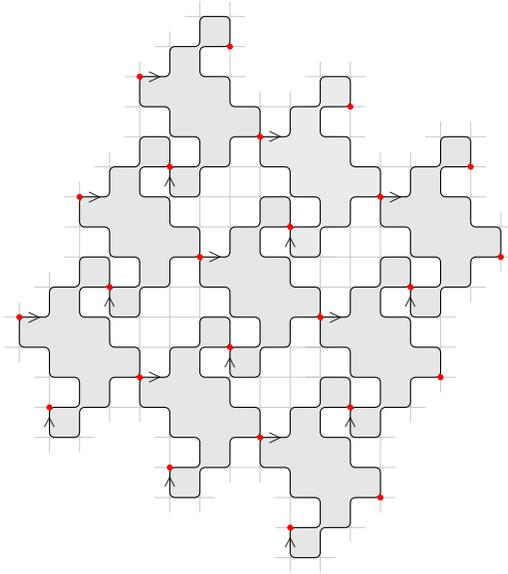
\begin{figure}[!h]
\centering
\begin{tikzpicture}[scale=.4,rounded corners=.8pt]
% draw grid
   \draw[black!20] (-4.5,.5) grid(-3.5,2.5)(-3.5,-2.5) grid (-1.5,3.5)
    (-2.5,3.5) grid (-1.5,6.5) (-1.5,-1.5)grid(-.5,7.5)(-0.5,-1.5) grid (.5,10.5)(.5,-4.5)grid(2.5,11.5)(2.5,-3.5) grid (4.5,9.5)(1.5,9.5)grid(3.5,12.5)(4.5,-6.5) grid (6.5,9.5)(5.5,-5.5)grid(7.5,10.5)(7.5,-4.5)grid(8.5,7.5)(8.5,-1.5) grid(9.5,7.5)(9.5,-.5)grid(10.5,8.5)(10.5,2.5)grid(11.5,8.5)(11.5,3.5)grid(12.5,5.5);
    \filldraw[fill=black!10] {(0.1,0) coordinate (P0) -- (0.90,0) --
          (1,0.10) -- (1,0.9) -- (1.1,1) -- (1.9,1) -- (2,1.1) --
    (2,1.9) -- (2.1,2) -- (2.9,2) --(3,1.9) -- (3,1.1) -- (2.9,1) -- (2.1,1) -- (2,0.9) --
          (2,0.1) -- (2.1,0) -- (2.9,0) -- (3,-0.1) -- (3,-0.9) --
          (3.1,-1) -- (3.9,-1) -- (4,-1.1) -- (4,-1.9) -- (3.9,-2.) --
          (3.1,-2) -- (3,-2.1) -- (3,-2.9) -- (2.9,-3) -- (2.1,-3.) --
          (2,-3.1) -- (2,-3.9) -- (1.9,-4) -- (1.1,-4) -- (1,-3.9) --
          (1,-3.1) -- (1.1,-3) -- (1.9,-3) -- (2,-2.9) -- (2,-2.1) --
          (1.9,-2.) -- (1.1,-2.) -- (1,-1.9) -- (1,-1.1) -- (0.90,-1) --
          (0.10,-1) -- (0,-0.9) -- (0,-0.1) -- (0.1,0)};
    \filldraw[fill=black!10,xshift=-2cm,yshift=6cm] {(0.1,0) coordinate (P1)--
          (0.90,0) -- (1,0.10) -- (1,0.9) -- (1.1,1) --
          (1.9,1) -- (2,1.1) -- (2,1.9) -- (2.1,2) -- (2.9,2) --
          (3,1.9) -- (3,1.1) -- (2.9,1) -- (2.1,1) -- (2,0.9) --
          (2,0.1) -- (2.1,0) -- (2.9,0) -- (3,-0.1) -- (3,-0.9) --
          (3.1,-1) -- (3.9,-1) -- (4,-1.1) -- (4,-1.9) -- (3.9,-2.) --
          (3.1,-2) -- (3,-2.1) -- (3,-2.9) -- (2.9,-3) -- (2.1,-3.) --
          (2,-3.1) -- (2,-3.9) -- (1.9,-4) -- (1.1,-4) -- (1,-3.9) --
          (1,-3.1) -- (1.1,-3) -- (1.9,-3) -- (2,-2.9) -- (2,-2.1) --
          (1.9,-2.) -- (1.1,-2.) -- (1,-1.9) -- (1,-1.1) -- (0.90,-1) --
          (0.10,-1) -- (0,-0.9) -- (0,-0.1) -- (0.1,0)};
    \filldraw[fill=black!10,xshift=-4cm,yshift=2cm] {(0.10,0) coordinate (P2) --
          (0.90,0) -- (1,0.10) -- (1,0.9) -- (1.1,1) --
          (1.9,1) -- (2,1.1) -- (2,1.9) -- (2.1,2) -- (2.9,2) --
          (3,1.9) -- (3,1.1) -- (2.9,1) -- (2.1,1) -- (2,0.9) --
          (2,0.1) -- (2.1,0) -- (2.9,0) -- (3,-0.1) -- (3,-0.9) --
          (3.1,-1) -- (3.9,-1) -- (4,-1.1) -- (4,-1.9) -- (3.9,-2.) --
          (3.1,-2) -- (3,-2.1) -- (3,-2.9) -- (2.9,-3) -- (2.1,-3.) --
          (2,-3.1) -- (2,-3.9) -- (1.9,-4) -- (1.1,-4) -- (1,-3.9) --
          (1,-3.1) -- (1.1,-3) -- (1.9,-3) -- (2,-2.9) -- (2,-2.1) --
          (1.9,-2.) -- (1.1,-2.) -- (1,-1.9) -- (1,-1.1) -- (0.90,-1) --
          (0.10,-1) -- (0,-0.9) -- (0,-0.1) -- (0.1,0)};
    \filldraw[fill=black!10,xshift=4cm,yshift=-2cm] {(0.1,0) coordinate (P3) --
          (0.90,0) -- (1,0.10) -- (1,0.9) -- (1.1,1) --
          (1.9,1) -- (2,1.1) -- (2,1.9) -- (2.1,2) -- (2.9,2) --
          (3,1.9) -- (3,1.1) -- (2.9,1) -- (2.1,1) -- (2,0.9) --
          (2,0.1) -- (2.1,0) -- (2.9,0) -- (3,-0.1) -- (3,-0.9) --
          (3.1,-1) -- (3.9,-1) -- (4,-1.1) -- (4,-1.9) -- (3.9,-2.) --
          (3.1,-2) -- (3,-2.1) -- (3,-2.9) -- (2.9,-3) -- (2.1,-3.) --
          (2,-3.1) -- (2,-3.9) -- (1.9,-4) -- (1.1,-4) -- (1,-3.9) --
          (1,-3.1) -- (1.1,-3) -- (1.9,-3) -- (2,-2.9) -- (2,-2.1) --
          (1.9,-2.) -- (1.1,-2.) -- (1,-1.9) -- (1,-1.1) -- (0.90,-1) --
          (0.10,-1) -- (0,-0.9) -- (0,-0.1) -- (0.1,0)};
    \filldraw[fill=black!10,xshift=2cm,yshift=4cm] {(0.1,0) coordinate (P4) --
          (0.90,0) -- (1,0.10) -- (1,0.9) -- (1.1,1) --
          (1.9,1) -- (2,1.1) -- (2,1.9) -- (2.1,2) -- (2.9,2) --
          (3,1.9) -- (3,1.1) -- (2.9,1) -- (2.1,1) -- (2,0.9) --
          (2,0.1) -- (2.1,0) -- (2.9,0) -- (3,-0.1) -- (3,-0.9) --
          (3.1,-1) -- (3.9,-1) -- (4,-1.1) -- (4,-1.9) -- (3.9,-2.) --
          (3.1,-2) -- (3,-2.1) -- (3,-2.9) -- (2.9,-3) -- (2.1,-3.) --
          (2,-3.1) -- (2,-3.9) -- (1.9,-4) -- (1.1,-4) -- (1,-3.9) --
          (1,-3.1) -- (1.1,-3) -- (1.9,-3) -- (2,-2.9) -- (2,-2.1) --
          (1.9,-2.) -- (1.1,-2.) -- (1,-1.9) -- (1,-1.1) -- (0.90,-1) --
          (0.10,-1) -- (0,-0.9) -- (0,-0.1) -- (0.1,0)};
    \filldraw[fill=black!10,xshift=6cm,yshift=2cm] {(0.1,0) coordinate (P5) --
          (0.90,0) -- (1,0.10) -- (1,0.9) -- (1.1,1) --
          (1.9,1) -- (2,1.1) -- (2,1.9) -- (2.1,2) -- (2.9,2) --
          (3,1.9) -- (3,1.1) -- (2.9,1) -- (2.1,1) -- (2,0.9) --
          (2,0.1) -- (2.1,0) -- (2.9,0) -- (3,-0.1) -- (3,-0.9) --
          (3.1,-1) -- (3.9,-1) -- (4,-1.1) -- (4,-1.9) -- (3.9,-2.) --
          (3.1,-2) -- (3,-2.1) -- (3,-2.9) -- (2.9,-3) -- (2.1,-3.) --
          (2,-3.1) -- (2,-3.9) -- (1.9,-4) -- (1.1,-4) -- (1,-3.9) --
          (1,-3.1) -- (1.1,-3) -- (1.9,-3) -- (2,-2.9) -- (2,-2.1) --
          (1.9,-2.) -- (1.1,-2.) -- (1,-1.9) -- (1,-1.1) -- (0.90,-1) --
          (0.10,-1) -- (0,-0.9) -- (0,-0.1) -- (0.1,0)};
    \filldraw[fill=black!8,xshift=4cm,yshift=8cm] {(0.1,0) coordinate (P6) --
          (0.90,0) -- (1,0.10) -- (1,0.9) -- (1.1,1) --
          (1.9,1) -- (2,1.1) -- (2,1.9) -- (2.1,2) -- (2.9,2) --
          (3,1.9) -- (3,1.1) -- (2.9,1) -- (2.1,1) -- (2,0.9) --
          (2,0.1) -- (2.1,0) -- (2.9,0) -- (3,-0.1) -- (3,-0.9) --
          (3.1,-1) -- (3.9,-1) -- (4,-1.1) -- (4,-1.9) -- (3.9,-2.) --
          (3.1,-2) -- (3,-2.1) -- (3,-2.9) -- (2.9,-3) -- (2.1,-3.) --
          (2,-3.1) -- (2,-3.9) -- (1.9,-4) -- (1.1,-4) -- (1,-3.9) --
          (1,-3.1) -- (1.1,-3) -- (1.9,-3) -- (2,-2.9) -- (2,-2.1) --
          (1.9,-2.) -- (1.1,-2.) -- (1,-1.9) -- (1,-1.1) -- (0.90,-1) --
          (0.10,-1) -- (0,-0.9) -- (0,-0.1) -- (0.1,0)};
    \filldraw[fill=black!10,xshift=8cm,yshift=6cm] {(0.1,0) coordinate (P7) --
          (0.90,0) -- (1,0.10) -- (1,0.9) -- (1.1,1) --
          (1.9,1) -- (2,1.1) -- (2,1.9) -- (2.1,2) -- (2.9,2) --
          (3,1.9) -- (3,1.1) -- (2.9,1) -- (2.1,1) -- (2,0.9) --
          (2,0.1) -- (2.1,0) -- (2.9,0) -- (3,-0.1) -- (3,-0.9) --
          (3.1,-1) -- (3.9,-1) -- (4,-1.1) -- (4,-1.9) -- (3.9,-2.) --
          (3.1,-2) -- (3,-2.1) -- (3,-2.9) -- (2.9,-3) -- (2.1,-3.) --
          (2,-3.1) -- (2,-3.9) -- (1.9,-4) -- (1.1,-4) -- (1,-3.9) --
          (1,-3.1) -- (1.1,-3) -- (1.9,-3) -- (2,-2.9) -- (2,-2.1) --
          (1.9,-2.) -- (1.1,-2.) -- (1,-1.9) -- (1,-1.1) -- (0.90,-1) --
          (0.10,-1) -- (0,-0.9) -- (0,-0.1) -- (0.1,0)};
    \filldraw[fill=black!10,xshift=0cm,yshift=10cm] {(0.1,0) coordinate (P8) --
          (0.90,0) -- (1,0.10) -- (1,0.9) -- (1.1,1) --
          (1.9,1) -- (2,1.1) -- (2,1.9) -- (2.1,2) -- (2.9,2) --
          (3,1.9) -- (3,1.1) -- (2.9,1) -- (2.1,1) -- (2,0.9) --
          (2,0.1) -- (2.1,0) -- (2.9,0) -- (3,-0.1) -- (3,-0.9) --
          (3.1,-1) -- (3.9,-1) -- (4,-1.1) -- (4,-1.9) -- (3.9,-2.) --
          (3.1,-2) -- (3,-2.1) -- (3,-2.9) -- (2.9,-3) -- (2.1,-3.) --
          (2,-3.1) -- (2,-3.9) -- (1.9,-4) -- (1.1,-4) -- (1,-3.9) --
          (1,-3.1) -- (1.1,-3) -- (1.9,-3) -- (2,-2.9) -- (2,-2.1) --
          (1.9,-2.) -- (1.1,-2.) -- (1,-1.9) -- (1,-1.1) -- (0.90,-1) --
          (0.10,-1) -- (0,-0.9) -- (0,-0.1) -- (0.1,0)};
    % Arrows in direction of ABABAD
    %Horizontal
    \draw[] (4.5,-2) node[] {\tiny$>$}; \draw[] (.5,0) node[] {\tiny$>$};
    \draw[] (-3.5,2) node[] {\tiny$>$}; \draw[] (6.5,2) node[] {\tiny$>$};
    \draw[] (2.5,4) node[] {\tiny$>$}; \draw[] (-1.5,6) node[] {\tiny$>$};\draw[] (8.5,6) node[] {\tiny$>$};
    \draw[] (4.5,8) node[] {\tiny$>$}; \draw[] (.5,10) node[] {\tiny$>$};
    % Vertical
    \draw[] (-3,-1.5) node[] {\tiny$\wedge$}; \draw[] (-1,2.5) node[] {\tiny$\wedge$};
    \draw[] (1,-3.5) node[] {\tiny$\wedge$}; \draw[] (1,6.5) node[] {\tiny$\wedge$};
    \draw[] (3,.5) node[] {\tiny$\wedge$};     \draw[] (5,-5.5) node[] {\tiny$\wedge$};\draw[] (5,4.5) node[] {\tiny$\wedge$};
    \draw[] (7,-1.5) node[] {\tiny$\wedge$}; \draw[] (9,2.5) node[] {\tiny$\wedge$};
    \node[fill=red,circle, inner sep=.8pt] at (0,10){};
    \node[fill=red,circle, inner sep=.8pt] at (3,11){};
    \node[fill=red,circle, inner sep=.8pt] at (-2,6){};
    \node[fill=red,circle, inner sep=.8pt] at (1,7){};
    \node[fill=red,circle, inner sep=.8pt] at (4,8){};
    \node[fill=red,circle, inner sep=.8pt] at (7,9){};
    \node[fill=red,circle, inner sep=.8pt] at (-4,2){};
    \node[fill=red,circle, inner sep=.8pt] at (-1,3){};
    \node[fill=red,circle, inner sep=.8pt] at (2,4){};
    \node[fill=red,circle, inner sep=.8pt] at (5,5){};
    \node[fill=red,circle, inner sep=.8pt] at (8,6){};
    \node[fill=red,circle, inner sep=.8pt] at (11,7){};
    \node[fill=red,circle, inner sep=.8pt] at (-3,-1){};
    \node[fill=red,circle, inner sep=.8pt] at (0,0){};
    \node[fill=red,circle, inner sep=.8pt] at (3,1){};
    \node[fill=red,circle, inner sep=.8pt] at (6,2){};
    \node[fill=red,circle, inner sep=.8pt] at (9,3){};
    \node[fill=red,circle, inner sep=.8pt] at (12,4){};
    \node[fill=red,circle, inner sep=.8pt] at (1,-3){};
    \node[fill=red,circle, inner sep=.8pt] at (4,-2){};
    \node[fill=red,circle, inner sep=.8pt] at (7,-1){};
    \node[fill=red,circle, inner sep=.8pt] at (10,0){};
    \node[fill=red,circle, inner sep=.8pt] at (5,-5){};
    \node[fill=red,circle, inner sep=.8pt] at (8,-4){};
    %\node[fill=blue,circle, inner sep=.8pt] at (0,0){};
\end{tikzpicture}
\caption{\small $\tha$-loops and anti-$\tha$-loops when $\tha(a)=ababad$.}\label{fig:D2loops}
\end{figure}

The regions bounded by the curves $K[abcd]$ and $K[dcba]$ and their  translations over even Gaussian integers correspond to the white, respectively, gray squares in Figure~\ref{fig:D2}. It follows that the curves $\rK[\tha(abcd)]$ and $\rK[\tha(dcba)]$ are  $\tha$-loops (white interior) and anti-$\tha$-loops (gray interior) as in the example in Figure~\ref{fig:D2loops} (we see in fact, in this figure part of the $K$-$\tha$-image $K[\tha](\DD)$ of $\DD$). Let $w$ be a $\DD$-word such that $\rK[w]$ is self-avoiding. Since the $\tha$-loops and anti-$\tha$-loops are self-avoiding, the same will be true for the curve $\rK[\tha(w)]$, which will be a subset of the full $\tha$-loops and anti-$\tha$-loops pattern
(after rounding off in the appropriate way). It thus follows by iteration of this observation that all $\rK[\tha^n(a)]$ will be self-avoiding, and so by definition $S$ is self-avoiding. \hfill\qed

\section{Plane-filling folding sequences}\label{sec:planefill}

A surprising property of the curve generated by 2-folding is that it is not only self-avoiding, but very compactly traverses large parts of the set of Gaussian integers (\cite{Gardner,DK}).

We call a folding string $S$ \emph{planefilling} if  the curve associated to the infinite sequence $S^{*\infty}$ passes twice through every grid point lying in some  arbitrary large circle.

It will be convenient to reformulate this to the equivalent property that this curve passes through
 all the segments connecting Gaussian integers forming a $(2k+1)$-cross in $\DD$ for every $k=1,2\ldots$. Here the $(2k+1)$-crosses are translations of $(2k+1)$-crosses \emph{centered at} $(1+i)/2$ defined as the set of edges in $\DD$ lying within the square through the points $k+1+q$, $(k+1)i+q$, $-k-1+q$ and $-(k+1)i+q$, where $q=(1+i)/2$. Note that as directed graphs there are actually two types of $(2k+1)$-crosses: those where the center loop runs clockwise, and those where the  center loop runs anti-clockwise. In the sequel we will not make this distinction explicitly.

\begin{figure}[!t]
\centering
\vspace*{-.5cm}
\begin{tikzpicture}[scale=.5]
% draw 3-cross without arrows
        \draw (-1,0)--(2,0)--(2,1)--(-1,1)--cycle;
        \draw (0,-1)--(1,-1)--(1,2)--(0,2)--cycle;
% draw lines around 3-cross
\draw[black!40, dashed, thin] (-1.6,.5)--(.5,2.6)--(2.6,.5)--(.5,-1.6)--cycle;
% Drawing the red points
    \foreach \x in {-1,0,1,2}{
        \foreach \y in {-1,0,1,2}{
        \node[fill=red,circle, inner sep=.8pt] at (\x,\y)  {};
        }    }
        % Removing the red points
    \foreach \x in {-1,2}{
        \foreach \y in {-1,2}{
        \node[fill=white,circle, inner sep=.9pt] at (\x,\y)  {};
        }    }
% draw arrows
        \draw[] (-.5,0) node[] {\tiny$<$}; \draw[] (0.5,1) node[] {\tiny$<$};
        \draw[] (0.5,-1) node[] {\tiny$<$}; \draw[] (1.5,0) node[] {\tiny$<$};
        \draw[] (-.5,1) node[] {\tiny$>$};
        \draw[] (0.5,2) node[] {\tiny$>$};\draw[] (0.5,0) node[] {\tiny$>$};
        \draw[] (1.5,1) node[] {\tiny$>$};
        \draw[] (0,.5) node[] {\tiny$\vee$};
        \draw[] (1,1.5) node[] {\tiny$\vee$};\draw[] (1,-.5) node[] {\tiny$\vee$};
        \draw[] (2,.5) node[] {\tiny$\vee$};
        \draw[] (-1,.5) node[] {\tiny$\wedge$};
        \draw[] (0,1.5) node[] {\tiny$\wedge$};\draw[] (0,-.5) node[] {\tiny$\wedge$};
        \draw[] (1,.5) node[] {\tiny$\wedge$};
% draw 5-cross without arrows
        \draw (4,0)--(9,0)--(9,1)--(4,1)--cycle;
        \draw (6,-2)--(6,3)--(7,3)--(7,-2)--cycle;
        \draw (5,-1)--(5,2)--(8,2)--(8,-1)--cycle;
% draw lines around 5-cross
\draw[black!40, dashed, thin] (3.4,.5)--(6.5,3.6)--(9.6,.5)--(6.5,-2.6)--cycle;
% Drawing the red points
    \foreach \x in {5,6,7,8}{
        \foreach \y in {-1,0,1,2}{
        \node[fill=red,circle, inner sep=.8pt] at (\x,\y)  {};
        }    }
        % More the red points
    \foreach \x in {4,9}{
        \foreach \y in {0,1}{
        \node[fill=red,circle, inner sep=.8pt] at (\x,\y)  {};
        }    }
          \foreach \x in {6,7}{
        \foreach \y in {-2,3}{
        \node[fill=red,circle, inner sep=.8pt] at (\x,\y)  {};
        }    }
% draw arrows
        \draw[] (6.5,3) node[] {\tiny$<$};
        \draw[] (5.5,2) node[] {\tiny$<$};\draw[] (7.5,2) node[] {\tiny$<$};
        \draw[] (4.5,1) node[] {\tiny$<$};\draw[] (6.5,1) node[] {\tiny$<$};\draw[] (8.5,1) node[] {\tiny$<$};
        \draw[] (5.5,0) node[] {\tiny$<$};\draw[] (7.5,0) node[] {\tiny$<$};
        \draw[] (6.5,-1) node[] {\tiny$<$};
        \draw[] (6.5,2) node[] {\tiny$>$};
        \draw[] (5.5,1) node[] {\tiny$>$};\draw[] (7.5,1) node[] {\tiny$>$};
        \draw[] (4.5,0) node[] {\tiny$>$};\draw[] (6.5,0) node[] {\tiny$>$};\draw[] (8.5,0) node[] {\tiny$>$};
        \draw[] (5.5,-1) node[] {\tiny$>$};\draw[] (7.5,-1) node[] {\tiny$>$};
        \draw[] (6.5,-2) node[] {\tiny$>$};
        \draw[] (4,0.5) node[] {\tiny$\vee$};
        \draw[] (5,1.5) node[] {\tiny$\vee$}; \draw[] (5,-.5) node[] {\tiny$\vee$};
        \draw[] (6,2.5) node[] {\tiny$\vee$}; \draw[] (6,.5) node[] {\tiny$\vee$}; \draw[] (6,-1.5) node[] {\tiny$\vee$};
        \draw[] (7,1.5) node[] {\tiny$\vee$}; \draw[] (7,-.5) node[] {\tiny$\vee$};
        \draw[] (8,.5) node[] {\tiny$\vee$};
        \draw[] (5,.5) node[] {\tiny$\wedge$};
        \draw[] (6,1.5) node[] {\tiny$\wedge$}; \draw[] (6,-.5) node[] {\tiny$\wedge$};
        \draw[] (7,2.5) node[] {\tiny$\wedge$}; \draw[] (7,.5) node[] {\tiny$\wedge$}; \draw[] (7,-1.5) node[] {\tiny$\wedge$};
        \draw[] (8,1.5) node[] {\tiny$\wedge$}; \draw[] (8,-.5) node[] {\tiny$\wedge$};
        \draw[] (9,.5) node[] {\tiny$\wedge$};
\end{tikzpicture}
\caption{\small The 3-cross and 5-cross graphs.}\label{fig:crosses}
\end{figure}
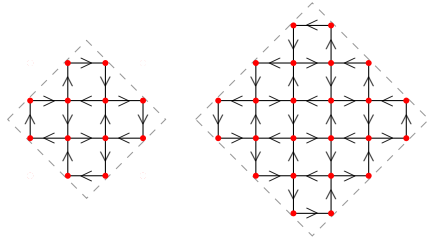

The goal of this section is to give a characterization of self-avoiding planefilling folding strings. In the next section we will see that there exist many self-avoiding planefilling folding strings.

\medskip

We call a string $S$  \emph{maximally simple}, if its $\theta$-loop  is simple and if the region the $\theta$-loop  bounds (by Jordan's Theorem), only contains grid points that are passed twice by the $\theta$-loop.

\begin{lemma}\label{lem:conv-max}
Let $S$ and $T$ be two maximally simple folding strings. Then $S*T$ is also maximally simple.
\end{lemma}

\noindent {\bf Proof:} This is similar to the proof of Lemma~\ref{lem:connect}.

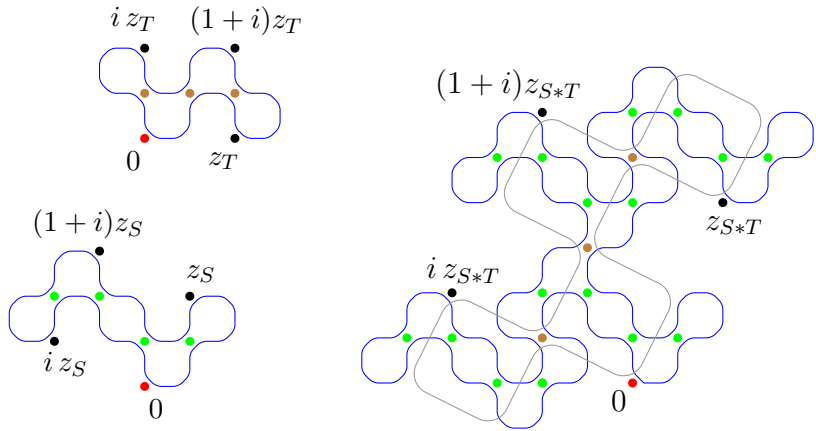
\begin{figure}[!h]
\centering
\begin{tikzpicture}[scale=.6,rounded corners=2pt]
% Fold T:
\draw[blue] (-2., 1.)--(-2.25, 0.75)--(-2.25,0.25)--(-2.5, 0.)-- (-3., 0.)--(-3.25, 0.25)--(-3.25,0.75)--
(-3.5, 1.)-- (-4., 1.)--(-4.25, 1.25)--(-4.25,1.75)--(-4.5, 2.)-- (-5., 2.)--(-5.25, 1.75)--(-5.25,1.25)--
(-5.5, 1.)-- (-6., 1.)--(-6.25, 0.75)--(-6.25,0.25)--(-6., 0.)-- (-5.5, 0.)--(-5.25, 0.25)--(-5.25,0.75)--
(-5., 1.)-- (-4.5, 1.)--(-4.25, 0.75)--(-4.25,0.25)--(-4., 0.)--(-3.5, 0.)--(-3.25,-0.25)--(-3.25,-0.75)--
(-3.,-1.)--(-2.5,-1.)--(-2.25,-0.75)--(-2.25,-0.25)--(-2., 0.)-- (-1.5, 0.)--(-1.25, 0.25)--(-1.25,0.75)--
(-1.5, 1.)--cycle;
% Fold S:
\draw[blue] (-3,4.5)--(-2.5,4.5)--(-2.25,4.75)--(-2.25,5.25)--(-2,5.5)--(-1.5,5.5)--(-1.25,5.25)--(-1.25,4.75)--
(-1,4.5)--(-.5,4.5)--(-.25,4.75)--(-.25,5.25)--(-.5,5.5)--(-1,5.5)--(-1.25,5.75)--(-1.25,6.25)--(-1.5,6.5)--
(-2,6.5)--(-2.25,6.25)--(-2.25,5.75)--(-2.5,5.5)--(-3,5.5)--(-3.25,5.75)--(-3.25,6.25)--(-3.5,6.5)--
(-4,6.5)--(-4.25,6.25)--(-4.25,5.75)--(-4,5.5)--(-3.5,5.5)--(-3.25,5.25)--(-3.25,4.75)-- cycle;
% \node  at (0,0) -> (-.25,0)
\node at (-3.,-1.5) {\small $0$};
\node at (-2,1.5) {\small $z_S$};
\node at (-4.5,2.6) {\small $(1+i)z_S$};
\draw (-5,-0.5) node[] {\small $i\,z_S$};
\node at (-3.5,4) {\small $0$};
\node at (-1.5, 4) {\small $z_T$};
\node at (-1,7.1) {\small $(1+i)z_T$};
\draw (-3.5,7.15) node[] {\small $i\,z_T$};
% Red dots
\node[fill=red,circle, inner sep=1.2pt] at (-3.25,-1) {};
\node[fill=black,circle, inner sep=1.2pt] at (-5.25,0) {};
\node[fill=black,circle, inner sep=1.2pt] at (-2.25,1) {};
\node[fill=black,circle, inner sep=1.2pt] at (-4.25,2) {};
\node[fill=red,circle, inner sep=1.2pt] at (-3.25,4.5) {};
\node[fill=black,circle, inner sep=1.2pt] at (-1.25,4.5) {};
\node[fill=black,circle, inner sep=1.2pt] at (-3.25,6.5) {};
\node[fill=black,circle, inner sep=1.2pt] at (-1.25,6.5) {};
% Brown dots
\node[fill=brown,circle, inner sep=1.2pt] at (-3.25,5.5) {};
\node[fill=brown,circle, inner sep=1.2pt] at (-2.25,5.5) {};
\node[fill=brown,circle, inner sep=1.2pt] at (-1.25,5.5) {};
% Green dots
\node[fill=green,circle, inner sep=1.2pt] at (-3.25,0) {};
\node[fill=green,circle, inner sep=1.2pt] at (-2.25,0) {};
\node[fill=green,circle, inner sep=1.2pt] at (-5.25,1) {};
\node[fill=green,circle, inner sep=1.2pt] at (-4.25,1) {};
\end{tikzpicture}\hspace*{.5cm}
\begin{tikzpicture}[scale=.6,rounded corners=2pt]
\draw[blue]
(0., 0.)--(0.5, 0.)--(0.75,0.25)--(0.75, 0.75)--(1.,1.)--(1.5,1.)--(1.75, 1.25)--(1.75, 1.75)--(1.5, 2.)--
(1., 2.)--(0.75,1.75)--(0.75,1.25)--(0.5,1.)--(0.,1.)--(-0.25,1.25)--(-0.25, 1.75)--(-0.5, 2.)--(-1.,2.)--
(-1.25, 2.25)--(-1.25, 2.75)--(-1., 3.)--(-0.5, 3.)--(-0.25, 3.25)--(-0.25, 3.75)--(0., 4.)--(0.5,4.)--
(0.75,4.25)--(0.75,4.75)--(0.5,5.)--(0.,5.)--(-0.25, 5.25)--(-0.25,5.75)--(0.,6.)--(0.5,6.)--(0.75,5.75)--
(0.75, 5.25)--(1.,5.)--(1.5,5.)--(1.75,4.75)--(1.75,4.25)--(2.,4.)--(2.5, 4.)--(2.75, 4.25)--(2.75,4.75)--
(3., 5.)--(3.5, 5.)--(3.75,5.25)--(3.75,5.75)--(3.5, 6.)--(3., 6.)--(2.75, 5.75)--(2.75,5.25)--(2.5, 5.)--
(2., 5.)--(1.75, 5.25)--(1.75, 5.75)--(1.5, 6.)--(1.,6.)--(0.75,6.25)--(0.75,6.75)--(0.5,7.)--(0.,7.)--
(-0.25, 6.75)--(-0.25, 6.25)--(-0.5, 6.)--     (-1., 6.)-- (-1.25, 5.75)--(-1.25,5.25)--
(-1., 5.)--(-0.5, 5.)-- (-0.25, 4.75)--(-0.25,4.25)--(-0.5, 4.)-- (-1., 4.)--(-1.25, 4.25)--(-1.25,4.75)--
(-1.5, 5.)-- (-2., 5.)--(-2.25, 5.25)--(-2.25,5.75)--(-2.5, 6.)-- (-3., 6.)--(-3.25, 5.75)--(-3.25,5.25)--
(-3.5, 5.)-- (-4., 5.)--(-4.25, 4.75)--(-4.25,4.25)--(-4., 4.)-- (-3.5, 4.)--(-3.25, 4.25)--(-3.25,4.75)--
(-3., 5.)-- (-2.5, 5.)--(-2.25, 4.75)--(-2.25,4.25)--(-2., 4.)-- (-1.5, 4.)--(-1.25, 3.75)--(-1.25,3.25)--
(-1.5, 3.)-- (-2., 3.)--(-2.25, 2.75)--(-2.25,2.25)--(-2.5, 2.)-- (-3., 2.)--(-3.25, 1.75)--(-3.25,1.25)--
(-3., 1.)-- (-2.5, 1.)--(-2.25, 0.75)--(-2.25,0.25)--(-2.5, 0.)-- (-3., 0.)--(-3.25, 0.25)--(-3.25,0.75)--
(-3.5, 1.)-- (-4., 1.)--(-4.25, 1.25)--(-4.25,1.75)--(-4.5, 2.)-- (-5., 2.)--(-5.25, 1.75)--(-5.25,1.25)--
(-5.5, 1.)-- (-6., 1.)--(-6.25, 0.75)--(-6.25,0.25)--(-6., 0.)-- (-5.5, 0.)--(-5.25, 0.25)--(-5.25,0.75)--
(-5., 1.)-- (-4.5, 1.)--(-4.25, 0.75)--(-4.25,0.25)--(-4., 0.)--(-3.5, 0.)--(-3.25,-0.25)--(-3.25,-0.75)--
(-3.,-1.)--(-2.5,-1.)--(-2.25,-0.75)--(-2.25,-0.25)--(-2., 0.)-- (-1.5, 0.)--(-1.25, 0.25)--(-1.25,0.75)--
(-1.5, 1.)--(-2., 1.)-- (-2.25, 1.25)--(-2.25,1.75)--(-2., 2.)--(-1.5, 2.)-- (-1.25,1.75)--(-1.25, 1.25)--
(-1.,1.)--(-0.5,1.)--(-0.25,0.75)--(-.25,.25)--(0,0);
\draw[black!40,rounded corners=8pt]
(-0.25,0)--++(1,2)--++(-2,1)--++(1,2)--++(2,-1)--++(1,2)--++(-2,1)--++(-1,-2)--++(-2,1)--++(-1,-2)--++(2,-1)
--++(-1,-2)--++(-2,1)--++(-1,-2)--++(2,-1)--++(1,2)--cycle;
% \node  at (0,0) -> (-.25,0)
\node at (-.55,-.4) {$0$}; \node at (2,3.5) {$z_{S*T}$};
\node at (-3,6.6) {$(1+i)z_{S*T}$}; \node at (-4,2.5) {$i\,z_{S*T}$};
% Red dots
\node[fill=red,circle, inner sep=1.2pt] at (-.25,0) {};
\node[fill=black,circle, inner sep=1.2pt] at (1.75,4) {};
\node[fill=black,circle, inner sep=1.2pt] at (-2.25,6) {};
\node[fill=black,circle, inner sep=1.2pt] at (-4.25,2) {};
% Blue dots
\node[fill=brown,circle, inner sep=1.2pt] at (-2.25,1) {};
\node[fill=brown,circle, inner sep=1.2pt] at (-1.25,3) {};
\node[fill=brown,circle, inner sep=1.2pt] at (-.25,5) {};
% Green dots
\node[fill=green,circle, inner sep=1.2pt] at (-3.25,0) {};
\node[fill=green,circle, inner sep=1.2pt] at (-2.25,0) {};
\node[fill=green,circle, inner sep=1.2pt] at (-2.25,0) {};
\node[fill=green,circle, inner sep=1.2pt] at (-5.25,1) {};
\node[fill=green,circle, inner sep=1.2pt] at (-4.25,1) {};
\node[fill=green,circle, inner sep=1.2pt] at (-.25,1) {};
\node[fill=green,circle, inner sep=1.2pt] at (.75,1) {};
\node[fill=green,circle, inner sep=1.2pt] at (-2.25,2) {};
\node[fill=green,circle, inner sep=1.2pt] at (-1.25,2) {};
\node[fill=green,circle, inner sep=1.2pt] at (-1.25,4) {};
\node[fill=green,circle, inner sep=1.2pt] at (-0.25,4) {};
\node[fill=green,circle, inner sep=1.2pt] at (-3.25,5) {};
\node[fill=green,circle, inner sep=1.2pt] at (-2.25,5) {};
\node[fill=green,circle, inner sep=1.2pt] at (1.75,5) {};
\node[fill=green,circle, inner sep=1.2pt] at (2.75,5) {};
\node[fill=green,circle, inner sep=1.2pt] at (-.25,6) {};
\node[fill=green,circle, inner sep=1.2pt] at (.75,6) {};
\end{tikzpicture}
\caption{\small The $\tha$-loops of the folds $S=DUDD$, $T=DUU$, and $S*T$.}\label{fig:conv-max-simp}
\end{figure}

We obtain the $\tha_{S*T}$-loop by replacing scaled unit-loops in the $\tha_T$-loop by $\tha_S$-loops. Grid points in the interior of the $\tha_{S*T}$-loop will either lie in these $\tha_S$-loops, or will be
in the interior of the scaled $\tha_T$-loop, cf. Figure~\ref{fig:conv-max-simp}. In both cases these points are visited twice, and so the $\tha_{S*T}$-loop is maximally simple.\qed

\medskip

Lemma \ref{lem:conv-max} forms the basis for the planefilling property. It implies that for a given string $S$ with maximally simple $\tha$-loop all its $\tha^n$-loops will also be maximally simple. However, there are some subtleties in actually proving that this implies that $S$ is planefilling.

We call a word $w$ from $A^*$ a \emph{loop word} if it is  a $\DD$-word, and $f(w)=0$, where $f$ is the homomorphism from Equation~(\ref{eq:f}).

We call a loop word \emph{maximally simple} if $\rK[w]$ is a self-avoiding curve that visits each
grid point in its interior twice. Examples are the 8 words
\begin{equation}\label{eq:squareword}
\sigma^k(abcd),\; \mathrm{and}\; \sigma^k(adcb)\quad \mathrm{for}\; k=0,1,2,3,
\end{equation}
which we call \emph{square} words.

\begin{lemma}\label{lem:loopword}
Let $w$ be a maximally simple loop word. Then $w$ contains at least one of the square words.
\end{lemma}

\begin{figure}[h]
\centering
\begin{tikzpicture}[scale=.6,rounded corners=2pt]
\draw[blue] (0,0)--(0.50, 0.)--(0.75, 0.25) --(0.75, 0.75) --
(0.5,1)--(0,1)--(-0.25,1.25)--(-0.25,1.75)--(0,2)--(0.5,2)--(0.75, 2.25) --(0.75, 2.75) --(1,3)--
(1.50, 3)--(1.75, 2.75)--(1.75, 2.25)--(1.5, 2)--(1,2)--(0.75,1.75)--(0.75,1.25)--(1., 1.) --(1.5, 1.) --(1.75, 1.25)--
(1.75, 1.75)--(2., 2.)--(2.5, 2.)--(2.75, 1.75)--(2.75, 1.25) --(3., 1.) --(3.5, 1.)--(3.75, 0.75) --(3.75, 0.25)--
(4, 0.)--(4.50, 0.)--(4.75, 0.25)--(4.75, 0.75)--(4.5,1)--(4,1)--(3.75, 1.25) --(3.75, 1.75)--(4,2)--(4.5,2)--(4.75, 2.25)
--(4.75, 2.75)--(4.5,3)-- (4,3)--(3.75, 3.25)--(3.75, 3.75)-- (3.5, 4.)--(3,4)--(2.75, 3.75)--(2.75, 3.25)--(3,3)--(3.5,3)
--(3.75, 2.75)--(3.75, 2.25)--(3.5, 2.) --(3., 2.) -- (2.75, 2.25) --(2.75, 2.75) --(2.5, 3.) --(2., 3.) --(1.75, 3.25)
--(1.75, 3.75) --(1.5, 4.)--(1., 4.)--(0.75, 3.75) --(0.75, 3.25) --(0.50, 3.) --(0., 3.) --(-0.25, 2.75) --(-0.25, 2.25)
 --(-0.50, 2.)--(-1., 2.0) --(-1.25, 1.75)--(-1.25, 1.25)--(-1., 1.)--(-0.50, 1.)--(-0.25, 0.75)--(-0.25, 0.25)--(0., 0.);
% De mid grid points:
\node[fill=green,circle, inner sep=1.2pt] at (0.25,0.5) {};
\node[fill=green,circle, inner sep=1.2pt] at (-0.75,1.5){};
\node[fill=green,circle, inner sep=1.2pt] at (1.25,1.5) {};
\node[fill=green,circle, inner sep=1.2pt] at (3.25,1.5) {};
\node[fill=green,circle, inner sep=1.2pt] at (0.25,2.5) {};
\node[fill=green,circle, inner sep=1.2pt] at (2.25,2.5) {};
\node[fill=green,circle, inner sep=1.2pt] at (4.25,2.5) {};
\node[fill=green,circle, inner sep=1.2pt] at (1.25,3.5) {};
\node[fill=green,circle, inner sep=1.2pt] at (3.25,3.5) {};
\node[fill=green,circle, inner sep=1.2pt] at (4.25,0.5) {};
% De oorsprong:
\node[fill=red,circle, inner sep=.8pt] at (-0.25,0)  {};
% De mid grid tree:
\draw[green] (0.25,0.5)--(-0.75,1.5)--(0.25,2.5)  --(1.25,3.5) --(2.25,2.5)--(1.25,1.5);
\draw[green] (2.25,2.5)--(3.25,1.5)--(4.25,2.5)  --(3.25,3.5);
\draw[green] (3.25,1.5)--(4.25,0.5);
\end{tikzpicture}
\caption{\small The curve $\rK[w]$ of a maximally simple  loop word of length $m=40$.}\label{fig:maxsimp}
\end{figure}
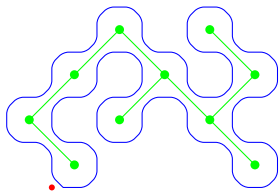

\noindent {\bf Proof:} A point $u+iv$ with $u=N+1/2$ and $v=M+1/2$, where $N$ and $M$ are integers will be  called a \emph{mid grid point}.
Maximally simple loop words $w$ of length $m$ are associated to particularly nice curves $\rK[w]$.
The curve $\rK[w]$ contains exactly $m/4$  mid grid points in its interior. Consider the graph $\cal{G}$ whose vertices are these mid grid points, and where there is an edge
if two vertices are neighbors (i.e., the corresponding mid grid points are at
distance $\sqrt{2}$), see also Figure~\ref{fig:maxsimp}.
If $m=4$, then $w$ \emph{is} a square word, so consider $m\ge 8$.
Since $\rK[w]$ is simple, $\cal{G}$ is connected, and since $\rK[w]$ is maximally simple, $\cal{G}$ has no circuits, and so $\cal{G}$ is a tree.  A mid grid point corresponding to any leaf of this tree is the center of a $\rK[\sigma^k(abcd)]$ or a $\rK[\sigma^k(adcb)]$ (ignoring rounding of). Since $m\ge 8$ the tree has at least two leaves. At least one of the corresponding two square words will occur in $w$ (one has to take into account the possibility that one of the two corresponding mid grid points is the point $(1+i)/2$).\qed

\medskip

\begin{lemma}\label{lem:squre-in-tha}
Let a folding morphism $\tha$  have a maximally simple $\tha$-loop. Then   a square word will occur in $\tha^4(a)$.
\end{lemma}

\noindent {\bf Proof:} If $m=2$, then $\tha(a)=ab$, and it is easily checked that $w_{\textsc{sq}}=bcda$ occurs in $\tha^4(a)$, so suppose $m\ge 3$. From Lemma~\ref{lem:loopword} we know that $w=\tha(abcd)$ contains a square word $w_{\textsc{sq}}$. As $m\ge 3$, $w_{\textsc{sq}}$ has to occur in at least one of $\tha(ab)$ or $\tha(bc)$, or $\tha(cd)$.

Obviously $ab$ occurs in $\tha^3(a)=ab\ldots$. We claim that also $bc$ occurs in $\tha^3(a)$. To see this, note that if $bc$ occurs in $\tha(a)$, then it occurs in
$\tha^2(a)$ and in $\tha^3(a)$, since $\tha(b)$ ends in $b$ and $\tha(c)$ starts with $c$. Now either
$\tha(a)=abc\ldots$, and then $bc$ will also appear in $\tha^3(a)$, or $\tha(a)=aba\ldots$, and then $bc$ will  appear in $\tha^2(a)$ (and as argued before, in $\tha^3(a)$), since in that case
$\tha(b)=\sigma\tau\tha(a)$ ends in $bcb$. Finally, $cd$ occurs in $\tha^3(a)$, since $c$ occurs in $\tha(b)$, and $\tha(c)=(\sigma\tau)^2\tha(a)=cd\ldots$.

Since  $ab$, $bc$ and $cd$  all
occur in $\tha^3(a)$, $w_{\textsc{sq}}$ will occur in $\tha^4(a)$.\qed

%\medskip

\begin{theorem}\label{th:planefill}
Let $S$ be a self-avoiding folding string. Then $S$ is a planefilling string if and only if the $\theta$-loop of $S$ is maximally simple.
\end{theorem}

\noindent {\bf Proof:} \emph{Part 1:} If the $\tha$-loop is not maximally simple, then at least one of the grid points in its interior is visited only once. It follows that a $(2k+1)$-cross can nowhere occur if $k$ is large enough.

\medskip

\emph{Part 2:} Note that the 8 square words have two different types, those associated to a $\tha$-loop: $abcd, bcda, cdab,$ and $dabc$, and the other four, which are associated to an anti-$\tha$-loop (cf.~the white, respectively gray squares in
Figure~\ref{fig:D2}). In terms of curves this means that the $\tha$-loop runs counter clockwise, and the anti-$\tha$-loop clockwise.

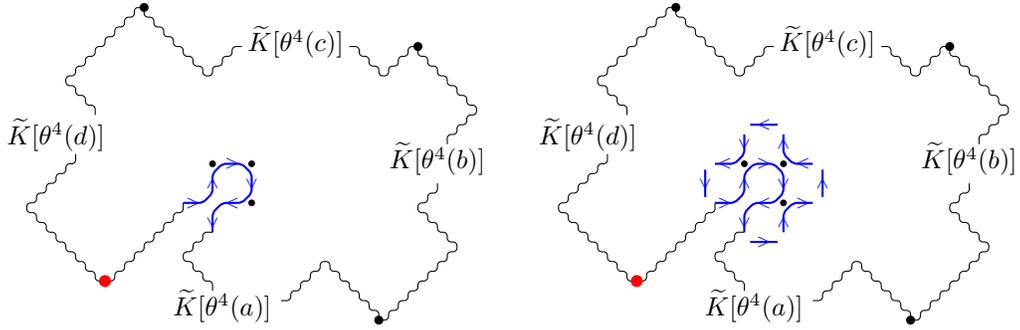
\begin{figure}[t!]
\centering
\begin{tikzpicture}[domain=0:1,scale=0.13,rounded corners=1.5pt]
% Begin van theta(a)--tot de grote loop
\draw (12,8) -- ++(1,0) -- ++(0,1) -- ++(1,0) -- ++(0,1) -- ++(1,0) -- ++(0,1) -- ++(1,0) -- ++(0,1) --
++(1,0) -- ++(0,1) -- ++(1,0) -- ++(0,1) -- ++(1,0) -- ++(0,1) -- ++(1,0) -- ++(0,1) ;
% de grote loop
\draw[blue,thick] (20,16) -- ++(2,0) -- ++(1,1) -- ++(0,2) -- ++(1,1) -- ++(2,0) -- ++(1,-1)
-- ++(0,-2) -- ++(-1,-1) -- ++(-2,0) -- ++(-1,-1) -- ++ (0,-2);
% Eind van theta(a)
\draw (23,13)--++(-1,0)--++(0,-1)--++(-1,0)--++(0,-1)--++(-1,0)-- ++(0,-1)%%hoekje om
--++(1,0)--++(0,-1)--++(1,0)--++(0,-1)--++(1,0)--++(0,-1)% nu zonder tekenen!:
++(1,0)++(0,-1)++(1,0)++(0,-1)++(1,0)++(0,-1)++(1,0)++(0,-1)++(1,0)++(0,1)++(1,0)++(0,1)++(1,0)++(0,1)% nu wel tekenen:
--++(1,0)--++(0,1)--++(1,0)--++(0,1)
--++(1,0)--++(0,1)--++(1,0)--++(0,1)%%hoekje om
-- ++(1,0) -- ++(0,-1) -- ++(1,0) -- ++(0,-1) -- ++(1,0) -- ++(0,-1) -- ++(1,0) -- ++(0,-1)
-- ++(1,0) -- ++(0,-1) -- ++(1,0) -- ++(0,-1);
% Begin theta(b): trappetje van (40,4) tot ?
\draw (40,4)--++(1,0)--++(0,1)--++(1,0)--++(0,1)--++(1,0)--++(0,1)--++(1,0) -- ++(0,1)-- ++(1,0) -- ++(0,1) -- ++(1,0) -- ++(0,1) -- ++(1,0) -- ++(0,1) -- ++(1,0) -- ++(0,1)%%hoekje om
--++(-1,0)--++(0,1)--++(-1,0)-- ++(0,1)-- ++(-1,0)--++(0,1)--++(-1,0)-- ++(0,1)%%hoekje om
-- ++(1,0) -- ++(0,1) -- ++(1,0) -- ++(0,1);
% Het einde van theta(b)
\draw (44,32) -- ++(0,-1)--++(1,0)--++(0,-1)-- ++(1,0)--++(0,-1)--++(1,0)--++(0,-1)--++(1,0)--++(0,-1)--++(1,0) -- ++(0,-1) -- ++(1,0) -- ++(0,-1)%%hoekje om
--++(-1,0)--++(0,-1)--++(-1,0)-- ++(0,-1)--++(-1,0)--++(0,-1);
% Het begin van theta(c)
\draw (44,32)--++(-1,0)--++(0,-1)--++(-1,0)--++(0,-1) -- ++(-1,0) -- ++(0,-1)--++(-1,0)%%hoekje om
--++(0,1)--++(-1,0)-- ++(0,1)--++(-1,0)--++(0,1)--++(-1,0);
% Het einde van theta(c)
\draw (16,36) -- ++(0,-1) -- ++(1,0) -- ++(0,-1) -- ++(1,0) -- ++(0,-1) -- ++(1,0) -- ++(0,-1)-- ++(1,0) -- ++(0,-1) -- ++(1,0) -- ++(0,-1) -- ++(1,0) -- ++(0,-1)%%hoekje om
--++(1,0)--++(0,1)--++(1,0)--++(0,1)--++(1,0)--++(0,1)--++(1,0);
% Het begin van theta(d)
\draw (16,36)--++(-1,0)--++(0,-1)--++(-1,0)--++(0,-1)--++(-1,0)--++(0,-1)--++(-1,0)--++(0,-1)-- ++(-1,0)--++(0,-1)--++(-1,0)--++(0,-1)--++(-1,0)--++(0,-1)--++(-1,0)--++(0,-1)%%hoekje om
--++(1,0)--++(0,-1)-- ++(1,0)--++(0,-1)--++(1,0)--++(0,-1);
% Het einde van theta(d)
\draw (12,8)--++(-1,0)--++(0,1)--++(-1,0)-- ++(0,1)--++(-1,0)--++(0,1)--++(-1,0)-- ++(0,1)--++(-1,0)--++(0,1)--++(-1,0)-- ++(0,1)--++(-1,0)--++(0,1)--++(-1,0)-- ++(0,1)%%hoekje om
--++(1,0)--++(0,1)--++(1,0)--++(0,1)--++(1,0)--++(0,1)--++(1,0)--++(0,1)--++(1,0)--++(0,1);
% Interne punten:
\node[fill=black,circle, inner sep=0.9pt] at (23,20) {};
\node[fill=black,circle, inner sep=0.9pt] at (27,20) {};
\node[fill=black,circle, inner sep=0.9pt] at (27,16) {};
% Pijlen
 \draw[blue] (25,16) node[] {\tiny$<$};
 \draw[blue] (25,20) node[] {\tiny$>$};
 \draw[blue] (27,18) node[] {\tiny$\vee$};
 \draw[blue] (23,18) node[] {\tiny$\wedge$};
 \draw[blue] (21,16) node[] {\tiny$>$};
 \draw[blue] (23,14) node[] {\tiny$\vee$};
 % Hoekpunten
\node[fill=red,circle, inner sep=1.6pt] at (12,8) {};
\node[fill=black,circle, inner sep=1.2pt] at (16,36) {};
\node[fill=black,circle, inner sep=1.2pt] at (40,4) {};
\node[fill=black,circle, inner sep=1.2pt] at (44,32) {};
% Namen:
 \draw (24,5.5) node[] {\footnotesize$\rK[\tha^4(a)]$};
 \draw (46,20.5) node[] {\footnotesize$\rK[\tha^4(b)]$};
 \draw (31.5,32.5) node[] {\footnotesize$\rK[\tha^4(c)]$};
 \draw (7,23) node[] {\footnotesize$\rK[\tha^4(d)]$};
\end{tikzpicture}\hspace*{0.3cm}
%%%%%%%%%%%%%%%%%%%%%%%%%%%%%%%%%%%%%%%%%%%%%%%%%%%%%%%%%%%%%%%%%%%%%%%%%%%%%%%%%%%%%%%%% PLAATJE RIGHT
\begin{tikzpicture}[domain=0:1,scale=0.13,rounded corners=1.5pt]
% Begin van theta(a)--tot de grote loop
\draw (12,8) -- ++(1,0) -- ++(0,1) -- ++(1,0) -- ++(0,1) -- ++(1,0) -- ++(0,1) -- ++(1,0) -- ++(0,1) --
++(1,0) -- ++(0,1) -- ++(1,0) -- ++(0,1) -- ++(1,0) -- ++(0,1) -- ++(1,0) -- ++(0,1) ;
% de grote loop
\draw[blue,thick] (20,16) -- ++(2,0) -- ++(1,1) -- ++(0,2) -- ++(1,1) -- ++(2,0) -- ++(1,-1)
-- ++(0,-2) -- ++(-1,-1) -- ++(-2,0) -- ++(-1,-1) -- ++ (0,-2);
% Eind van theta(a)
\draw (23,13)--++(-1,0)--++(0,-1)--++(-1,0)--++(0,-1)--++(-1,0)-- ++(0,-1)%%hoekje om
--++(1,0)--++(0,-1)--++(1,0)--++(0,-1)--++(1,0)--++(0,-1)% nu zonder tekenen!:
++(1,0)++(0,-1)++(1,0)++(0,-1)++(1,0)++(0,-1)++(1,0)++(0,-1)++(1,0)++(0,1)++(1,0)++(0,1)++(1,0)++(0,1)% nu wel tekenen:
--++(1,0)--++(0,1)--++(1,0)--++(0,1)
--++(1,0)--++(0,1)--++(1,0)--++(0,1)%%hoekje om
-- ++(1,0) -- ++(0,-1) -- ++(1,0) -- ++(0,-1) -- ++(1,0) -- ++(0,-1) -- ++(1,0) -- ++(0,-1)
-- ++(1,0) -- ++(0,-1) -- ++(1,0) -- ++(0,-1);
% Begin theta(b): trappetje van (40,4) tot ?
\draw (40,4)--++(1,0)--++(0,1)--++(1,0)--++(0,1)--++(1,0)--++(0,1)--++(1,0) -- ++(0,1)-- ++(1,0) -- ++(0,1) -- ++(1,0) -- ++(0,1) -- ++(1,0) -- ++(0,1) -- ++(1,0) -- ++(0,1)%%hoekje om
--++(-1,0)--++(0,1)--++(-1,0)-- ++(0,1)-- ++(-1,0)--++(0,1)--++(-1,0)-- ++(0,1)%%hoekje om
-- ++(1,0) -- ++(0,1) -- ++(1,0) -- ++(0,1);
% Het einde van theta(b)
\draw (44,32) -- ++(0,-1)--++(1,0)--++(0,-1)-- ++(1,0)--++(0,-1)--++(1,0)--++(0,-1)--++(1,0)--++(0,-1)--++(1,0) -- ++(0,-1) -- ++(1,0) -- ++(0,-1)%%hoekje om
--++(-1,0)--++(0,-1)--++(-1,0)-- ++(0,-1)--++(-1,0)--++(0,-1);
% Het begin van theta(c)
\draw (44,32)--++(-1,0)--++(0,-1)--++(-1,0)--++(0,-1) -- ++(-1,0) -- ++(0,-1)--++(-1,0)%%hoekje om
--++(0,1)--++(-1,0)-- ++(0,1)--++(-1,0)--++(0,1)--++(-1,0);
% Het einde van theta(c)
\draw (16,36) -- ++(0,-1) -- ++(1,0) -- ++(0,-1) -- ++(1,0) -- ++(0,-1) -- ++(1,0) -- ++(0,-1)-- ++(1,0) -- ++(0,-1) -- ++(1,0) -- ++(0,-1) -- ++(1,0) -- ++(0,-1)%%hoekje om
--++(1,0)--++(0,1)--++(1,0)--++(0,1)--++(1,0)--++(0,1)--++(1,0);
% Het begin van theta(d)
\draw (16,36)--++(-1,0)--++(0,-1)--++(-1,0)--++(0,-1)--++(-1,0)--++(0,-1)--++(-1,0)--++(0,-1)-- ++(-1,0)--++(0,-1)--++(-1,0)--++(0,-1)--++(-1,0)--++(0,-1)--++(-1,0)--++(0,-1)%%hoekje om
--++(1,0)--++(0,-1)-- ++(1,0)--++(0,-1)--++(1,0)--++(0,-1);
% Het einde van theta(d)
\draw (12,8)--++(-1,0)--++(0,1)--++(-1,0)-- ++(0,1)--++(-1,0)--++(0,1)--++(-1,0)-- ++(0,1)--++(-1,0)--++(0,1)--++(-1,0)-- ++(0,1)--++(-1,0)--++(0,1)--++(-1,0)-- ++(0,1)%%hoekje om
--++(1,0)--++(0,1)--++(1,0)--++(0,1)--++(1,0)--++(0,1)--++(1,0)--++(0,1)--++(1,0)--++(0,1);
% Interne punten:
\node[fill=black,circle, inner sep=0.9pt] at (23,20) {};
\node[fill=black,circle, inner sep=0.9pt] at (27,20) {};
\node[fill=black,circle, inner sep=0.9pt] at (27,16) {};
% Pijlen
 \draw[blue] (25,16) node[] {\tiny$<$};
 \draw[blue] (25,20) node[] {\tiny$>$};
 \draw[blue] (27,18) node[] {\tiny$\vee$};
 \draw[blue] (23,18) node[] {\tiny$\wedge$};
 \draw[blue] (21,16) node[] {\tiny$>$};
 \draw[blue] (23,14) node[] {\tiny$\vee$};
 % Hoekpunten
\node[fill=red,circle, inner sep=1.6pt] at (12,8) {};
\node[fill=black,circle, inner sep=1.2pt] at (16,36) {};
\node[fill=black,circle, inner sep=1.2pt] at (40,4) {};
\node[fill=black,circle, inner sep=1.2pt] at (44,32) {};
% Namen:
 \draw (24,5.5) node[] {\footnotesize$\rK[\tha^4(a)]$};
 \draw (46,20.5) node[] {\footnotesize$\rK[\tha^4(b)]$};
 \draw (31.5,32.5) node[] {\footnotesize$\rK[\tha^4(c)]$};
 \draw (7,23) node[] {\footnotesize$\rK[\tha^4(d)]$};
% het 3-CROSS!
\draw[blue,thick] (19,16.6)--++(0,2.8)  (20,20)--++(2,0)--++(1,1)--++(0,2)   (23.6,24)--++(2.8,0)
(27,23)--++(0,-2)--++(1,-1)--++(2,0)   (31,19.4)--++(0,-2.8)   (30,16)--++(-2,0)--++(-1,-1)--++(0,-2)   (23.6,12)--++(2.8,0);
% Pijlen in  het 3-CROSS!
 \draw[blue] (25,24) node[] {\tiny$<$};
 \draw[blue] (25,12) node[] {\tiny$>$};
 \draw[blue] (19,18) node[] {\tiny$\vee$};
 \draw[blue] (31,18) node[] {\tiny$\wedge$};
 \draw[blue] (21,20) node[] {\tiny$<$};
 \draw[blue] (23,22) node[] {\tiny$\vee$};
 \draw[blue] (27,22) node[] {\tiny$\wedge$};
 \draw[blue] (29,20) node[] {\tiny$<$};
 \draw[blue] (27,14) node[] {\tiny$\wedge$};
 \draw[blue] (29,16) node[] {\tiny$>$};
\end{tikzpicture}
\caption{\small Left: $\rK[w_{\textsc{sq}}]$ in $\tha^4$-loop. Right: forced segments in the $\tha^4$-loop.}\label{fig:opptype}
\end{figure}

Let  $w_{\textsc{sq}}$ be the square word occurring in $\tha^4(a)$ according to Lemma~\ref{lem:squre-in-tha}. Suppose that $w_{\textsc{sq}}$ is of the second (clockwise) type (otherwise we replace below everywhere the $\tha^4$-loop by the  anti-$\tha^4$-loop). Now consider the curve $\rK[\tha^4(a)]$ in  the $\tha^4$-loop. Since the type of this loop is \emph{opposite} to the type of $w_{\textsc{sq}}$, the 3 grid points passed by $\rK[w_{\textsc{sq}}]$ must be \emph{interior} points of this loop, cf. Figure~\ref{fig:opptype}, left. This implies that the 10 segments indicated in Figure~\ref{fig:opptype}, right must all be part of the $\tha^4$-loop. In other words: the $\tha^4$-loop contains a 3-cross.
Since $w_{\textsc{sq}}$ occurs in $\tha^4(a)$ this implies that $\rK[\tha^8(a)]$ contains a 3-cross.

\begin{wrapfigure}{r}{15mm}%%%%was 15 mm
\vspace*{-.6cm}
%%%%%\vspace*{-0.4cm}% Deze als je net onder een figuur zit!
\begin{center}
\begin{tikzpicture}[scale=.5]
% draw 3-cross without arrows
        \draw (0,1)--(1,1);  \draw (0,0)--(1,0)--(1,2)--(0,2)--cycle;
% Drawing the red points
 \foreach \x in {0,1}{\foreach \y in {0,1,2}{\node[fill=red,circle, inner sep=.8pt] at (\x,\y) {};}}
% draw arrows
 \draw[] (0.5,1) node[] {\tiny$>$}   (0.5,2) node[] {\tiny$<$} (0.5,0) node[] {\tiny$<$};
 \draw[] (0,.5) node[] {\tiny$\wedge$}   (1,1.5) node[] {\tiny$\wedge$};
 \draw[] (0,1.5) node[] {\tiny$\vee$}    (1,.5) node[] {\tiny$\vee$};
 \draw[] (0,1) node[left] {\tiny$0$}    (1,1) node[right] {\tiny\;$1$};
 \draw[] (0,0) node[left] {\tiny$-i$}   (1,0) node[right] {\tiny$1\!-\!i$};
 \draw[] (0,2) node[left] {\tiny$i$}    (1,2) node[right] {\tiny$1+i$};
\end{tikzpicture}
\end{center}
\vspace*{-.9cm}
\end{wrapfigure}

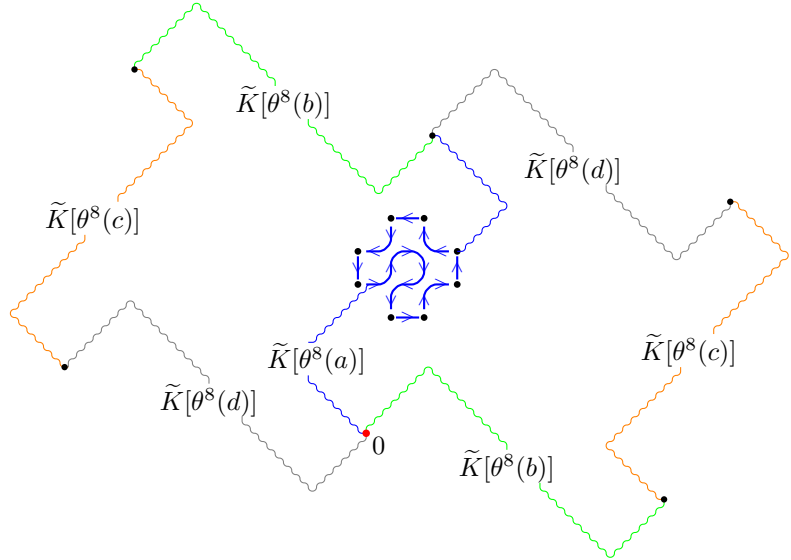
\begin{figure}[h!]
\centering
\begin{tikzpicture}[domain=0:1,scale=0.11,rounded corners=1.5pt]
\path(1,0)  coordinate (a);
\path(0,1)  coordinate (b);
\path(-1,0) coordinate (c);
\path(0,-1) coordinate (d);
% Begin van theta(a):
\draw[blue] (13,9)--++(a)--++(b)--++(a)--++(b)--++(a)--++(b)--++(a)--++(b)--++(a)--++(b)--++(a)--++(b)--++(a)--++(b)--++(a); % trappetje van (11,7) tot (20,-2)
\draw[blue] (13,5)--++(d)--++(a)--++(d)--++(a)--++(d)--++(a)--++(d)--++(a)--++(d)--++(a)--++(d)--++(a)--++(d)--++(a);
% de grote loop
\draw[blue,thick] (20,16)--++(2,0)--++(1,1)--++(0,2)--++(1,1)--++(2,0)--++(1,-1)--++(0,-2)--++(-1,-1)--++(-2,0)--++(-1,-1)--++ (0,-2);
% het 3-CROSS!
\draw[blue,thick] (19,16.6)--++(0,2.8) (20,20)--++(2,0)--++(1,1)--++(0,2)  (23.6,24)--++(2.8,0)
(27,23)--++(0,-2)--++(1,-1)--++(2,0) (31,19.4)--++(0,-2.8) (30,16)--++(-2,0)--++(-1,-1)--++(0,-2) (23.6,12)--++(2.8,0);
% Pijlen in  het 3-CROSS!
 \draw[blue] (25,24) node[] {\tiny$<$};
 \draw[blue] (25,12) node[] {\tiny$>$};
 \draw[blue] (19,18) node[] {\tiny$\vee$};
 \draw[blue] (31,18) node[] {\tiny$\wedge$};
 \draw[blue] (21,20) node[] {\tiny$<$};
 \draw[blue] (23,22) node[] {\tiny$\vee$};
 \draw[blue] (27,22) node[] {\tiny$\wedge$};
 \draw[blue] (29,20) node[] {\tiny$<$};
 \draw[blue] (27,14) node[] {\tiny$\wedge$};
 \draw[blue] (29,16) node[] {\tiny$>$};
% Interne punten:
\node[fill=black,circle, inner sep=0.9pt] at (31,20) {};
\node[fill=black,circle, inner sep=0.9pt] at (31,16) {};
\node[fill=black,circle, inner sep=0.9pt] at (27,12) {};
\node[fill=black,circle, inner sep=0.9pt] at (19,20) {};
\node[fill=black,circle, inner sep=0.9pt] at (19,16) {};
\node[fill=black,circle, inner sep=0.9pt] at (23,12) {};
\node[fill=black,circle, inner sep=0.9pt] at (27,24) {};
\node[fill=black,circle, inner sep=0.9pt] at (23,24) {};
% Pijlen
 \draw[blue] (25,16) node[] {\tiny$<$};
 \draw[blue] (25,20) node[] {\tiny$>$};
 \draw[blue] (27,18) node[] {\tiny$\vee$};
 \draw[blue] (23,18) node[] {\tiny$\wedge$};
 \draw[blue] (21,16) node[] {\tiny$>$};
 \draw[blue] (23,14) node[] {\tiny$\vee$};
% Eind van theta(a):
\draw[blue] (31,20)--++(a)--++(b)--++(a)--++(b)--++(a)--++(b)--++(a)--++(b)--++(a)--++(b)--++(a)%nu de hoek om
--++(b)--++(c)--++(b)--++(c)--++(b)--++(c)--++(b)--++(c)--++(b)--++(c)--++(b)--++(c)--++(b)--++(c)--++(b)--++(c)--++(b)--++(c);
% hele theta(b)
\draw[green] (20,-2)--++(b)--++(a)--++(b)--++(a)--++(b)--++(a)--++(b)--++(a)--++(b)--++(a)--++(b)--++(a)--++(b)--++(a)--++(b)--++(a)%hoekje om
--++(d)--++(a)--++(d)--++(a)--++(d)--++(a)--++(d)--++(a)--++(d)--++(a)--++(d)--++(a)--++(d)--++(a)--++(d)--++(a)--++(d)--++(a)--++(d)++(a)++(d)++(a)++(d)++(a)
++(d)%verlenging in het midden
++(d)++(a)--++(d)--++(a)--++(d)--++(a)--++(d)--++(a)--++(d)--++(a)--++(d)--++(a)--++(d)--++(a)--++(d)--++(a)--++(d)--++(a)--++(d)--++(a)%hoekje om
--++(b)--++(a)--++(b)--++(a)--++(b)--++(a)--++(b)--++(a)--++(b)--++(a)--++(b)--++(a)--++(b);
% hele theta(c)
\draw[orange] (56,-10)--++(c)--++(b)--++(c)--++(b)--++(c)--++(b)--++(c)--++(b)--++(c)--++(b)--++(c)--++(b)--++(c)%hoekje om
--++(b)--++(a)--++(b)--++(a)--++(b)--++(a)--++(b)--++(a)--++(b)--++(a)--++(b)--++(a)--++(b)--++(a)--++(b)--++(a)--++(b)--++(a)
--++(b)++(a)++(b)++(a)++(b)++(a)++(b)++(a)++(b)--++(b)--++(a)%verlenging in het midden
--++(b)--++(a)--++(b)--++(a)--++(b)--++(a)--++(b)--++(a)--++(b)--++(a)--++(b)--++(a)--++(b)--++(a)--++(b)--++(a)--++(b)%hoekje om
--++(c)--++(b)--++(c)--++(b)--++(c)--++(b)--++(c)--++(b)--++(c)--++(b)--++(c)--++(b)--++(c);
% hele theta(d)
\draw[gray] (28,34)--++(b)--++(a)--++(b)--++(a)--++(b)--++(a)--++(b)--++(a)--++(b)--++(a)--++(b)--++(a)--++(b)--++(a)--++(b)--++(a)%hoekje om
--++(d)--++(a)--++(d)--++(a)--++(d)--++(a)--++(d)--++(a)--++(d)--++(a)--++(d)--++(a)--++(d)--++(a)--++(d)--++(a)--++(d)--++(a)--++(d)++(a)++(d)++(a)++(d)++(a)
++(d)%verlenging in het midden
++(d)++(a)--++(d)--++(a)--++(d)--++(a)--++(d)--++(a)--++(d)--++(a)--++(d)--++(a)--++(d)--++(a)--++(d)--++(a)--++(d)--++(a)--++(d)--++(a)%hoekje om
--++(b)--++(a)--++(b)--++(a)--++(b)--++(a)--++(b)--++(a)--++(b)--++(a)--++(b)--++(a)--++(b);
% hele theta(b) nummer 2
\draw[green] (-8,42)--++(b)--++(a)--++(b)--++(a)--++(b)--++(a)--++(b)--++(a)--++(b)--++(a)--++(b)--++(a)--++(b)--++(a)--++(b)--++(a)%hoekje om
--++(d)--++(a)--++(d)--++(a)--++(d)--++(a)--++(d)--++(a)--++(d)--++(a)--++(d)--++(a)--++(d)--++(a)--++(d)--++(a)--++(d)--++(a)--++(d)++(a)++(d)++(a)++(d)++(a)
++(d)%verlenging in het midden
++(d)++(a)--++(d)--++(a)--++(d)--++(a)--++(d)--++(a)--++(d)--++(a)--++(d)--++(a)--++(d)--++(a)--++(d)--++(a)--++(d)--++(a)--++(d)--++(a)%hoekje om
--++(b)--++(a)--++(b)--++(a)--++(b)--++(a)--++(b)--++(a)--++(b)--++(a)--++(b)--++(a)--++(b);
% hele theta(c) nummer 2
\draw[orange] (-16,6)--++(c)--++(b)--++(c)--++(b)--++(c)--++(b)--++(c)--++(b)--++(c)--++(b)--++(c)--++(b)--++(c)%hoekje om
--++(b)--++(a)--++(b)--++(a)--++(b)--++(a)--++(b)--++(a)--++(b)--++(a)--++(b)--++(a)--++(b)--++(a)--++(b)--++(a)--++(b)--++(a)
--++(b)++(a)++(b)++(a)++(b)++(a)++(b)++(a)++(b)--++(b)--++(a)%verlenging in het midden
--++(b)--++(a)--++(b)--++(a)--++(b)--++(a)--++(b)--++(a)--++(b)--++(a)--++(b)--++(a)--++(b)--++(a)--++(b)--++(a)--++(b)%hoekje om
--++(c)--++(b)--++(c)--++(b)--++(c)--++(b)--++(c)--++(b)--++(c)--++(b)--++(c)--++(b)--++(c);
% hele theta(d) nummer 2
\draw[gray] (-16,6)--++(b)--++(a)--++(b)--++(a)--++(b)--++(a)--++(b)--++(a)--++(b)--++(a)--++(b)--++(a)--++(b)--++(a)--++(b)--++(a)%hoekje om
--++(d)--++(a)--++(d)--++(a)--++(d)--++(a)--++(d)--++(a)--++(d)--++(a)--++(d)--++(a)--++(d)--++(a)--++(d)--++(a)--++(d)--++(a)--++(d)++(a)++(d)++(a)++(d)++(a)
++(d)%verlenging in het midden
++(d)++(a)--++(d)--++(a)--++(d)--++(a)--++(d)--++(a)--++(d)--++(a)--++(d)--++(a)--++(d)--++(a)--++(d)--++(a)--++(d)--++(a)--++(d)--++(a)%hoekje om
--++(b)--++(a)--++(b)--++(a)--++(b)--++(a)--++(b)--++(a)--++(b)--++(a)--++(b)--++(a)--++(b);
% Hoekpunten
\node[fill=black,circle, inner sep=0.85pt,left] at (-16,6) {};
\node[fill=black,circle, inner sep=0.85pt] at (-8,42) {};
\node[fill=red,circle, inner sep=1pt] at (20,-2) {};
\node[fill=black,circle, inner sep=0.85pt] at (28,34) {};
\node[fill=black,circle, inner sep=0.85pt] at (56,-10) {};
\node[fill=black,circle, inner sep=0.85pt] at (64,26) {};
% Namen:
 \draw (14,7) node[] {\footnotesize$\rK[\tha^8(a)]$};
 \draw (10,38) node[] {\footnotesize$\rK[\tha^8(b)]$};
 \draw (-13,24) node[] {\footnotesize$\rK[\tha^8(c)]$};
 \draw (1,2) node[] {\footnotesize$\rK[\tha^8(d)]$};
 \draw (37,-6) node[] {\footnotesize$\rK[\tha^8(b)]$};
 \draw (59,8) node[] {\footnotesize$\rK[\tha^8(c)]$};
 \draw (45,30) node[] {\footnotesize$\rK[\tha^8(d)]$};
% Origin:
 \draw (21.5,-3.5) node[] {\footnotesize $0$};
\end{tikzpicture}
\caption{\small A 3-cross in a figure eight configuration at level 8.}\label{fig:eight}
\end{figure}

With a \emph{figure-8} configuration we denote the subgraph of $\DD$ with 6 vertices and 7 edges pictured on the right.\\
Let $Q_\eight=[0,1]\cup[1,1+i]\cup\ldots\cup[-i,0]$ be the corresponding set of segments.\\
Now consider $K[\tha^8]([0,1])=K[\tha^8(a)]$ in  $K[\tha^8](Q_\eight)$, see Figure~\ref{fig:eight}.
Since  all the segments of $K[\tha^8]([0,1])$ lie in the interior of the (rounded off) closed curve corresponding to $K[\tha^8](Q_\eight)\setminus K[\tha^8]([0,1])$, the eight boundary points of the 3-cross are interior points of this curve, and so it must contain a 5-cross.

 Now note that $K[\tha^8](Q_\eight)$ (rounded off) appears in $\rK[\tha^{16}(a)]$, since $\rK[\tha^8(a)]$ contains a figure-8 configuration (as a subgraph of the 3-cross, both for the clockwise and for the anti-clockwise type).
Thus  $\rK[\tha^{16}(a)]$ has to contain a 5-cross. The same type of argument yields that $\rK[\tha^{24}(a)]$ contains a 7-cross, and more generally that $\rK[\tha^{8k}(a)]$ will contain a $(2k+1)$-cross, for each natural number $k$.\qed

\begin{theorem}\label{th:folstarfol}
Let $S$ and $T$ be two self-avoiding planefilling folding strings. Then $S*T$ is a self-avoiding planefilling folding string.
\end{theorem}

\noindent {\bf Proof:} This follows directly from Lemma~\ref{lem:conv-max} and Theorem~\ref{th:planefill}.\qed

 \section{Existence of planefilling folds}\label{sec:exist}

Let $z=f(\tha(a))$ be the endpoint of the curve associated to a $m$-fold. By Lemma \ref{lem:dir} in Section \ref{sec:classify} the sets $\rK[\tha(a)]$ and $\rK[\tha(c)]+(1+i)z$ are congruent, and so are $\rK[\tha(b)]+z$ and $\rK[\tha(d)]+iz$. This implies that the interior of the $\tha$-loop and the square through $0,z,(1+i)z$ and $iz$ have the same area. For a maximally simple $\tha$-loop this area is equal to $m$ (ignoring second order effects caused by rounding off), since all $4m$ segments of such a $\tha$-loop contribute one side to a unit square.
It thus follows from Theorem~\ref{th:planefill} that self-avoiding planefilling $m$-folds can only exist if $m=|z|^2$ is a sum of two squares.

\begin{theorem}\label{exist}
For any integer $m>1$ which is a sum of two squares there exists a self-avoiding and planefilling folding string of length $m$.
\end{theorem}

\noindent {\bf Proof:} We will prove in fact more, namely that for every Gaussian integer $z=u+iv$ with $u^2+v^2>1$ there is a self-avoiding   planefilling folding string whose curve ends in $z$.  By Theorem~\ref{th:folstarfol} and Lemma \ref{endstar}, it suffices to prove this  for those $z=u+iv$ which are  Gaussian prime numbers. Moreover, we can assume, by rotating if necessary, that $u>0, v\ge 0$. We already know that there exists a planefilling folding string to the Gaussian prime number $z=1+i$. Moreover, $z=u+iv$ is certainly not prime if both $u$ and $v$ are even, and since
$$u+iv=\big(1+i\big)\Big(\frac{u+v}{2}+i\,\frac{v-u}{2}\Big),$$
$z$ is also not prime if $u$ and $v$ are both odd (and not both equal to 1). So we may assume in the sequel that $z=u+iv$ is an odd Gaussian integer, i.e., $u+v$ is odd.
Moreover, by exchanging $u$ and $v$ if necessary, we can assume that $u$ is odd and $v$ even.

First consider the case $v\ne 0$.
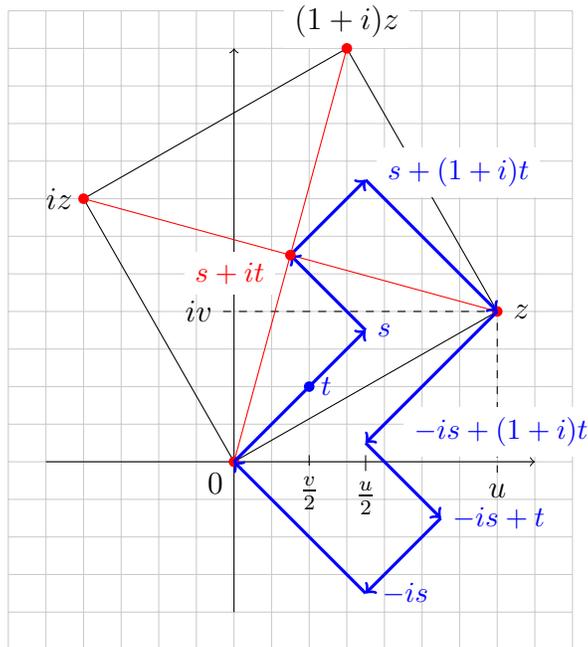
\begin{figure}[!h]
\centering
\begin{tikzpicture}[scale=0.5]
%[->]
% draw grid and axes
\draw[black!20] (-6,-5) grid (9,12); \draw[->] (-5,0)--(8,0) node (xaxis) [below right] {$$};
\draw[->] (0,-4)--(0,11) node (yaxis) [left] {$$};
% draw square and circle nodes 0 and z
\draw (0,0) node[below left] {$0$} coordinate (N) -- (-4,7) node[left] {$iz$} coordinate (A) -- (3,11) node[above, fill=white] {$(1+i)z$} coordinate (B) -- (7,4) node[right] {$\,z$} coordinate (z) -- (0,0);
\draw[dashed] (yaxis |- z) node[left] {$iv$} -| (xaxis -| z) node[below] {$u$};
\fill[red] (z)  circle (4pt); \fill[red] (0,0)  circle (4pt); \fill[blue] (2,2)  circle (4pt);
% calculate and plot intersection of red lines 's+it' (spit)
\coordinate (spit) at (intersection of N--B and A--z);
% draw red diagonals
\draw [color=red](0,0) -- (3,11); \draw [color=red](-4,7) -- (7,4);
% plot octagon
    \draw [blue,very thick] (N) -- (2,2) coordinate (t) node[right] {\small$t$};
    \draw [blue,->,very thick] (t) -- (3.5,3.5) coordinate (s) node[right] {\small$s$};
    \draw [blue,->,very thick] (s) -- (spit);
    \draw [blue,->,very thick] (spit) -- (3.5,7.5);
    \draw [blue] (3.8,7.7) node[right,fill=white] {\small$s+(1+i)t$};
    \draw [blue,->,very thick] (3.5,7.5) -- (z);
    \draw [blue,->,very thick] (z) -- (3.5,0.5);
    \draw [blue] (4.5,0.8) node[right,fill=white] {\small $-is+(1+i)t$}; % fill=white HAALT  DE x-AS WEG!
    \draw [blue,->,very thick] (3.5,0.5) -- (5.5,-1.5) node[right] {\small$-is+t$};
    \draw [blue,->,very thick] (5.5,-1.5) -- (3.5,-3.5) node[right] {\small$\,-is$};
    \draw [blue,->,very thick] (3.5,-3.5) -- (N);
% \draw[dashed] (yaxis |- s) node[left] {$$}
%        -| (xaxis -| s) node[below] {$\frac{a}{2}$};
    \draw (2,5pt) -- (2,-3pt)  node[anchor=north] {$\frac{v}{2}$};
    \draw (3.5,5pt) -- (3.5,-7pt)  node[anchor=north] {$\frac{u}{2}$};
%   plot 2 red circles
\fill[red] (-4,7) circle (4pt); \fill[red] (3,11) circle (4pt);
% plot intersection of red lines 's+it' (spit)
\fill[red] (spit) circle (4pt); \draw[red] (1.1,5.6) node[below left,fill=white] {\small$s+it$};
\end{tikzpicture}
\caption{\small The octagon $\Omega$ associated to $z=u+iv$ when $u=7$ and $v=4$.} \label{fig:oct}
\end{figure}

We will define an octagon $\Omega$,  based on the \midg \;$s=\frac12u(1+i)$ and the grid point $t=\frac12v(1+i)$ (see Figure~\ref{fig:oct}). The octagon is given by its vertices:
$$\Omega \thickapprox (0,\, s,\, s+it,\, s+(1+i)t,\, z,\, -is+(1+i)t,\, -is+t,\, -is).$$
The regular way in which $\Omega$ is constructed ensures that it is simply connected.
Now consider the finite directed subgraph of $\DD$ whose edges lie in  $\Omega$. It is easily checked that all its vertices have in-degree 1 and out-degree 1 or in-degree 2 and out-degree 2, except the vertices 0 and $z$, which have out-degree 1, respectively in-degree 1 (here we use again that $u+v$ is odd).

\begin{figure}[!h]% Octagon All
\centering
% 0, 90, 180, 270 degrees rotated to the left
\begin{tikzpicture}[scale=.3]
% draw grid and axes
    \draw[black!20] (-9,-5) grid (12,17);
    \draw[->] (-8,0)--(10,0) node (xaxis) [below right] {$$};
    \draw[->] (0,-4)--(0,15) node (yaxis) [left] {$$};
% draw square
    \draw (0,0) node[below left] {$0$} coordinate (N) -- (-4,7) node[above left] {$iz$} coordinate (A) -- (3,11) node[above right] {$(1+i)z$} coordinate (B) -- (7,4) node[right] {$z$} coordinate (z) -- (0,0);
    \draw[dashed] (yaxis |- z) node[left] {}
        -| (xaxis -| z) node[below] {};
% draw red diagonals
    \draw [color=red](0,0) -- (3,11);
    \draw [color=red](-4,7) -- (7,4);
% calculate and plot intersection of red lines 's+it' (spit)
    \coordinate (spit) at (intersection of N--B and A--z);
    \fill[red] (spit) circle (2pt) node[right] {$s+it$};
% plot octagon1
    \draw[blue,->,very thick] (0,0) -- ++(2,2) coordinate (t) node[right] {$$} -- ++(1.5,1.5) coordinate (s) node[right] {$$} -- ++(-2,2) -- ++(2,2) -- ++(3.5,-3.5) -- ++(-3.5,-3.5) -- ++(2,-2) -- ++(-2,-2) -- ++(-3.5,3.5);
% plot octagon2
   \pgftransformrotate{90};
   \pgftransformxshift{4cm};
   \pgftransformyshift{-7cm};
    \draw[blue,->,very thick] (0,0) -- ++(2,2) coordinate (t) node[right] {$$} -- ++(1.5,1.5) coordinate (s) node[right] {$$} -- ++(-2,2) -- ++(2,2) -- ++(3.5,-3.5) -- ++(-3.5,-3.5) -- ++(2,-2) -- ++(-2,-2) -- ++(-3.5,3.5);
% plot octagon3
   \pgftransformrotate{90};
   \pgftransformxshift{4cm};
   \pgftransformyshift{-7cm};
    \draw[blue,->,very thick] (0,0)--++(2,2) coordinate (t) node[right] {$$}--++(1.5,1.5) coordinate (s)
    node[right] {$$}--++(-2,2)--++(2,2)--++(3.5,-3.5)--++(-3.5,-3.5)--++(2,-2)--++(-2,-2)--++(-3.5,3.5);
% plot octagon2
   \pgftransformrotate{90};
   \pgftransformxshift{4cm};
   \pgftransformyshift{-7cm};
    \draw[blue,->,very thick] (0,0) -- ++(2,2) coordinate (t) node[right] {$$} -- ++(1.5,1.5) coordinate (s) node[right] {$$} -- ++(-2,2) -- ++(2,2) -- ++(3.5,-3.5) -- ++(-3.5,-3.5) -- ++(2,-2) -- ++(-2,-2) -- ++(-3.5,3.5);
\end{tikzpicture}
\caption{\small Four octagons forming a $\tha$-loop region.} \label{fig:oct4}
\end{figure}
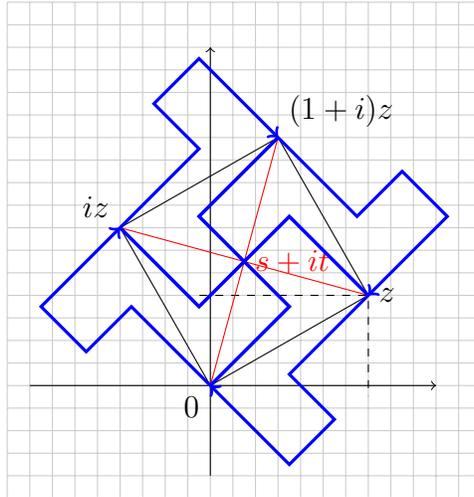

 Now remember that an directed Euler path is a path along a connected directed graph that connects all the vertices and that traverses every edge of the graph only once, and that
    the degree conditions above yield  the existence of a directed Euler path  from 0 to $z$ (\cite{Euler}).
    But such a path is nothing else than the curve given by a self-avoiding folding string $S$ of length $u^2+v^2$, since the octagon $\Omega$ has half the area of  the square $\Xi$ with vertices 0, $z, (1+i)z$, and  $iz$.
    Moreover, the curve obtained by joining the four curves produced by rotating thrice over $90^o$ around
    the center $s+it=(1+i)z/2$ of $\Xi$ is the $\theta$-loop of $S$, and it visits each of its internal
    points exactly twice (cf.~Figure~\ref{fig:oct4}). It follows from Theorem~\ref{th:planefill} that the
    string $S$ generates a planefilling curve.

    When $v=0$ the octagon $\Omega$ degenerates into a square, but the arguments will remain the same.       \qed

%Mid-grid borders can be much more complicated than the octagon $\Omega$ defined in the proof above.
 It can be seen that actually one can find a \emph{symmetric} fold to $u+iv$ for any odd Gaussian integer. The idea here is that Euler's theorem will guarantee a $\DD$-word $w$ of length $(m-1)/2$ with its curve filling half of $\Omega$ after the central mid-grid square has been removed. The symmetric $m$-fold will then be given by $w\,a\,\tau(w)$ if $u+v-1$ is divisible by 4, or else by $w\,c\,\tau(w)$.

\section{Perfect sequences: the carrousel theorem}\label{sec:classify}

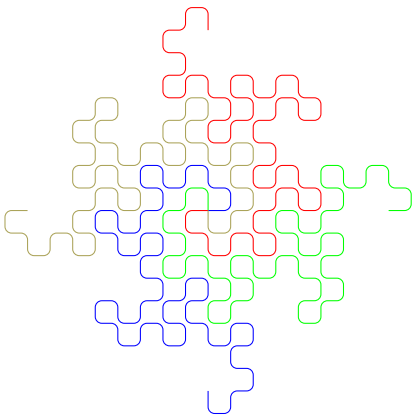
\begin{figure}[!h]% Carrousel
\centering
\begin{tikzpicture}[scale=.3,rounded corners=2.5pt]
\path(1,0)  coordinate (a);
\path(0,1)  coordinate (b);
\path(-1,0) coordinate (c);
\path(0,-1) coordinate (d);
%thna:
\draw[blue] (0,0)--++(a)--++(b)--++(c)--++(b)--++(c)--++(d)--++(c)--++(b)--++(c)--++(d)--++(a)--++(d)--++(c)--++(d)--++
(c)--++(b)--++(c)--++(d)--++(a)--++(d)--++(a)--++(b)--++(a)--++(d)--++(c)--++(d)--++(a)--++(d)--++(c)--++(d)--++(c)--++
(b)--++(c)--++(d)--++(a)--++(d)--++(a)--++(b)--++(a)--++(d)--++(a)--++(b)--++(c)--++(b)--++(a)--++(b)--++(a)--++(d)--++
(c)--++(d)--++(a)--++(d)--++(a)--++(b)--++(a)--++(d)--++(c)--++(d)--++(a)--++(d)--++(c)--++(d)--++(c)--++(b);
%thnb:
\draw[green] (0,0)--++(b)--++(c)--++(d)--++(c)--++(d)--++(a)--++(d)--++(c)--++(d)--++(a)--++(b)--++(a)--++(d)--++(a)--++
(d)--++(c)--++(d)--++(a)--++(b)--++(a)--++(b)--++(c)--++(b)--++(a)--++(d)--++(a)--++ (b)--++(a)--++(d)--++(a)--++(d)--++
(c)--++(d)--++(a)--++(b)--++(a)--++(b)--++(c)--++(b)--++(a)--++(b)--++(c)--++(d)--++(c)--++(b)--++(c)--++(b)--++(a)--++
(d)--++(a)--++(b)--++(a)--++(b)--++(c)--++(b)--++(a)--++(d)--++(a)--++(b)--++(a)--++(d)--++(a)--++(d)--++(c);
%thnc:
\draw[red] (0,0)--++(c)--++(d)--++(a)--++(d)--++(a)--++(b)--++(a)--++(d)--++(a)--++(b)--++(c)--++(b)--++(a)--++(b)--++
(a)--++(d)--++ (a)--++(b)--++(c)--++(b)--++(c)--++(d)--++(c)--++(b)--++ (a)--++(b)--++(c)--++(b)--++(a)--++(b)--++(a)--++
(d)--++(a)--++(b)--++(c)--++(b)--++(c)--++(d)--++(c)--++(b)--++(c)--++(d)--++(a)--++(d)--++(c)--++(d)--++(c)--++(b)--++
(a)--++(b)--++(c)--++(b)--++(c)--++(d)--++(c)--++(b)--++(a)--++(b)--++(c)--++(b)--++(a)--++(b)--++(a)--++(d);
%thnd:
\draw[yellow!60!black] (0,0)--++(d)--++(a)--++(b)--++(a)--++(b)--++(c)--++(b)--++(a)--++(b)--++(c)--++(d)--++(c)--++(b)--++(c)--++ (b)--++(a)--++(b)--++(c)--++ (d)--++(c)--++(d)--++(a)--++(d)--++(c)--++(b)--++(c)--++(d)--++
(c)--++(b)--++(c)--++(b)--++(a)--++(b)--++(c)--++(d)--++(c)--++(d)--++(a)--++(d)--++(c)--++(d)--++ (a)--++(b)--++(a)--++(d)--++(a)--++(d)--++ (c)--++(b)--++(c)--++(d)--++(c)--++(d)--++(a)--++(d)--++(c)--++(b)--++(c)--++(d)--++(c)--++(b)--++(c)--++(b)--++(a);
\end{tikzpicture}
\caption{\small The carrousel of order 6 for $S=D$.} \label{fig:carr}
\end{figure}

A  string $S$ is called \emph{perfect} if it is self-avoiding, and the four curves $K[S^{*\infty}]$, $iK[S^{*\infty}]$, $-K[S^{*\infty}]$, and $-iK[S^{*\infty}]$ together cover every edge of $\DD$ exactly once. Davis and Knuth
 prove in \cite{DK} that the 2-fold $S=D$ has this property, stating that this is ``probably the most interesting property of the dragon curve". See Figure~\ref{fig:carr}.

 Actually it is remarkable that a fold $S$ can be  perfect, since the set above, which we call the \emph{carrousel}
associated to $S$, is incompatible with the directed graph structure of $\DD$---note, e.g., that the out-degree of the origin is four for the carrousel. Nevertheless, perfect folding sequences abound, and the following characterization holds.

\begin{theorem}\label{carrousel}
Let $S$ be a folding string that generates a self-avoiding planefilling curve. Then $S$ is perfect if and only if the curve $K[\tha(a)]$ of $S$ does \emph{not} traverse the segment from $(1+i)z-i-1$ to $(1+i)z-i$, where $z=f(\tha(a))$ is the endpoint of the curve of $S$.
\end{theorem}

We first present a lemma which resolves the incompatibility observed above.

\begin{lemma}\label{lem:dir}
Let $S$ be a folding string with $z=f(\tha(a))$. Then
$$iK[\tha(a)]=K[\tha(d)]+iz, \;i^2K[\tha(a)]=K[\tha(c)], \;i^3K[\tha(a)]=K[\tha(b)]+i^3z.$$
\end{lemma}

\noindent {\bf Proof:} Note that these are equalities of \emph{sets}, which we prove by going through the segments making up the sets in a certain order. To see that the first one holds, traverse
$iK[\tha(a)]=K[\sigma\tha(a)]$ from $iz$ to 0. This amounts to applying $\sigma^2\tau$ to the word $\sigma\tha(a)$ (since we reversed the direction the segments in the curve are rotated over $180^o$). But since $\sigma$ and $\tau$ commute, and by definition of the folding morphism we have
$$\sigma^2\tau(\sigma\tha(a))=(\sigma\tau)^3\tha(a)=\tha(d).$$
For the second equality, note that $i^2K[\tha(a)]=K[\sigma^2\tha(a)]$, and that $$\sigma^2\tha(a)=(\sigma\tau)^2\tha(a)=\tha(c).$$
The third equality follows as the first one.\qed

\medskip

The segment $[1+i,i]$ will play a crucial role in the proof of Theorem \ref{carrousel}.

\begin{lemma}\label{lem:segment}
Let $S$ be a self-avoiding planefilling string. Then the $\tha$-loop of $S$ contains the
segment $[1+i,i]$.
\end{lemma}

\noindent {\bf Proof:}
The statement is true when $m=2$, and for $m\ge 3$ we either have $\tha(a)=abc\dots$,
and then $K[\tha(a)]$ already contains $[1+i,i]$, or $\tha(a)=aba\dots$, and then the point $1+i$ is in the interior of the $\tha$-loop, what implies that the segment $[1+i,i]$ must be traversed by it.
\qed

\medskip

\noindent{\bf Proof of  Theorem~\ref{carrousel}:}

Let $\tha$ be the morphism associated to $S$, with $\tha(a)=ab\ldots$.
 First we will translate the mysterious condition on the curve $K[\tha(a)]$ of $S$ to a slightly less mysterious condition on the $\tha$-loop of $S$.
Note that (by Lemma~\ref{lem:dir}) the $K[\tha(c)]$ part of the $\tha$-loop (ignoring the rounding off)
can be written as $K[\tha(c)]+(1+i)z=-K[\tha(a)]+(1+i)z,$
where $z=f(\tha(a))$ is the endpoint of $K[\tha(a)]$.
\begin{wrapfigure}{r}{35mm}
\vspace*{-0.8cm}
\begin{center}
\begin{tikzpicture}[scale=.4]
% variables
\def\tangle{45}
\def\rot{90}
% draw grid and axes
    \draw[black!20] (-3,-1) grid (5,7);
% nodes
\node[fill=black,circle, inner sep=.8pt] at (0,0){};
\draw (0,0) node[below left] {\small$0$}  (4,2) node[right] {\small$z$}  (2,6) node[above right] {\small$(1+i)z$};
\node[fill=black,circle, inner sep=.8pt] at (-2,4){};
\node[fill=black,circle, inner sep=.8pt] at (4,2){};
\node[fill=black,circle, inner sep=.8pt] at (2,6){};
% Arrows
\draw[red,->,very thick] (1,5) -- ++(1,0);
\node[fill=blue,circle, inner sep=.8pt] at (1,5){};
\node[fill=blue,circle, inner sep=.8pt] at (2,5){};
\draw[red,->,very thick] (1,1) -- ++(-1,0);
\node[fill=blue,circle, inner sep=.8pt] at (1,1){};
\node[fill=blue,circle, inner sep=.8pt] at (0,1){};
\end{tikzpicture}
\end{center}
\vspace*{-.7cm}
%\caption{\footnotesize{Translating a condition}}%\label{fig:mystery}
\end{wrapfigure}
 It follows that $K[\tha(a)]$ does go through the segment from $(1+i)z-i-1$ to $(1+i)z-i$ if and only if $K[\tha(c)]+(1+i)z$ goes through the segment from $1+i$ to $i$ (see also the figure on the right).

In the sequel, let $$Q_n:=K[\tha^n(a)]\,\cup\, iK[\tha^n(a)]\,\cup\, -K[\tha^n(a)]\,\cup\, -iK[\tha^n(a)]$$
be the carrousel of order $n$. Since we can replace the 2-fold $D$ by the 4-fold $D\ast D=DDU$, we will assume henceforth that $m\ge 3$.

\noindent \emph{Part 1:} Suppose $[1+i,i]$ is traversed by $K[\tha(c)]+(1+i)z$ in the $\tha$-loop.
\begin{wrapfigure}{l}{18mm}
\vspace*{-1.1cm}
\begin{center}
%    %%%%\includegraphics{toucan.eps}
\begin{tikzpicture}[domain=0:1,scale=0.4]\label{fig:segment}
\foreach \angle in {0,90,180,-90}
{  \path(0,0) coordinate (P0);
  \path(P0) ++(0+\angle:0.8 cm) coordinate (P1);
  \path(P1) ++(45+\angle:0.28284 cm) coordinate (P2);
  \path(P2) ++(90+\angle:0.6 cm) coordinate (P3);
  \path(P3) ++(45+\angle:0.28284 cm) coordinate (P4);
  \path(P4) ++(0+\angle:0.8 cm) coordinate (P5);
  \draw    (P0) -- (P1)
           arc (-90+\angle:\angle:0.2 cm) -- (P2)
           (P2) -- (P3)
           arc (180+\angle:90+\angle:0.2 cm) -- (P4)
           (P4) -- (P5);             }
\node  at (1.2,0.1) {\scriptsize{$1$}};
\node  at (0.9,1.35) {\scriptsize{$1\!+\!i$}};
\node  at (0.05,1.35) {\scriptsize{$i$}};
\node  at (1,0) {{$\cdot$}}; \node  at (1,1) {{$\cdot$}}; \node  at (0,1) {{$\cdot$}};
% draw arrows
        \draw[] (-1.5,-1) node[] {\tiny$<$};
        \draw[] (0.5,-1) node[] {\tiny$<$};
        \draw[] (-0.5,0) node[] {\tiny$<$};
        \draw[] (-0.5,1) node[] {\tiny$>$};
        \draw[] (0.5,0) node[] {\tiny$>$};
        \draw[] (1.5,1) node[] {\tiny$>$};
        \draw[] (-1,1.5) node[] {\tiny$\vee$};
        \draw[] (0,.5) node[] {\tiny$\vee$};
        \draw[] (-1,-.5) node[] {\tiny$\vee$};
        \draw[] (1,-1.5) node[] {\tiny$\wedge$};
        \draw[] (0,-.5) node[] {\tiny$\wedge$};%OK
        \draw[] (1,.5) node[] {\tiny$\wedge$};%OK
\end{tikzpicture}
\end{center}
\vspace*{-1.2cm}
%\caption{\footnotesize{The segment $[1+i,i]$}}
\end{wrapfigure}
Then $\tha(a)$ must start with $aba$, so $\tha(d)$ must end in $dad$,  and in $Q_1$ the segment $[1+i,i]$ is absent (cf.~figure on the left). In $Q_2$, the segment $[1+i,i]$ is still missing, since it could only be traversed by $K[\tha(c)]+(1+i)z$,
which is missing from the incomplete $\tha$-loop with $i$ as one of its internal points.
 Continuing in this fashion we see that the segment $[1+i,i]$ will be missing in all $Q_n$, and so not all edges of $\DD$ are visited.

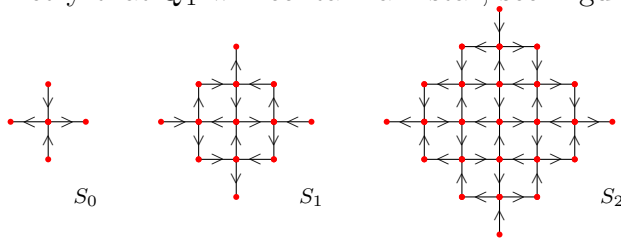
\begin{figure}[!h]
\centering\vspace*{-.5cm}
\begin{tikzpicture}[scale=.5]
\path(1,0)  coordinate (a); \path(0,1)  coordinate (b);
\path(-1,0) coordinate (c); \path(0,-1) coordinate (d);
%==========================================================
\def\arrowa((#1,#2)){\draw (#1,#2)--++(a);
\node[fill=red,circle, inner sep=.8pt] at (#1,#2)  {};
\node[fill=red,circle, inner sep=.8pt] at ({#1+1},#2)  {};
\draw[] ({#1+0.5},#2) node[] {\tiny$>$};}
\def\arrowb((#1,#2)){\draw (#1,#2)--++(b);
\node[fill=red,circle, inner sep=.8pt] at (#1,#2)  {};
\node[fill=red,circle, inner sep=.8pt] at (#1,{#2+1})  {};
\draw[] ({#1},{#2+0.5}) node[] {\tiny $\wedge$};}
\def\arrowc((#1,#2)){\draw (#1,#2)--++(c);
\node[fill=red,circle, inner sep=.8pt] at (#1,#2)  {};
\node[fill=red,circle, inner sep=.8pt] at ({#1-1},#2) {};
\draw[] ({#1-0.5},#2) node[] {\tiny $<$};}
\def\arrowd((#1,#2)){\draw (#1,#2)--++(d);
\node[fill=red,circle, inner sep=.8pt] at (#1,#2)  {};
\node[fill=red,circle, inner sep=.8pt] at (#1,{#2-1}) {};
\draw[] (#1,{#2-0.5}) node[] {\tiny $\vee$};}%
\def\starn((#1,#2)){\arrowa(({#1-1},#2));\arrowb((#1,#2));\arrowc(({#1+1},#2));\arrowd((#1,#2));}
\def\starp((#1,#2)){\arrowa((#1,#2));\arrowb((#1,{#2-1}));\arrowc((#1,#2));\arrowd((#1,{#2+1}));}
% 0-star:
\starp((0,0));
% 1-star:
\starn((4,0));\starn((6,0));\starn((5,1));\starn((5,-1));
% 2-star:
\starp((10,0));\starp((12,0));\starp((14,0));\starp((11,1));\starp((13,1));
\starp((12,2));\starp((11,-1));\starp((13,-1));\starp((12,-2));
\node  at (1,-2) {\scriptsize{$S_0$}};\node  at (7,-2) {\scriptsize{$S_1$}};\node  at (15,-2) {\scriptsize{$S_2$}};
\end{tikzpicture}
\caption{\small The 0-star, 1-star and 2-star graphs.}\label{fig:squares}
\end{figure}

\noindent \emph{Part 2:} We now  suppose $[1+i,i]$ is \emph{not} traversed by $K[\tha(c)]+(1+i)z$ in the $\tha$-loop.
 Note that either $\tha(a)$ starts with $abc$, or with $aba$. The latter case will be considered at the end of the proof. If $\tha(a)$ starts with $abc$, then the 3$^{rd}$ segment of $K[\tha(a)]$ traverses $[1+i,i]$. Since $K[\tha(a)]$ does not end in the point $i$, it will then also traverse the segment $[i,2i]$. It thus follows from rotational symmetry that $Q_1$ will contain a 1-star, see Figure~\ref{fig:squares}.

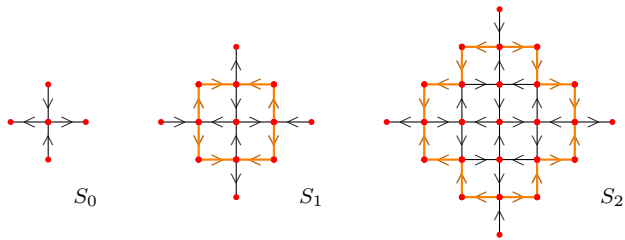
\begin{figure}[!h]
\centering
\begin{tikzpicture}[scale=.5]
\path(1,0)  coordinate (a); \path(0,1)  coordinate (b);
\path(-1,0) coordinate (c); \path(0,-1) coordinate (d);
%==========================================================
\def\arrowa((#1,#2)){\draw (#1,#2)--++(a);
\node[fill=red,circle, inner sep=.8pt] at (#1,#2){};
\node[fill=red,circle, inner sep=.8pt] at ({#1+1},#2){};
\draw[] ({#1+0.5},#2) node[] {\tiny$>$};}%
\def\arrowb((#1,#2)){\draw (#1,#2)--++(b);
\node[fill=red,circle, inner sep=.8pt] at (#1,#2)  {};
\node[fill=red,circle, inner sep=.8pt] at (#1,{#2+1})  {};
\draw[] ({#1},{#2+0.5}) node[] {\tiny $\wedge$};}%
\def\arrowc((#1,#2)){\draw (#1,#2)--++(c);
\node[fill=red,circle, inner sep=.8pt] at (#1,#2)  {};
\node[fill=red,circle, inner sep=.8pt] at ({#1-1},#2) {};
\draw[] ({#1-0.5},#2) node[] {\tiny $<$};}%
\def\arrowd((#1,#2)){\draw (#1,#2)--++(d);
\node[fill=red,circle, inner sep=.8pt] at (#1,#2)  {};
\node[fill=red,circle, inner sep=.8pt] at (#1,{#2-1}) {};
\draw[] (#1,{#2-0.5}) node[] {\tiny $\vee$};}%
% Nu met een andere kleur
\def\bouna((#1,#2)){\draw[orange,thick] (#1,#2)--++(a);
\node[fill=red,circle, inner sep=.8pt] at (#1,#2){};
\node[fill=red,circle, inner sep=.8pt] at ({#1+1},#2){};
\draw[orange] ({#1+0.5},#2) node[] {\tiny$>$};}%
\def\bounb((#1,#2)){\draw[orange,thick] (#1,#2)--++(b);
\node[fill=red,circle, inner sep=.8pt] at (#1,#2)  {};
\node[fill=red,circle, inner sep=.8pt] at (#1,{#2+1})  {};
\draw[orange] ({#1},{#2+0.5}) node[] {\tiny $\wedge$};}%
\def\bounc((#1,#2)){\draw[orange,thick] (#1,#2)--++(c);
\node[fill=red,circle, inner sep=.8pt] at (#1,#2)  {};
\node[fill=red,circle, inner sep=.8pt] at ({#1-1},#2) {};
\draw[orange] ({#1-0.5},#2) node[] {\tiny $<$};}%
\def\bound((#1,#2)){\draw[orange,thick] (#1,#2)--++(d);
\node[fill=red,circle, inner sep=.8pt] at (#1,#2)  {};
\node[fill=red,circle, inner sep=.8pt] at (#1,{#2-1}) {};
\draw[orange] (#1,{#2-0.5}) node[] {\tiny $\vee$};}%

\def\starn((#1,#2)){\arrowa(({#1-1},#2));\arrowb((#1,#2));\arrowc(({#1+1},#2));\arrowd((#1,#2));}
\def\starp((#1,#2)){\arrowa((#1,#2));\arrowb((#1,{#2-1}));\arrowc((#1,#2));\arrowd((#1,{#2+1}));}
% 0-star:
\starp((0,0));
% 1-star:
\starn((4,0));\starn((6,0));\starn((5,1));\starn((5,-1));
% 2-star:
\starp((10,0));\starp((12,0));\starp((14,0));\starp((11,1));\starp((13,1));
\starp((12,2));\starp((11,-1));\starp((13,-1));\starp((12,-2));
% boundary of 1-star:
\bouna((4,1));\bouna((4,-1));\bounb((4,0));\bounb((6,0));
\bounc((6,1));\bounc((6,-1));\bound((4,0));\bound((6,0));
% boundary of 2-star:
\bounb((10,-1));\bound((10,1));\bounb((14,-1));\bound((14,1));
\bounb((11,-2));\bound((11,2));\bounb((13,-2));\bound((13,2));
\bouna((12,2));\bouna((12,-2));\bouna((13,1));\bouna((13,-1));
\bounc((11,1));\bounc((11,-1));\bounc((12,2));\bounc((12,-2));
\node  at (1,-2) {\scriptsize{$S_0$}};\node  at (7,-2) {\scriptsize{$S_1$}};\node  at (15,-2) {\scriptsize{$S_2$}};
\end{tikzpicture}
\caption{\small The boundaries in the 1-star and 2-star graphs.}\label{fig:starbound}
\end{figure}

In general the \emph{$k$-star graph} is the subgraph $S_k$ of $\DD$ consisting of all the vertices and edges lying in the square through $k+1$, $(k+1)i$, $-k\!-\!1$ and $-(k+1)i$. It will be important later that we have the following property of the star graphs: one obtains $S_{k+1}$ from $S_k$ by adding all edges (and corresponding vertices) to the vertices that are not yet of degree 4 in $S_k$. The edges and vertices that are added, except for the outer 4, will be called the boundary $B_k$ of $S_k$ (see also Figure~\ref{fig:starbound}).

\begin{figure}[!h]
\centering
\begin{tikzpicture}[scale=.3,rounded corners=3pt]
\path(1,0)  coordinate (a); \path(0,1)  coordinate (b);
\path(-1,0) coordinate (c); \path(0,-1) coordinate (d);
%==========================================================
\def\arrowa((#1,#2)){\draw (#1,#2)--++(a);
\node[fill=red,circle, inner sep=.4pt] at (#1,#2)  {};
\node[fill=red,circle, inner sep=.4pt] at ({#1+1},#2)  {};}%
\def\arrowb((#1,#2)){\draw (#1,#2)--++(b);
\node[fill=red,circle, inner sep=.4pt] at (#1,#2)  {};
\node[fill=red,circle, inner sep=.4pt] at (#1,{#2+1})  {};}%
\def\arrowc((#1,#2)){\draw (#1,#2)--++(c);
\node[fill=red,circle, inner sep=.4pt] at (#1,#2)  {};
\node[fill=red,circle, inner sep=.4pt] at ({#1-1},#2) {};}%
\def\arrowd((#1,#2)){\draw (#1,#2)--++(d);
\node[fill=red,circle, inner sep=.4pt] at (#1,#2)  {};
\node[fill=red,circle, inner sep=.4pt] at (#1,{#2-1}) {};}
\def\starn((#1,#2)){\arrowa(({#1-1},#2));\arrowb((#1,#2));\arrowc(({#1+1},#2));\arrowd((#1,#2));}
\def\starp((#1,#2)){\arrowa((#1,#2));\arrowb((#1,{#2-1}));\arrowc((#1,#2));\arrowd((#1,{#2+1}));}
% 1-star:
\starn((-1,7));\starn((1,7));\starn((0,8));\starn((0,6));
% Q_1:
\draw[blue] (0,7)--++(a)--++(b)--++(c)--++(b)--++(a)--++(d)--++(a)--++(b);
\draw[blue] (0,7)--++(b)--++(c)--++(d)--++(c)--++(b)--++(a)--++(b)--++(c);
\draw[blue] (0,7)--++(c)--++(d)--++(a)--++(d)--++(c)--++(b)--++(c)--++(d);
\draw[blue] (0,7)--++(d)--++(a)--++(b)--++(a)--++(d)--++(c)--++(d)--++(a);
\node  at (0,0) {\scriptsize{$Q_1$}};
\end{tikzpicture}\hspace{2cm}
\begin{tikzpicture}[scale=.3,rounded corners=3pt]
\path(1,0)  coordinate (a); \path(0,1)  coordinate (b);
\path(-1,0) coordinate (c); \path(0,-1) coordinate (d);
%==========================================================
\def\arrowa((#1,#2)){\draw (#1,#2)--++(a);
\node[fill=red,circle, inner sep=.4pt] at (#1,#2)  {};
\node[fill=red,circle, inner sep=.4pt] at ({#1+1},#2)  {};}%
\def\arrowb((#1,#2)){\draw (#1,#2)--++(b);
\node[fill=red,circle, inner sep=.4pt] at (#1,#2)  {};
\node[fill=red,circle, inner sep=.4pt] at (#1,{#2+1})  {};}%
\def\arrowc((#1,#2)){\draw (#1,#2)--++(c);
\node[fill=red,circle, inner sep=.4pt] at (#1,#2)  {};
\node[fill=red,circle, inner sep=.4pt] at ({#1-1},#2) {};}%
\def\arrowd((#1,#2)){\draw (#1,#2)--++(d);
\node[fill=red,circle, inner sep=.4pt] at (#1,#2)  {};
\node[fill=red,circle, inner sep=.4pt] at (#1,{#2-1}) {};}
\def\starn((#1,#2)){\arrowa(({#1-1},#2));\arrowb((#1,#2));\arrowc(({#1+1},#2));\arrowd((#1,#2));}
% Q_2:
\draw[blue] (0,0)--++(a)--++(b)--++(c)--++(b)--++(a)--++(d)--++(a)--++(b)--++(c)--++(b)--++(a)--++(b)--++(c)--++(d)--++(c) --++(b)--++(c)--++(d)--++(a)--++(d)--++(c)--++(b)--++(c)--++(d)--++(c)--++(b)--++(a)--++(b)--++(c) --++(d)--++(c)--++(b)--++(a)--++(b)--++(c)--++(b) --++(a)--++(d)--++(a)--++(b)--++(a)--++(d)--++(c)        --++(d)--++(a)--++(b)--++(a)--++(d) --++(a)--++(b)--++(c) --++(b)--++(a)--++(d)--++(a)--++(b)--++(c)--++(b)--++(a)--++(b)--++(c)--++(d) --++(c)--++(b)
(0,0)--++(b)--++(c)--++(d)--++(c)--++(b)--++(a)--++(b)--++(c)--++(d)--++(c)--++(b)--++(c)--++(d)--++(a)--++(d) --++(c)--++(d)--++(a)--++(b)--++(a)--++(d)--++(c)        --++(d)--++(a)--++(d)--++(c)--++(b)--++(c)--++(d)--++(a)--++(d)--++(c)--++(b)--++(c)--++(d) --++(c)--++(b)--++(a)--++(b)--++(c)--++(b) --++(a)--++(d)        --++(a)--++(b)--++(c)--++(b)--++(a)--++(b)--++(c)--++(d)--++(c)--++(b)--++(a)--++(b)--++(c)--++(d) --++(c)--++(b)--++(c)--++(d)--++(a)--++(d)--++(c)
(0,0)--++(c)--++(d)--++(a)--++(d)--++(c)--++(b)--++(c)--++(d)--++(a)--++(d)--++(c)--++(d) --++(a)--++(b)--++(a)--++(d)--++(a)--++(b)--++(c)--++(b) --++(a)--++(d)        --++(a)--++(b)--++(a)--++(d)--++(c)--++(d)--++(a)--++(b)--++(a)--++(d)--++(c)--++(d)--++(a)--++(d) --++(c)--++(b)--++(c)--++(d)--++(c)--++(b) --++(a)--++(b)--++(c)--++(d) --++(c)--++(b)--++(c)--++(d)--++(a)--++(d)--++(c)--++(b)--++(c)--++(d)--++(a)--++(d)--++(c)--++(d) --++(a)--++(b)--++(a)--++(d)
(0,0)--++(d)--++(a)--++(b)--++(a)--++(d)--++(c)--++(d)--++(a)--++(b)--++(a)--++(d)--++(a)--++(b)--++(c) --++(b)
--++(a)--++(b)--++(c)--++(d)--++(c)--++(b)--++(a)        --++(b)--++(c)--++(b)--++(a)--++(d) --++(a)--++(b)
--++(c)--++(b)--++(a)--++(d)--++(a)--++(b)--++(a)--++(d)--++(c)--++(d)--++(a)--++(d) --++(c)--++(b)--++(c)
--++(d)--++(a)--++(d)--++(c)--++(d)--++(a)--++(b)--++(a)--++(d)--++(c)--++(d) --++(a)--++(b)--++(a)--++(d)
--++(a)--++(b)--++(c)--++(b)--++(a);
% 1-star:
\starn((-1,0));\starn((1,0));\starn((0,1));\starn((0,-1));
% Nu met een andere kleur: oranje de boundary image
\def\thaba((#1,#2)){\draw[orange,thick] (#1,#2)--++(c)--++(b)--++(a)--++(b)--++(c)--++(d)--++(c)--++(b)
--++(a)--++(b)--++(c)--++(b)--++(a)--++(d)--++(a)--++(b); }
\def\thabc((#1,#2)){\draw[orange,thick] (#1,#2)--++(c)--++(b)--++(a)--++(b)--++(c)--++(d)--++(c)--++(b)
--++(c)--++(d)--++(a)--++(d)--++(c)--++(b)--++(c)--++(d); }
\def\thadc((#1,#2)){\draw[orange,thick] (#1,#2)--++(a)--++(d)--++(c)--++(d)--++(a)--++(b)--++(a)--++(d)
--++(c)--++(d)--++(a)--++(d)--++(c)--++(b)--++(c)--++(d); }
\def\thada((#1,#2)){\draw[orange,thick] (#1,#2)--++(a)--++(d)--++(c)--++(d)--++(a)--++(b)--++(a)--++(d)
--++(a)--++(b)--++(c)--++(b)--++(a)--++(d)--++(a)--++(b); }
\thaba((-2,-2));  \thabc((2,2));  \thadc((2,2)); \thada((-2,-2));
\node  at (-4,-7) {\scriptsize{$Q_2$}};
\end{tikzpicture}
\caption{\small The case $\tha(a)=abcbadab$; the carrousel $Q_1$ and $Q_2$ with the 1-star. On the right the $K$-$\tha$-image of $B_1$ is indicated in orange.}\label{fig:stargrowth}
\end{figure}
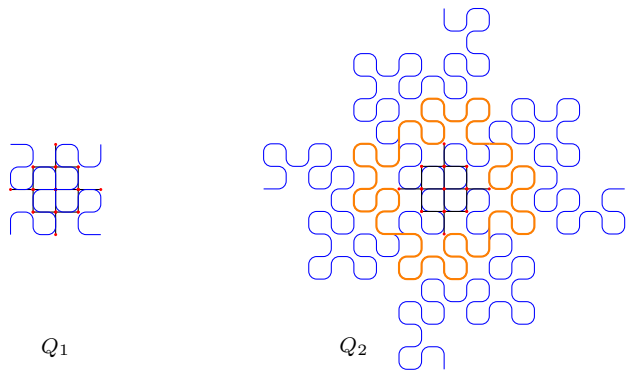

Now note that $B_1$ is a closed simple loop in $S_1$, which has all the segments of $S_0$ in its interior, and that we know that $Q_1$ contains $S_1$. It follows (recall Equation (\ref{eq:Kthimage}) and (\ref{eq:Kthimageunion})) that  $Q_2$ contains the $K$-$\tha$-image of $B_1$, which will be a closed (simple) loop, that contains $S_1$ in its interior (except for the 8 extremal points), since the $K$-$\tha$-image of $S_0$ is $Q_1$. Now we get to the main idea of the argument: since the $\tha$-loop is maximally simple, all grid points within the $K$-$\tha$-image of $B_1$ will be visited twice by $Q_2$.
But then $S_1$ has to grow\footnote{There is one special case ($m=4$ and $\tha(a)=abcb$) where this is not true, due to the fact that $B_1$ does not correspond to a $\DD$-word and because for this particular $\tha$ the $K$-$\tha$-image of $S_0$ \emph{equals} $S_1$.} till $S_2$, and so $Q_2$ will contain $S_2$. Continuing in this fashion we find that $Q_k$ contains $S_k$ for all $k$, and so the $Q_k$ increase to  cover all of $\Z_2$. See Figure~\ref{fig:stargrowth} for an example with a self-avoiding and planefilling 8-fold.

We still have to consider the possibility that $\tha(a)$ starts with $aba$. By Lemma~\ref{lem:segment} one of the appropriate translates of $K[\tha(e)],\; e=a,b$ or $d$  traverses $[1+i,i]$.
First, suppose that $[1+i,i]$ is in  $K[\tha(a)]$ or $K[\tha(d)]+iz$. Then obviously $Q_1$ contains the 1-star graph.

Second, if $[1+i,i]$ is in  $K[\tha(b)]+z$, then the segment $[1+i,i]$ will not be in $Q_1$. However, since  $\tha(a)$ starts with $aba$, $\tha^2(a)$ starts with $\tha(a)\tha(b)$, and thus the segment $[1+i,i]$ will be contained in $Q_2$. In all cases $Q_2$ will contain the 1-star graph. Thus for the string$S*S$  the carrousel $Q_1$ will contain the 1-star graph, and we can copy the proof of the $abc$ case.
\qed

\medskip

Genetically speaking, non-perfectness is recessive, and perfectness is dominant.

\begin{theorem}\label{th:genetics}
Let $S$ and $T$ be  self-avoiding planefilling folding strings. If neither  $S$ nor $T$ is perfect, then $S*T$ is not perfect. If either  $S$ or $T$ is perfect, then $S*T$ is perfect.
\end{theorem}

\noindent {\bf Proof:} We use the characterization of perfectness based on the segment $[1+i,i]$. We write $z_S=f(\tha_S(a))$, $z_T=f(\tha_T(a))$, and $z_{S*T}=z_Sz_T$.

\noindent\emph{Part 1:}  Suppose neither  $S$ nor $T$ is perfect. Then $K[\tha_T(c)]+(1+i)z_T$ traverses $[1+i,i]$, and so part of  the curve $K[\tha_S(\tha_T(c))]+(1+i)z_{S*T}$ will go from the point $(1+i)z_S$ to the point $iz_S$. But this part must be $K[\tha_S(c)]+(1+i)z_S$ in the $\tha_S$-loop, and since $\tha_S$ is not perfect, this part traverses $[1+i,i]$.
Conclusion: $K[\tha_S\tha_T(c)]+iz_{S*T}$ traverses $[1+i,i]$, and thus $S*T$ is not perfect.

\begin{figure}[!b]
\centering
\begin{tikzpicture}[scale=.5,rounded corners=5pt]
\path(1,0)  coordinate (a); \path(0,1)  coordinate (b);
\path(-1,0) coordinate (c); \path(0,-1) coordinate (d);
% The  loop from S:
\draw[blue](-2.5,8)--++(a)--++(b)--++(a)--++(b)--++(c)--++(b)--++(a)--++(d)--++(a)--++(d)--++(a)--++(b)
--++(a)--++(b)--++(c)--++(d)--++(c) --++(b)--++(c)--++(b)--++(c)--++(d)--++(c)--++(d)--++(a) --++(d)--++(c)--++(b)--++(c)--++(b)--++(c)--++(d)--++(c) --++(d)--++(a)--++(b)--++(a)--++(d)--++(a)--++(d)--cycle;
\node[fill=red,circle, inner sep=0.8pt] at (-2.5,8)  {};
\node[fill=black,circle, inner sep=0.8pt] at (0.5,9)  {};
\node[fill=black,circle, inner sep=0.8pt] at (-0.5,12)  {};
\node[fill=black,circle, inner sep=0.8pt] at (-3.5,11)  {};
\node at (-2.7,7.5) {\small $0$};
\node at (0.4,8.5) {\small $z_S$};
\node at (-0.5,12.5) {\small $(1+i)z_S$};
\draw (-4,11.5) node[] {\small $i\,z_S$};
\draw[red,->, thick] (-1.55,8.95)--(-2.45,8.95);
% The  loop from T:
\draw[blue](15,9)--++(a)--++(b)--++(a)--++(b)--++(c)--++(d)--++(c)--++(b)--++(c)--++(b)--++(c)--++(d)
--++(c)--++(d)--++(a)--++(b)--++(a)--++(d)--++(a)--++(d)--cycle;
\node[fill=red,circle, inner sep=0.8pt] at (15,9)  {};
\node[fill=black,circle, inner sep=0.8pt] at (16,11)  {};
\node[fill=black,circle, inner sep=0.8pt] at (14,12)  {};
\node[fill=black,circle, inner sep=0.8pt] at (13,10)   {};
\node at (14.7,8.5) {\small $0$};
\node at (15.9,11.5) {\small $z_T$};
\node at (13.9,12.5) {\small $(1+i)z_T$};
\draw (12.9,9.5) node[] {\small $i\,z_T$};
\draw[red,->, thick] (15.95,9.95)--(15.05,9.95);
%==========================================================
\path(3,1)  coordinate (A); \path(-1,3)  coordinate (B);
\path(-3,-1) coordinate (C); \path(1,-3) coordinate (D); \path(12,4) coordinate (zSc);
%==========================================================
% The  loop from S*T:
\draw[blue](10,0)--++(a)--++(b)--++(a)--++(b)--++(c)--++(b)--++(a)--++(d)--++(a)--++(d)--++(a)--++(b)--++(a)--++(b)--++(c)--++(d)
--++(c)--++(b)--++(c)--++(b)--++(a)--++(b)--++(a)--++(b)--++(c)--++(b)--++(a)--++(d)--++(a)--++(d)--++(a)--++(b)--++(a)--++(b)
--++(c)--++(d)--++(c)--++(b)--++(c)--++(b)--++(c)--++(d)--++(c)--++(d)--++(a)--++(d)--++(c)--++(b)--++(c)--++(b)--++(c)--++(d)
--++(c)--++(d)--++(a)--++(b)--++(a)--++(d)--++(a)--++(d)--++(c)--++(d)--++(c)--++(d)--++(a)--++(d)--++(c)--++(b)--++(c)--++(b)
--++(a)--++(b)--++(a)--++(b)--++(c)--++(d)--++(c)--++(b)--++(c)--++(b)--++(c)--++(d)--++(c)--++(d)--++(a)--++(d)--++(c)--++(b)
--++(c)--++(b)--++(a)--++(b)--++(a)--++(b)--++(c)--++(d)--++(c)--++(b)--++(c)--++(b)--++(c)--++(d)--++(c)--++(d)--++(a)--++(d)
--++(c)--++(b)--++(c)--++(b)--++(c)--++(d)--++(c)--++(d)--++(a)--++(b)--++(a)--++(d)--++(a)--++(d)--++(c)--++(d)--++(c)--++(d)
--++(a)--++(d)--++(c)--++(b)--++(c)--++(b)--++(c)--++(d)--++(c)--++(d)--++(a)--++(b)--++(a)--++(d)--++(a)--++(d)--++(a)--++(b)
--++(a)--++(b)--++(c)--++(b)--++(a)--++(d)--++(a)--++(d)--++(a)--++(b)--++(a)--++(b)--++(c)--++(d)--++(c)--++(b)--++(c)--++(b)
--++(a)--++(b)--++(a)--++(b)--++(c)--++(b)--++(a)--++(d)--++(a)--++(d)--++(c)--++(d)--++(c)--++(d)--++(a)--++(b)--++(a)--++(d)
--++(a)--++(d)--++(a)--++(b)--++(a)--++(b)--++(c)--++(b)--++(a)--++(d)--++(a)--++(d)--++(c)--++(d)--++(c)--++(d)--++(a)--++(b)
--++(a)--++(d)--++(a)--++(d)--cycle;
% The  loop from T scaled:
\draw[black!40,rounded corners=10pt, thick]
(10,0)--++(A)--++(B)--++(A)--++(B)--++(C)--++(D)--++(C)--++(B)--++(C)--++(B)
--++(C)--++(D)--++(C)--++(D)--++(A)--++(B)--++(A)--++(D)--++(A)--++(D);
% Grid points and [1+i,i]
\node[fill=red,circle, inner sep=1.2pt] at (10,0)  {};
\node[fill=black,circle, inner sep=1.2pt] at (11,7)  {};
\node[fill=black,circle, inner sep=1.2pt] at (3,1)  {};
\node[fill=black,circle, inner sep=1.2pt] at (4,8) {};
\node[fill=black,circle, inner sep=0.8pt] at (12,4)  {};
\node[fill=black,circle, inner sep=0.8pt] at (9,3)  {};
\node[fill=black,circle, inner sep=0.8pt] at (13,1)   {};
\node at (9.6,-.4) {$0$}; \node at (10.9,7.6) {$z_{S*T}$};
\node at (4.5,8.6) {$(1+i)z_{S*T}$}; \node at (3.5,0.4) {$i\,z_{S*T}$};
\node at (12.9,0.5) {\small $z_S$};
% wit maken van de achtergrond: rectangle met LO hoek en RB hoek
\draw[color=white,fill=white](7.5,3.23) rectangle(9.7,3.75);
\draw (8.5,3.5) node[] {\small $i\,z_S$};
\draw[color=white,fill=white](11.2,3.23) rectangle(12.9,3.85);
\node at (12.9,3.5) {\small $(1+i)z_S$};
\draw[red,->, thick] (10.95,.95)--(10.05,.95);
\end{tikzpicture}
\vspace*{-1.2cm}\caption{\small The $\tha$-loops of the non-perfect fold $S=DUDDUUUDU$, the perfect fold $T=DUDD$, and the perfect 50-fold $S*T$.}\label{fig:50-fold}
\end{figure}
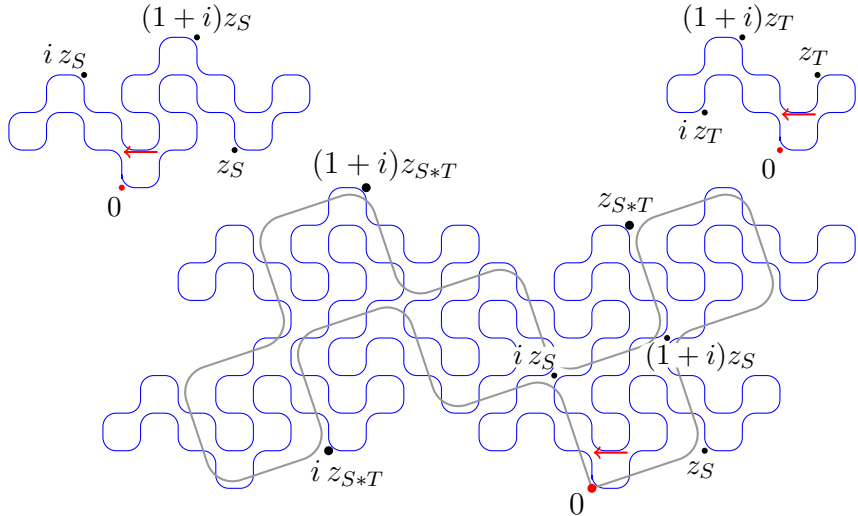

\noindent\emph{Part 2:}
\emph{a)} First suppose $S$ is perfect, $T$ perfect or not. There is the special case that $S=D$, the 2-fold. Then $\tha_S(a)=ab$ and $\tha_S(b)=cb$, so
$$\tha_S\tha_T(a)=\tha_S(ab\dots)=abcb\dots$$
starts with $abc$, and hence $S*T$ is perfect. The following argument is illustrated by Figure~\ref{fig:conv-max-simp}.
If $m_S\ge 3$, then either $\tha_S(a)=abc\dots$, and then $[1+i,i]$ is traversed by $K[\tha_S\tha_T(a)]$, or $\tha_S(a)=aba\dots$, in which case $1+i$ is in the interior of the $\tha_S\tha_T$-loop, but also in the  interior of the $\tha_S$-loop (actually the loop corresponding to the $K$-$\tha_S$-image of the unit square, since the segments will not be traversed in the direct order). By perfectness of $S$, $[1+i,i]$ is traversed by either $K[\tha_S(ab)]$, which is the beginning part of $K[\tha_S\tha_T(a)]$, or by $K[\tha_S(d)]+i\,z_S$ which is the last part of $K[\tha_S\tha_T(d)]+i\,z_{S*T}$. Therefore $[1+i,i]$ can not be traversed by $K[\tha_S\tha_T(c)]+(1+i)z_{S*T}$, and so $S*T$ is perfect.

\emph{b)} This argument is illustrated by Figure~\ref{fig:50-fold}. If $T$ is perfect, then (by Lemma~\ref{lem:segment}) part of  one of the curves $K[\tha_S(\tha_T(a))]$, $K[\tha_S(\tha_T(b))]+z_Sz_T$ or $K[\tha_S(\tha_T(d))]+i\,z_Sz_T$ will go from the point $(1+i)z_S$ to the point $i\,z_S$, and this part must be $K[\tha_S(c)]+(1+i)z_S$ in the $\tha_S$-loop (or rather the $K$-$\tha_S$-image of the unit square). Then if $S$ is \emph{not} perfect, this $K[\tha_S(c)]+(1+i)z_S$ traverses $[1+i,i]$.
So if $T$ is perfect, and $S$ not, then $[1+i,i]$ is not traversed by $K[\tha_S\tha_T(c)]+(1+i)z_{S*T}$, and so $S*T$ is perfect.\qed

\begin{figure}[!b]
\centering
%\begin{tikzpicture}[scale=.4,rounded corners=1pt]
%\draw (0,0) circle (8cm); \node at (0,0) {TOO LONG};
%\end{tikzpicture}
\begin{tikzpicture}[scale=.06,rounded corners=1pt]
\path(1,0)  coordinate (a); \path(0,1)  coordinate (b);
\path(-1,0) coordinate (c); \path(0,-1) coordinate (d);
%========================================================
\path(.5,0)  coordinate (ha); \path(0,.5)  coordinate (hb);
\path(-.5,0) coordinate (hc); \path(0,-.5) coordinate (hd);
%=========================================================
\color{blue!90!white}
\def\tha{\draw (P)++(ha)--++(ha)--++(b)--++(a)--++(b)--++(c)--++(b)--++(a)--++(d)--++(a)--++(d)--++(a)--++(b)--++(a)--++(b)--++(c)--++(d)
 --++(c)--++(b)--++(c)--++(b)--++(a)--++(b)--++(a)--++(b)--++(c)--++(b)--++(a)--++(d)--++(a)--++(d)--++(a)--++(b)--++(a)--++(b) --++(c)--++(d)--++(c)--++(b)--++(c)--++(b)--++(c)--++(d)--++(c)--++(d)--++(a)--++(d)--++(c)--++(b)
--++(c)--++(hb)++(hb) coordinate(P) };
\def\thb{\draw (P)++(hc)--++(hc)--++(d)--++(c)--++(d)--++(a)--++(b)--++(a)--++(d)--++(a)--++(d)--++(c)--++(d)--++(c)--++(d)--++(a)
--++(d)--++(c)--++(b)--++(c)--++(b)--++(a)--++(b)--++(a)--++(b)--++(c)--++(d)--++(c)--++(b)--++(c)--++(b)--++(c)--++(d) --++(c)--++(d)--++(a)--++(d)--++(c)--++(b)--++(c)--++(b)--++(a)--++(b)--++(a) --++(b)--++(c)--++(d)--++(c)--++(b)--++(c)--++(hb)++(hb) coordinate(P)}
\def\thc{\draw (P)++(hc)--++(hc)--++(d)--++(c)--++(d)--++(a)--++(d)--++(c)--++(b)--++(c)--++(b)--++(c)--++(d)--++(c)
--++(d)--++(a)--++(b)--++(a)--++(d)--++(a)--++(d)--++(c)--++(d)--++(c)--++(d)--++(a)--++(d)--++(c)
--++(b)--++(c)--++(b)--++(c)--++(d)--++(c)--++(d)--++(a)--++(b)--++(a)--++(d)--++(a)--++(d)--++(a)
--++(b)--++(a)--++(b)--++(c)--++(b)--++(a)--++(d)--++(a)--++(hd)++(hd) coordinate(P)}
\def\thd{\draw
(P)++(ha)--++(ha)--++(b)--++(a)--++(b)--++(c)--++(d)--++(c)--++(b)--++(c)--++(b)--++(a)--++(b)--++(a)
--++(b)--++(c)--++(b)--++(a)--++(d)--++(a)--++(d)--++(c)--++(d)--++(c)--++(d)--++(a) --++(b)--++(a)--++(d)--++(a)--++(d)--++(a)--++(b)--++(a)--++(b)--++(c)--++(b)--++(a)--++(d) --++(a)--++(d)--++(c)--++(d)--++(c)--++(d)--++(a)--++(b)--++(a)--++(d)--++(a)--++(hd)++(hd) coordinate(P)}
\def\thab{\draw (Q)++(a)++(b)++(a)++(b)++(c)++(b)++(a)++(d)++(a)++(d)++(a)++(b)++(a)++(b)++(c)++(d)
 ++(c)++(b)++(c)++(b)++(a)++(b)++(a)++(b)++(c)++(b)++(a)++(d)++(a)++(d)++(a)++(b)++(a)++(b) ++(c)++(d)++(c)++(b)++(c)++(b)++(c)++(d)++(c)++(d)++(a)++(d)++(c)++(b)
++(c)++(hb)--++(hb) coordinate(Q) --++(hc)};
\def\thad{\draw (Q)++(a)++(b)++(a)++(b)++(c)++(b)++(a)++(d)++(a)++(d)++(a)++(b)++(a)++(b)++(c)++(d)
 ++(c)++(b)++(c)++(b)++(a)++(b)++(a)++(b)++(c)++(b)++(a)++(d)++(a)++(d)++(a)++(b)++(a)++(b) ++(c)++(d)++(c)++(b)++(c)++(b)++(c)++(d)++(c)++(d)++(a)++(d)++(c)++(b)
++(c)++(hb)--++(hb) coordinate(Q) --++(ha)};
\def\thbc{\draw (Q)++(c)++(d)++(c)++(d)++(a)++(b)++(a)++(d)++(a)++(d)++(c)++(d)++(c)++(d)++(a)
++(d)++(c)++(b)++(c)++(b)++(a)++(b)++(a)++(b)++(c)++(d)++(c)++(b)++(c)++(b)++(c)++(d) ++(c)++(d)++(a)++(d)++(c)++(b)++(c)++(b)++(a)++(b)++(a) ++(b)++(c)++(d)++(c)++(b)++(c)++(hb)--++(hb) coordinate(Q)--++(hc)};
\def\thbx{\draw (Q)++(c)++(d)++(c)++(d)++(a)++(b)++(a)++(d)++(a)++(d)++(c)++(d)++(c)++(d)++(a)
++(d)++(c)++(b)++(c)++(b)++(a)++(b)++(a)++(b)++(c)++(d)++(c)++(b)++(c)++(b)++(c)++(d) ++(c)++(d)++(a)++(d)++(c)++(b)++(c)++(b)++(a)++(b)++(a) ++(b)++(c)++(d)++(c)++(b)++(c)--++(hb)--++(hb)};
\def\thba{\draw (Q)++(c)++(d)++(c)++(d)++(a)++(b)++(a)++(d)++(a)++(d)++(c)++(d)++(c)++(d)++(a)
++(d)++(c)++(b)++(c)++(b)++(a)++(b)++(a)++(b)++(c)++(d)++(c)++(b)++(c)++(b)++(c)++(d) ++(c)++(d)++(a)++(d)++(c)++(b)++(c)++(b)++(a)++(b)++(a) ++(b)++(c)++(d)++(c)++(b)++(c)++(hb)--++(hb) coordinate(Q)--++(ha)};
\def\thcd{\draw (Q)++(c)++(d)++(c)++(d)++(a)++(d)++(c)++(b)++(c)++(b)++(c)++(d)++(c)
++(d)++(a)++(b)++(a)++(d)++(a)++(d)++(c)++(d)++(c)++(d)++(a)++(d)++(c)
++(b)++(c)++(b)++(c)++(d)++(c)++(d)++(a)++(b)++(a)++(d)++(a)++(d)++(a)
++(b)++(a)++(b)++(c)++(b)++(a)++(d)++(a)++(hd)--++(hd) coordinate(Q) --++(ha)};
\def\thcb{\draw (Q)++(c)++(d)++(c)++(d)++(a)++(d)++(c)++(b)++(c)++(b)++(c)++(d)++(c)
++(d)++(a)++(b)++(a)++(d)++(a)++(d)++(c)++(d)++(c)++(d)++(a)++(d)++(c)
++(b)++(c)++(b)++(c)++(d)++(c)++(d)++(a)++(b)++(a)++(d)++(a)++(d)++(a)
++(b)++(a)++(b)++(c)++(b)++(a)++(d)++(a)++(hd)--++(hd) coordinate(Q) --++(hc)};
\def\thcx{\draw (Q)++(c)++(d)++(c)++(d)++(a)++(d)++(c)++(b)++(c)++(b)++(c)++(d)++(c)++(d)++(a)++(b)++(a)++(d)++(a)
++(d)++(c)++(d)++(c)++(d)++(a)++(d)++(c)++(b)++(c)++(b)++(c)++(d)++(c)++(d)++(a)++(b)++(a)++(d)
++(a)++(d)++(a)++(b)++(a)++(b)++(c)++(b)++(a)++(d)++(a)++(hd)++(hd)};
\def\thda{\draw
(Q)++(a)++(b)++(a)++(b)++(c)++(d)++(c)++(b)++(c)++(b)++(a)++(b)++(a)++(b)++(c)++(b)++(a)++(d)++(a)++(d)++(c)
++(d)++(c)++(d)++(a) ++(b)++(a)++(d)++(a)++(d)++(a)++(b)++(a)++(b)++(c)++(b)++(a)++(d) ++(a)++(d)++(c)++(d)++(c)++(d)++(a)++(b)++(a)++(d)++(a)++(hd)--++(hd) coordinate(Q) --++(ha)};
\def\thdc{\draw
(Q)++(a)++(b)++(a)++(b)++(c)++(d)++(c)++(b)++(c)++(b)++(a)++(b)++(a)
++(b)++(c)++(b)++(a)++(d)++(a)++(d)++(c)++(d)++(c)++(d)++(a) ++(b)++(a)++(d)++(a)++(d)++(a)++(b)++(a)++(b)++(c)++(b)++(a)++(d) ++(a)++(d)++(c)++(d)++(c)++(d)++(a)++(b)++(a)++(d)++(a)++(hd)--++(hd) coordinate(Q) --++(hc)};
\def\thaa{\tha;\thb;\tha;\thb;\thc;\thb;\tha;\thd;\tha;\thd;\tha;\thb;\tha;\thb;\thc;\thd;
 \thc;\thb;\thc;\thb;\tha;\thb;\tha;\thb;\thc;\thb;\tha;\thd;\tha;\thd;\tha;\thb;\tha;\thb; \thc;\thd;\thc;\thb;\thc;\thb;\thc;\thd;\thc;\thd;\tha;\thd;\thc;\thb;\thc;\thb}
\def\thaabb{\thab;\thba;\thab;\thbc;\thcb;\thba;\thad;\thda;\thad;\thda;\thab;\thba;\thab;\thbc;\thcd;
\thdc;\thcb;\thbc;\thcb;\thba;\thab;\thba;\thab;\thbc;\thcb;\thba;\thad;\thda;\thad;\thda;\thab;\thba;\thab; \thbc;\thcd;\thdc;\thcb;\thbc;\thcb;\thbc;\thcd;\thdc;\thcd;\thda;\thad;\thdc;\thcb;\thbc;\thcb;\thbx;};
% initialiseren
\path(0,0)  coordinate (P); \path(0,0)  coordinate (Q); \draw (0,0)--(ha); \thaa; \thaabb;
\color{red!90!white}
\path(1,0)  coordinate (b); \path(0,1)  coordinate (c);
\path(-1,0) coordinate (d); \path(0,-1) coordinate (a);
%=========================================================
\path(.5,0)  coordinate (hb); \path(0,.5)  coordinate (hc);
\path(-.5,0) coordinate (hd); \path(0,-.5) coordinate (ha);
%=========================================================
\path(0,0)  coordinate (P);  \path(0,0)  coordinate (Q); \draw (0,0)--(ha); \thaa; \thaabb;
\color{green!90!white}
\path(1,0)  coordinate (c); \path(0,1)  coordinate (d);
\path(-1,0) coordinate (a); \path(0,-1) coordinate (b);
%=========================================================
\path(.5,0)  coordinate (hc); \path(0,.5)  coordinate (hd);
\path(-.5,0) coordinate (ha); \path(0,-.5) coordinate (hb);
%=========================================================
\path(0,0)  coordinate (P);  \path(0,0)  coordinate (Q); \draw (0,0)--(ha); \thaa; \thaabb;
\color{black!90!white}
\path(1,0)  coordinate (d); \path(0,1)  coordinate (a);
\path(-1,0) coordinate (b); \path(0,-1) coordinate (c);
%=========================================================
\path(.5,0)  coordinate (hd); \path(0,.5)  coordinate (ha);
\path(-.5,0) coordinate (hb); \path(0,-.5) coordinate (hc);
%=========================================================
\path(0,0)  coordinate (P);  \path(0,0)  coordinate (Q); \draw (0,0)--(ha); \thaa; \thaabb;
\end{tikzpicture}
\vspace*{-.2cm}\caption{\small The carrousel $Q_2$ of the 50-fold $S*T$ of Figure~\ref{fig:50-fold}.}\label{fig:50-fold-Q2}
\end{figure}
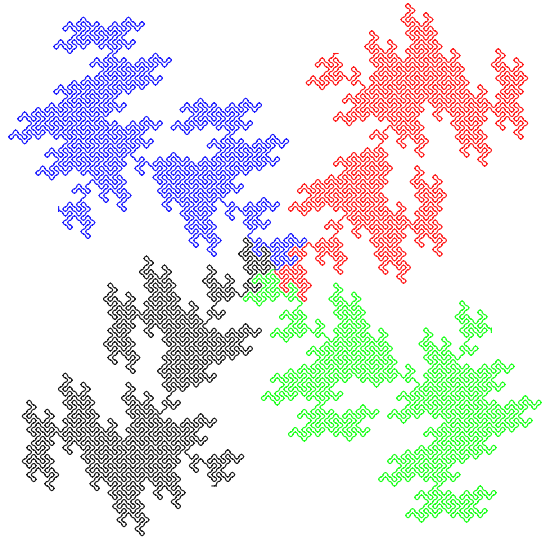

The simplest example  of a non-perfect fold, already given by Davis and Knuth (\cite{DK}), is  the 4-fold $S=DUU$. Here  $K[\tha^n(a)]$ fills a $45^o-45^o-90^o$ triangle for all $n$, and the carrousel fills only half of the plane.

Remarkably, the non-perfect 10-fold from Figure~\ref{fig:50-fold}  can be written as $S=DUDDUUUDU=DUDD*U$, and $DUDD$ and $U$ are both perfect folds. This does not contradict Theorem~\ref{th:genetics}, since we always require $S$ and $T$ to start with $D$.
The $Q_2$ of this perfect 50-fold is depicted in Figure~\ref{fig:50-fold-Q2}.

\section{Boundaries}

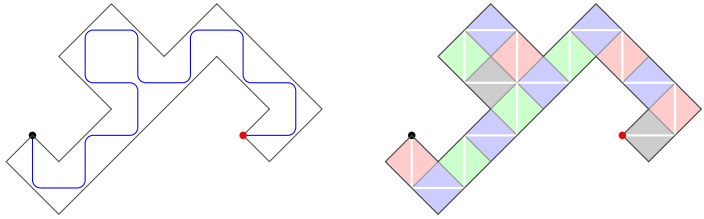
\begin{figure}[h!]
\hspace*{-1cm}
\centering
\begin{tikzpicture}[scale=.7,rounded corners=3pt]
\path(1,0)  coordinate (a);
\path(0,1)  coordinate (b);
\path(-1,0) coordinate (c);
\path(0,-1) coordinate (d);
% draw thna:
\draw[blue] (0,0)--++(a)--++(b)--++(c)--++(b)--++(c)--++(d)--++(c)--++(b)--++(c)--++(d)--++(a)--++(d)--++(c)--++(d)--++
(c)--++(b);
% Drawing the red points
\node[fill=red,circle, inner sep=1pt] at (0,0)  {};
\node[fill=black,circle, inner sep=1pt] at (-4,0)  {};
\path(0.5,0.5)  coordinate (NE);
\path(0.5,-0.5)  coordinate (SE);
\path(-0.5,0.5) coordinate (NW);
\path(-0.5,-0.5) coordinate (SW);
\draw[black!70, sharp corners] (0,0)--++(NE)--++(NW)--++(NW)--++(SW)--++(SW)--++(SW)--++(SW)--++(SW)--++(SW)--++(NW)
--++(NW)--++(NE)--++(SE)--++(NE)--++(NE)--++(NW)--++(NW)--++(NE)--++(NE)--++(SE)--++(SE)--++(NE)--++(NE)
--++(SE)--++(SE)--++(SE)--++(SE)--++(SW)--++(SW)--++(NW);
\end{tikzpicture}\qquad
\begin{tikzpicture}[scale=.7,rounded corners=3pt]
\path(1,0)  coordinate (a);
\path(0,1)  coordinate (b);
\path(-1,0) coordinate (c);
\path(0,-1) coordinate (d);
\def\Diama((#1,#2)){\draw[fill=black!20,draw=black!40,rounded corners=0pt](#1,#2)--++(1/2,1/2)--++(1/2,-1/2)--++(-1/2,-1/2)--++(-1/2,1/2); \draw[white, thick](#1,#2)--++(a)}
\def\Diamb((#1,#2)){\draw[fill=red!20,draw=black!40,rounded corners=0pt](#1,#2)--++(1/2,1/2)--++(-1/2,1/2)--++(-1/2,-1/2)--++(1/2,-1/2); \draw[white, thick](#1,#2)--++(b)}
\def\Diamc((#1,#2)){\draw[fill=blue!20,draw=black!40,rounded corners=0pt](#1,#2)--++(-1/2,1/2)--++(-1/2,-1/2)--++(1/2,-1/2)--++(1/2,1/2); \draw[white, thick](#1,#2)--++(c)}
\def\Diamd((#1,#2)){\draw[fill=green!20,draw=black!40,rounded corners=0pt](#1,#2)--++(-1/2,-1/2)--++(1/2,-1/2)--++(1/2,1/2)--++(-1/2,1/2); \draw[white, thick](#1,#2)--++(d)}
% Diamonds of the 16-fold
\Diama((0,0));\Diamb((1,0));\Diamc((1,1));\Diamb((0,1));
\Diamc((0,2));\Diamd((-1,2));\Diamc((-1,1));\Diamb((-2,1));\Diamc((-2,2));\Diamd((-3,2));
\Diama((-3,1));\Diamd((-2,1));\Diamc((-2,0));\Diamd((-3,0));\Diamc((-3,-1));\Diamb((-4,-1));
%% Drawing the end points
\node[fill=red,circle, inner sep=1pt] at (0,0)  {};
\node[fill=black,circle, inner sep=1pt] at (-4,0)  {};
\path(0.5,0.5)  coordinate (NE);
\path(0.5,-0.5)  coordinate (SE);
\path(-0.5,0.5) coordinate (NW);
\path(-0.5,-0.5) coordinate (SW);
\draw[black!70, sharp corners] (0,0)--++(NE)--++(NW)--++(NW)--++(SW)--++(SW)--++(SW)--++(SW)--++(SW)--++(SW)--++(NW)
--++(NW)--++(NE)--++(SE)--++(NE)--++(NE)--++(NW)--++(NW)--++(NE)--++(NE)--++(SE)--++(SE)--++(NE)--++(NE)
--++(SE)--++(SE)--++(SE)--++(SE)--++(SW)--++(SW)--++(NW);
\end{tikzpicture}
\caption{Left: Mid grid  border of  the curve $\rK[abcbcdcbcdadcdcb]$. Right: the polyomino $K_\diamond[abcbcdcbcdadcdcb]$. }\label{fig:midgrid}
\end{figure}

Let $w$ be a $\DD$-word. The curve $K[w]$ starts in 0, and ends in $f(w)$. The smallest simple closed polygonal  curve that contains $K[w]$ and has vertices 0, $f(w)$ and further only mid grid points will be called the \emph{mid grid border} of the curve $K[w]$, see Figure \ref{fig:midgrid}, Left. We have already encountered such borders: the octagon $\Omega$ from Section \ref{sec:exist} is the mid grid border of the curve $K[\tha(a)]$ generated by any of the morphisms $\tha$ that are constructed in the existence proof.

There is an alternative and attractive way to define this border, given by Knuth in the Addendum of \cite{Fun} (see also the recent \cite{K2010}): the mid grid border of $K[w]$ is the topological border of the polyomino generated by $w$. In the context of Section \ref{sec:rec} this polyomino is obtained via the diamond map  $K_\diamond[\cdot]$ mapping words to compact subsets of $\C$, given by
$$K_\diamond[a]=B(\tfrac12,\tfrac12),\; K_\diamond[b]=B(\tfrac{i}2,\tfrac12),\; K_\diamond[c]=B(-\tfrac{1}2,\tfrac12),\; K_\diamond[d]=B(-\tfrac{i}2,\tfrac12),
$$
where $B(z,r)$ is the ball with center $x$ and radius $r$ in the $L_1$-metric. Equation (\ref{eq:Ktha}) then defines the polyomino $K_\diamond[w]$ associated to any word $w$.

Let $S$ be a folding string with associated morphism $\tha$. Gaussian integers that are visited twice by the curve $K[\tha(a)]$ will be called \emph{internal points}, and those visited once \emph{boundary points}. Obviously, if $S$ codes an $m$-fold, and we denote by $i(S)$ and $r(S)$ the number of internal and boundary points of $S$, then
$$m=2\,i(S)+r(S)-1.$$

\begin{proposition}\label{prop:bound}
Let $r(S)$ be the number of Gaussian integer boundary points of the curve associated to a fold $S$. Let $S$ and $T$ be the strings of two self-avoiding folds. Suppose that $S$ is planefilling  and symmetric ($S=\overline{S}$). Then $$r(S*T)=\frac12\,r(S)r(T).$$
\end{proposition}

\begin{figure}[!h]
\centering
\begin{tikzpicture}[scale=.4,rounded corners=5pt]
\path(1,0)  coordinate (a); \path(0,1)  coordinate (b);
\path(-1,0) coordinate (c); \path(0,-1) coordinate (d);
%==========================================================
\def\octagon1((#1,#2)){\draw[fill=black!4,draw=black!70,rounded corners=0pt]
(#1,#2)--++(1.5,-1.5)--++(-2,-2)--++(2,-2)--++(1.5,1.5)--++(-1.5,1.5)--++(2,2)--++(-2,2)--++(-1.5,-1.5);
\node[fill=black,circle, inner sep=1.5pt] at (#1,#2)  {};}
\octagon1((0,20)); \octagon1((-4,17)); \octagon1((7,19)); \octagon1((3,16));
\def\octagon2((#1,#2)){\draw[fill=black!1,draw=black!70,rounded corners=0pt]
(#1,#2)--++(1.5,1.5)--++(2,-2)--++(2,2)--++(-1.5,1.5)--++(-1.5,-1.5)--++(-2,2)--++(-2,-2)--++(1.5,-1.5);
\node[fill=black,circle, inner sep=1.5pt] at (#1,#2)  {};}
\octagon2((-1,13)); \octagon2((3,16)); \octagon2((-8,14)); \octagon2((6,12));
% The curve:
\draw[blue](0,20)--++(a)--++(b)--++(a)--++(d)--++(a)--++(d)--++(c)--++(b)--++(c)--++(d)--++(a)--++(d)--++(c)--++(d)--++(a)--++(d)--++(c)--++(b)--++(c)--++(d)--++(a)--++(d)--++(a)
--++(b)--++(a)--++(b)--++(c)--++(b)--++(a)--++(b)--++(a)--++(d)--++(c)--++(d)--++(a)--++(b)--++(a)--++(d)--++(a)--++(b)--++(a)--++(d)--++(c)--++(d)--++(a)--++(b)
--++(a)--++(b)--++(c)--++(b)--++(a)--++(b)--++(a)--++(d)--++(a)--++(d)--++(c)--++(b)--++(c)--++(d)--++(a)--++(d)--++(c)--++(d)--++(a)--++(d)--++(c)--++(b)--++(c)
--++(d)--++(a)--++(d)--++(a)--++(b)--++(a)--++(d)--++(a)--++(d)--++(c)--++(d)--++(c)--++(b)--++(a)--++(b)--++(c)--++(d)--++(c)--++(b)--++(c)--++(d)--++(c)--++(b)
--++(a)--++(b)--++(c)--++(d)--++(c)--++(d)--++(a)--++(d)--++(c)--++(d)--++(c)--++(b)--++(c)--++(b)--++(a)--++(d)--++(a)--++(b)--++(c)--++(b)--++(a)--++(b)--++(c)
--++(b)--++(a)--++(d)--++(a)--++(b)--++(c)--++(b)--++(c)--++(d)--++(c)--++(d)--++(a)--++(d)--++(c)--++(d)--++(c)--++(b)--++(a)--++(b)--++(c)--++(d)--++(c)--++(b)
--++(c)--++(d)--++(c)--++(b)--++(a)--++(b)--++(c)--++(d)--++(c)--++(d)--++(a)--++(d)--++(c)--++(d)--++(c)--++(b)--++(c)--++(b)--++(a)--++(d)--++(a)--++(b)--++(c)
--++(b)--++(a)--++(b)--++(c)--++(b)--++(a)--++(d)--++(a)--++(b)--++(c)--++(b)--++(c)--++(d)--++(c)--++(d)--++(a)--++(d)--++(c)--++(d)--++(c)--++(b)--++(a)--++(b)
--++(c)--++(d)--++(c)--++(b)--++(c)--++(d)--++(c)--++(b)--++(a)--++(b) --++(c)--++(d)--++(c)--++(d)--++(a)--++(d);
\node[fill=red,circle, inner sep=1.5pt] at (0,20)  {};
\node[fill=black,circle, inner sep=1.5pt] at (-4,17) {};
\node[fill=black,circle, inner sep=1.5pt] at (10,15) {};
\node[fill=black,circle, inner sep=1.5pt] at (-8,14) {};
\node[fill=black,circle, inner sep=1.5pt] at (7,19) {};
% The quarterboundary points
\node [circle,draw,double,fill=orange!80, inner sep=1pt] at (4,17) {};
\node [circle,draw,double,fill=orange!80, inner sep=1pt] at (5,17) {};
\node [circle,draw,double,fill=orange!80, inner sep=1pt] at (6,16) {};
\node [circle,draw,double,fill=green!80, inner sep=1pt] at (5,13) {};
\node [circle,draw,double,fill=green!80, inner sep=1pt] at (5,14) {};
\node [circle,draw,double,fill=green!80, inner sep=1pt] at (6,15) {};
\node [circle,draw,double,fill=brown!80, inner sep=1pt] at (3,12) {};
\node [circle,draw,double,fill=brown!80, inner sep=1pt] at (4,11) {};
\node [circle,draw,double,fill=brown!80, inner sep=1pt] at (5,11) {};
\node [circle,draw,double,fill=yellow!80, inner sep=1pt] at (4,15) {};
\node [circle,draw,double,fill=yellow!80, inner sep=1pt] at (4,14) {};
\node [circle,draw,double,fill=yellow!80, inner sep=1pt] at (3,13) {};
\end{tikzpicture}
\vspace*{-4.5cm}
\caption{\small The quarter-borders of the symmetric self-avoiding planefilling 25-fold $\tha_S(a)=abadadcbcdadcdadcbcdadaba$ in  $\rK[\tha_S(c)]+2z_S$ in $\rK[\tha_S\tha_T(a)]$, where $T$ is the 8-fold with $\tha_T(a)=abadcdcd$. In this example $k_S=3$.}\label{fig:quarterboundary}
\end{figure}
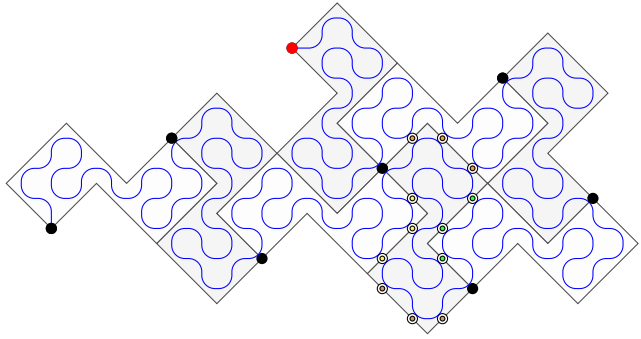

\noindent {\bf Proof:}
Let $\tha_S$ and $\tha_T$ of length $m_S$ and $m_T$ be  the folding morphisms of $S$ and $T$. It can be deduced from the symmetry that the  mid grid border of $K[\tha_S(e)],\; e=a,b,c,d$ consists of four congruent quarter-borders, each bordering say, $k_S$ boundary points (cf.~Figure~\ref{fig:quarterboundary}). Since $S$ is planefilling, there are no other boundary points, except for the first and the last point of $K[\tha_S(e)]$. Hence $r(S)=4k_S+2$. The curve $K[\tha_S\tha_T(a)]$ generated by $S*T$, is obtained by replacing the $m_T$ scaled segments of $L_S(K[\tha_T(a)])$ by translations of curves $K[\tha_S(e)],\; e=a,b,c$ or $d$, where $L_S$ is multiplication by $z_S=f(\tha_S(a))$.

In this process we pass $i(T)$ (scaled) internal points twice, where we lose  eight times
$k_S$ points, and four endpoints. Further we pass $r(T)$ (scaled) boundary points once, where we lose (except at the begin and endpoint of $L_S(K[\tha_T(a)])$) one endpoint and two times $k_S$ points. We thus obtain that
\begin{eqnarray*}
r(S*T)&=&m_T\,r(S)-i(T)(8k_S+4)-(r(T)-2)(2k_S+1)\\
&=&(2i(T)+r(T)-1)r(S)-2i(T)r(S)-(r(T)-2)r(S)/2\\
&=&r(S)r(T)/2.\hspace*{8cm}\qed
\end{eqnarray*}

A slightly surprising consequence of Proposition~\ref{prop:bound} is that if $S$ is any symmetric planefilling string with \emph{no} internal points, then the fold $S*S$ has $(m-1)^2/4$ internal points.

\section{Tiles}\label{sec:tiles}

 By Theorem (2.4) of \cite{Rec}  the sequence of sets $L_\tha^{-n}(K[\tha^n(abcd)])$ converges to a compact set ${\mathcal{T}}_S$, which we call the \emph{tile} generated by $S$. Now if $S$ is a self-avoiding and planefilling $m$-fold, then  $\rK[\tha(abcd)]$ is maximally simple, and this implies that there exists a set of $m$ Gaussian integers $\mathfrak{D}$  such that
 $$K[\tha(abcd)]=\bigcup_{\mathfrak{d}\in\mathfrak{D}} \big(K[abcd]+\mathfrak{d}\big).$$
 But then also for all $n=1,2,\dots$
 $$K[\tha^{n+1}(abcd)]=\bigcup_{\mathfrak{d}\in\mathfrak{D}} \big(K[\tha^n(abcd)]+L_\tha^n(\mathfrak{d})\big),$$
 and so for all $n=1,2,\dots$
  $$L_\tha L_\tha^{-n-1}K[\tha^{n+1}(abcd)]=\bigcup_{\mathfrak{d}\in\mathfrak{D}} \big( L_\tha^{-n}K[\tha^n(abcd)]+\mathfrak{d}\big).$$
  Taking the limit as $n\rightarrow\infty$ we obtain the \emph{self-similarity equation}
 \begin{equation}\label{eq:IFS}L_\tha({\mathcal{T}}_S)=\bigcup_{\mathfrak{d}\in\mathfrak{D}} \big({\mathcal{T}}_S+\mathfrak{d}\big).
 \end{equation}
Applying $L_\tha^{-1}$ on both sides we see that ${\mathcal{T}}_S$ is the unique attractor of a contracting iterated function system\footnote{This connection between recurrent sets as in \cite{Rec},  and language directed IFS's, also called recurrent IFS's, or  graph directed IFS's has been known for a long time, see \cite{Morcrette}.}.

A set ${\mathcal{T}}$ satisfying (\ref{eq:IFS}) is  called a \emph{self-affine} tile . Such tiles have been extensively studied, mostly under the condition that $L_\tha(\mathbb{G})\subset \mathbb{G}$ and that the digit set $\mathfrak{D}$ is a complete residue system for $\mathbb{G}/L_\tha(\mathbb{G})$, where we write $\mathbb{G}$ for the Gausssian integers. In this case $\mathfrak{D}$ is called a \emph{standard} digit set. The authors of  \cite{DKV} and \cite{SW} determine under the standard digit set condition  the Hausdorff dimension of the boundary $\partial {\mathcal{T}}_S$ of ${\mathcal{T}}_S$ , with a rather involved algorithm in both papers (see \cite{DW} for an even more involved algorithm that does not require the standard digit set assumption).

For our folding tiles ${\mathcal{T}}_S$ it is obvious that they will not have a standard digit set when $m$ is even (all digits will be even Gaussian integers). For the case $m=2$ we will still obtain a tiling of the plane by a scaling, translation and rotation of the well known twin-dragon tile (it has been observed by Bandt (\cite{Bandt}, p.550) that this holds for any  2 element digit set $\{0,\mathfrak{d}\}$).

For symmetric  self-avoiding planefilling folds $S$ we present a particularly simple formula for $\dim_{\rm H}(\partial {\mathcal{T}}_S)$, where $\dim_{\rm H}(\cdot)$ denotes the Hausdorff dimension of a set. (Actually \cite{DKV} and \cite{SW} do apply in this case, as it is easy to show that for symmetric $S$ the digit set will be standard.)

\begin{theorem}\label{hdim}
Let  $S$ be a self-avoiding planefilling symmetric $m$-fold with r boundary points. Let ${\mathcal{T}}_S$ be the tile generated by $S$ with topological boundary $\partial {\mathcal{T}}_S$. Then
$$\dim_{\rm H}(\partial {\mathcal{T}}_S)=\frac{2\log(r)-\log(4)}{\log(m)}.$$
\end{theorem}

\noindent {\bf Proof:}
The tile ${\mathcal{T}}_S$ is a union of 4 sets: the limits of $L_\tha^{-n}(K[\tha^n(e)])+f(e)$, where  $e=a,b,c,d.$
By Lemma~\ref{lem:dir}, the boundaries of these sets all have the same Hausdorff dimension. So it follows
with the usual counting and scaling argument (cf.~\cite{Falc}), that
$$\dim_{\rm H}(\partial {\mathcal{T}}_S)=\lim_{n\rightarrow\infty}\frac{\log r(S^{*n})}{\log(|z|^n)}.$$
where $z$ is the endpoint of $K[\tha(a)].$ Since $S$ is self-avoiding and planefilling, we have $|z|=\sqrt{m}$, and
from Proposition~\ref{prop:bound} we obtain
$$r(S^{*n})=\Big(\frac12\Big)^{n-1}\,r^n=2\Big(\frac{r}2\Big)^{n}.$$
This leads to the formula in the theorem.\qed

\medskip

We give an example: for $m=5$ we have the symmetric fold $S=DDUU$ with $\tha(a)=abcba$. Since there are no internal points,  the boundary of the tile $\mathcal{T}_S$ has Hausdorff dimension $$\frac{2\log(6)-\log(4)}{\log(5)}=\frac{\log(9)}{\log(5)}.$$

\section{Other angles, more morphisms}\label{sec:final}

The strip of paper can be opened at other angles than $90^o$. In their paper (\cite{DK}) Davis and Knuth observe that the fold $S=DU$, when opened at $60^o$
angles yields an interesting curve, which they coined the \emph{terdragon}. They prove, using a ternary number system, that 6 infinite terdragons rotated over successive multiples of $60^o$ pass exactly once through every edge of the triangular lattice. Actually, the theory developed in this paper can give a completely different proof of this property.
We should now use a 3 letter alphabet, say $\{a,c,e\}$. The letters correspond to segments in the complex plane by the map $f$, now defined by
$$f(a)=1, \quad f(c)=\omega, \quad f(e)=\omega^2, \quad {\rm with }\; \omega:=\me^{2\pi i/3}.$$
The anti-morphism $\tau$ is as before, and the morphism $\sigma$ is now defined by $\sigma(a)=c, \sigma(c)=e, \sigma(e)=a.$
Folding morphisms again are those that commute with $\sigma\tau$, and, moreover, do not contain $aa$, $cc$ or $ee$. The basic example is $\tha$ with $\tha(a)=aca$, which generates the terdragon.

Defining the $\tha$-loop by the curve $\rK[\tha(ace)]$, and the anti-$\tha$-loop by the curve $\rK[\tha(aec)]$, then Theorem \ref{th:avoid} remains true \emph{provided} we replace `the $\tha$-loop' by `the $\tha$-loop \emph{and} the anti-$\tha$-loop'. The difference is caused by the different symmetries of the triangular lattice.
%%%(Example that necessary??).

\begin{figure}[!h]
\centering
\begin{tikzpicture}[scale=.2,rounded corners=3.5pt]
\path(2,0)  coordinate (a);
\path(1,{sqrt(3)})  coordinate (b);
\path(-1,{sqrt(3)}) coordinate (c);
\path(-2,0) coordinate (d);
\path(-1,{-sqrt(3)}) coordinate (e);
\path(1,{-sqrt(3)}) coordinate (f);
% Midgrid boundary:
\path(1,{sqrt(3)/3})  coordinate (m1); \path(0,{2*sqrt(3)/3})  coordinate (m2);
\path(-1,{sqrt(3)/3})   coordinate (m3); \path(1,{-sqrt(3)/3}) coordinate (m6);
\def\midgrid{\draw[black!70, rounded corners=0pt] (0,0)--++(m1)--++(m2)--++(m1)--++(m2)--++(m1)--++(m2)--++(m1)--++(m2)--++(m1)--++(m2)
--++(m1)--++(m6)--++(m1)--++(m6)--++(m1)--++(m6)--++(m1)--++(m6)--++(m1)--++(m6)
(0,0)--++(m6)--++(m1)--++(m6)--++(m1)--++(m6)--++(m1)--++(m6)--++(m1)--++(m6)--++(m1)
--++(m2)--++(m1)--++(m2)--++(m1)--++(m2)--++(m1)--++(m2)--++(m1)--++(m2)--++(m1)};
\midgrid; \draw[black!70, rounded corners=0pt] (0,{10*{sqrt(3)}})--++(m6)--++(m1)--++(m6)--++(m1)--++(m6)--++(m1)--++(m6)--++(m1)--++(m6)--++(m1)
(15,{5*{sqrt(3)}})--++(m2)--++(m3)--++(m2)--++(m3)--++(m2)--++(m3)--++(m2)--++(m3)--++(m2)--++(m3);
\path(0,{2*sqrt(3)/3}) coordinate (m1);\path(-1,{sqrt(3)/3}) coordinate (m2);\path(1,{sqrt(3)/3}) coordinate (m6);
\midgrid;

\draw[blue,thick] (0,0)--++
(a)--++(c)--++(a)--++(e)--++(a)--++(c)--++(a)--++(e)--++(a)--++(c)--++(a)--++(e)--++(a)--++(c)--++(a)--++(e)--++(a)--++(c)--++(a)--++(c)--++(e)--++(c)--++(e)--++(c)--++(e)--++(c)--++(e)--++ (c)--++(a)--++(c)--++(a)--++(e)--++(a)--++(c)--++(a)--++(e)--++(a)--++(c)--++(a)--++(e)--++(a)--++(c)--++(a)--++(e)--++(a)--++(c)--++(a)--++(c)--++(e)--++(c)--++(e)--++(c)--++(e)--++(c)--++(e)--++ (c)--++(a)--++(c)--++(a)--++(e)--++(a)--++(c)--++(a)--++(e)--++(a)--++(c)--++(a)--++(e)--++(a)--++(c)--++(a)--++(e)--++(a)--++(c)--++(a)--++(c)--++(e)--++(c)--++(a)--++(c)--++(e)--++ (c)--++(a)--++(c)--++(e)--++(c)--++(a)--++ (c)--++(e)--++(c)--++(a)--++(c)--++(e)--++(c)--++(e)--++(a)--++(e)--++(a)--++(e)--++(a)--++(e)--++(a)--++(e)--++(c)--++(e)--++(c)--++(a)--++(c)--++(e)--++ (c)--++(a)--++(c)--++(e)--++ (c)--++(a)--++(c)--++(e)--++(c)--++(a)--++(c)--++(e)--++(c)--++(e)--++(a)--++(e)--++(a)--++(e)--++(a)--++(e)--++(a)--++(e)--++(c)--++(e)--++(c)--++(a)--++(c)--++(e)--++(c)--++(a)--++(c)--++(e)--++ (c)--++(a)--++(c)--++(e)--++(c)--++(a)--++(c)--++(e)--++(c)--++(e)--++(a)--++(e)--++(c)--++(e)--++(a)--++(e)--++(c)--++(e)--++(a)--++(e)--++(c)--++(e)--++(a)--++(e)--++(c)--++(e)--++(a)--++(e)--++ (a)--++(c)--++(a)--++(c)--++(a)--++(c)--++(a)--++(c)--++(a)--++(e)--++(a)--++(e)--++(c)--++(e)--++(a)--++(e)--++(c)--++(e)--++(a)--++(e)--++(c)--++(e)--++(a)--++(e)--++(c)--++(e)--++(a)--++(e)--++ (a)--++(c)--++(a)--++(c)--++(a)--++(c)--++ (a)--++(c)--++(a)--++(e)--++(a)--++(e)--++(c)--++(e)--++(a)--++(e)--++(c)--++(e)--++(a)--++(e)--++(c)--++(e)--++(a)--++(e)--++(c)--++(e)--++(a)--++(e);
% De oorsprong:
\node[fill=red,circle, inner sep=1.4pt] at (0,0)  {};
% Eindpunt van tha(a):
\node[fill=black,circle, inner sep=1.0pt] at (15,{5*{sqrt(3)}})  {};
% Eindpunt van tha(ac):
\node[fill=black,circle, inner sep=1.0pt] at (0,{10*{sqrt(3)}})  {};
\end{tikzpicture}\hspace*{2cm}
\begin{tikzpicture}[scale=.2,rounded corners=3.5pt]
\path(2,0)  coordinate (a); \path(1,{sqrt(3)})  coordinate (b);
\path(-1,{sqrt(3)}) coordinate (c); \path(-2,0) coordinate (d);
\path(-1,{-sqrt(3)}) coordinate (e); \path(1,{-sqrt(3)}) coordinate (f);
% Midgrid boundary:
\path(1,{sqrt(3)/3})  coordinate (m1); \path(0,{2*sqrt(3)/3})  coordinate (m2);
\path(-1,{sqrt(3)/3})   coordinate (m3); \path(1,{-sqrt(3)/3}) coordinate (m6);
\def\midgrid{\draw[black!70, rounded corners=0pt] (0,0)--++(m1)--++(m2)--++(m1)--++(m2)--++(m1)--++(m2)--++(m1)--++(m2)--++(m1)--++(m2)
--++(m1)--++(m6)--++(m1)--++(m6)--++(m1)--++(m6)--++(m1)--++(m6)--++(m1)--++(m6)
(0,0)--++(m6)--++(m1)--++(m6)--++(m1)--++(m6)--++(m1)--++(m6)--++(m1)--++(m6)--++(m1)
--++(m2)--++(m1)--++(m2)--++(m1)--++(m2)--++(m1)--++(m2)--++(m1)--++(m2)--++(m1)};
\midgrid; \draw[black!70, rounded corners=0pt] (0,{10*{sqrt(3)}})--++(m6)--++(m1)--++(m6)--++(m1)--++(m6)--++(m1)--++(m6)--++(m1)--++(m6)--++(m1)
(15,{5*{sqrt(3)}})--++(m2)--++(m3)--++(m2)--++(m3)--++(m2)--++(m3)--++(m2)--++(m3)--++(m2)--++(m3);
\path(0,{2*sqrt(3)/3}) coordinate (m1);\path(-1,{sqrt(3)/3}) coordinate (m2);\path(1,{sqrt(3)/3}) coordinate (m6);
\midgrid;

\draw[blue,thick](0,0)--++(a)--++(c)--++(a)--++(e)--++(a)--++(c)--++(a)--++(e)--++(a)--++(c)--++(a)--++(e)
--++(a)--++(c)--++(a)--++(e)--++(a)--++(c)--++(a)--++(c)--++(e)--++(c)--++(e)--++(c)--++(a)--++(c)--++(a)
--++(e)--++(a)--++(c)--++(a)--++(e)--++(a)--++(c)--++(a)--++(c)--++(e)--++(c)--++(e)--++(c)--++(e)--++(c)
--++(e)--++(a)--++(e)--++(c)--++(e)--++(a)--++(e)--++(c)--++(e)--++(c)--++(a)--++(c)--++(a)--++(c)--++(a)
--++(c)--++(a)--++(e)--++(a)--++(c)--++(a)--++(e)--++(a)--++(c)--++(a)--++(e)--++(a)--++(c)--++(a)--++(e)
--++(a)--++(c)--++(a)--++(c)--++(e)--++(c)--++(a)--++(c)--++(e)--++(c)--++(a)--++(c)--++(e)--++(c)--++(a)
--++(c)--++(e)--++(c)--++(a)--++(c)--++(e)--++(c)--++(e)--++(c)--++(e)--++(c)--++(e)--++(a)--++(e)--++(a)
--++(c)--++(a)--++(e)--++(a)--++(c)--++(a)--++(e)--++(a)--++(e)--++(a)--++(e)--++(a)--++(e)--++(c)--++(e)
--++(c)--++(a)--++(c)--++(e)--++(c)--++(a)--++(c)--++(e)--++(c)--++(e)--++(a)--++(e)--++(a)--++(e)--++(c)
--++(e)--++(c)--++(a)--++(c)--++(e)--++(c)--++(a)--++(c)--++(e)--++(c)--++(a)--++(c)--++(e)--++(c)--++(a)
--++(c)--++(e)--++(c)--++(e)--++(a)--++(e)--++(c)--++(e)--++(a)--++(e)--++(c)--++(e)--++(a)--++(e)--++(c)
--++(e)--++(a)--++(e)--++(c)--++(e)--++(a)--++(e)--++(a)--++(c)--++(a)--++(c)--++(a)--++(e)--++(a)--++(e)
--++(c)--++(e)--++(a)--++(e)--++(c)--++(e)--++(a)--++(e)--++(a)--++(c)--++(a)--++(c)--++(a)--++(c)--++(a)
--++(c)--++(e)--++(c)--++(a)--++(c)--++(e)--++(c)--++(a)--++(c)--++(a)--++(e)--++(a)--++(e)--++(a)--++(e)
--++(a)--++(e)--++(c)--++(e)--++(a)--++(e)--++(c)--++(e)--++(a)--++(e)--++(c)--++(e)--++(a)--++(e)--++(c)--++(e)--++ (a)--++(e);
% De oorsprong:
\node[fill=red,circle, inner sep=1.4pt] at (0,0)  {};
% Eindpunt van tha(a):
\node[fill=black,circle, inner sep=1.0pt] at (15,{5*{sqrt(3)}})  {};
% Eindpunt van tha(ac):
\node[fill=black,circle, inner sep=1.0pt] at (0,{10*{sqrt(3)}})  {};
\end{tikzpicture}
\caption{\small The $\tha$-loop of a symmetric and an asymmetric planefilling 75-fold.\label{fig:loop75}}
\end{figure}
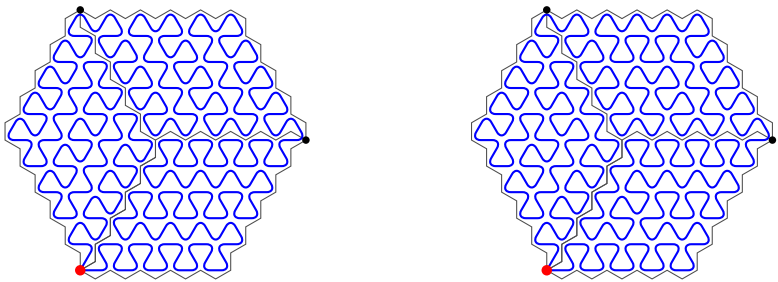

Theorem \ref{th:planefill} and \ref{th:folstarfol} will remain true with the obvious definitions of self-avoiding and planefilling, but the proof of the analogue of the existence Theorem \ref{exist} fails miserably. The analogue of Theorem \ref{exist} would be that there is a self-avoiding planefilling curve on the triangular lattice generated by a $m$-fold for every $m>1$ that can be written as
$$m=r^2+s^2-rs,\quad r,s \in \Z^+,$$
since $|z|=\sqrt{r^2+s^2-rs}$ is the norm of an Eisenstein integer $z=r+s\,\omega$, and we should have $|z|^2=m$. This set consists of the so called L\"oschian numbers larger than 1
$${\cal{L}}:=\{3,4,7,9,12,13,16,19,21,\dots\}.$$
The  reason that the proof fails is that not all paths in the directed triangular lattice correspond to allowed folding words (as, e.g., the word $aa$).
 In fact not only does the proof fail, but the result is also not true, since one can show  that there is no self-avoiding planefilling $m$-fold on the triangular lattice when $m$ is even. (The proof considers the union of the $\tha$-loop and the anti-$\tha$-loop: if we remove from this set the set $\rK[\tha(a)]$ then one obtains a simple curve that (as an application of Pick's Theorem shows) has $m-1$ Eisenstein integers in its interior. By symmetry of the midgrid boundaries of the the $\tha$-loop and the anti-$\tha$-loop, these integers have to be distributed evenly over the two loops, which is only possible if $m-1$ is even.)

Nevertheless, self-avoiding planefilling curves generated by folds opened at $60^o$ abound. For small $m$ these all have to be symmetric, but for large $m$ they can be asymmetric (see e.g., Figure~\ref{fig:loop75}, right). An infinite family is given, for example, by all $m$-folds with $m=3k^2$, $k=1,3,5,7,\ldots$. The folds here end in
$2k+k\,\omega$, and trace a paralellogram-shaped midgrid border in a zig-zag way (see Figure \ref{fig:loop75}, left for the $\tha$-loop of this fold when $k=5$.) Here the midgrid border is defined as in the $90^o$ case, by a diamond map $K_\diamond$, where now, for example, $K_\diamond[a]$ is the quadrangle  with vertices 0, $\frac12+\frac16\sqrt{3}i,\frac12-\frac16\sqrt{3}i$,  and 1.

\begin{figure}[!h]
\centering
\begin{tikzpicture}[scale=.08,rounded corners=0.3pt]
\path(2,0)  coordinate (a);
\path(1,{sqrt(3)})  coordinate (b);
\path(-1,{sqrt(3)}) coordinate (c);
\path(-2,0) coordinate (d);
\path(-1,{-sqrt(3)}) coordinate (e);
\path(1,{-sqrt(3)}) coordinate (f);
\draw[blue,thick] (0,20)--++(b)--++(a)--++(f)--++(a)--++(b)--++(c)--++(d)--++(c)--++(b)--++(c)--++(d)--++(e)--++(f)--++(e)--++(d)--++(e)--++(f)--++(a)--cycle;
% De oorsprong, en andere vouwpunten:
\node[fill=red,circle, inner sep=.8pt] at (0,20)  {};
\node[fill=black,circle, inner sep=.8pt] at (4,20)  {};
\node[fill=black,circle, inner sep=.8pt] at (6,{2*sqrt(3)+20})  {};
\node[fill=black,circle, inner sep=.8pt] at (4,{4*sqrt(3)+20})  {};
\node[fill=black,circle, inner sep=.8pt] at (0,{4*sqrt(3)+20})  {};
\node[fill=black,circle, inner sep=.8pt] at (-2,{2*sqrt(3)+20})  {};
\end{tikzpicture}\hspace*{2cm}
\begin{tikzpicture}[scale=.08,rounded corners=0.1pt]
\path(2,0)  coordinate (a);
\path(1,{sqrt(3)})  coordinate (b);
\path(-1,{sqrt(3)}) coordinate (c);
\path(-2,0) coordinate (d);
\path(-1,{-sqrt(3)}) coordinate (e);
\path(1,{-sqrt(3)}) coordinate (f);
\draw[blue,thick] (0,0)--++
(b)--++(a)--++(f)--++(a)--++(b)--++(c)--++(d)--++(c)--++(b)--++(a)--++(b)--++(c)--++(b)--++(a)--++(f)--++(e)--++(f)--++(a)--++(f)--++(e)--++(d)--++(e)--++ (f)--++(a)--++(b)--++(a)--++(f)--++(a)--++(b)--++(c)--++(b)--++(a)--++(f)--++(e)--++ (f)--++(a)--++(b)--++(a)--++(f)--++(a)--++(b)--++(c)--++(d)--++(c)--++(b)--++(c)--++(d)--++(e)--++ (d)--++(c)--++(b)--++(a)--++(b)--++(c)--++(d)--++(c)--++(b)--++(c)--++(d)--++(e)--++(f)--++(e)--++(d)--++(c)--++(d)--++(e)--++(d)--++ (c)--++(b)--++(a)--++(b)--++(c)--++(b)--++(a)--++(f)--++(a)--++(b)--++(c)--++(d)--++(c)--++(b)--++(a)--++(b)--++(c)--++(b)--++(a)--++(f)--++(e)--++ (f)--++(a)--++(b)--++(a)--++(f)--++(a)--++(b)--++(c)--++(d)--++(c)--++(b)--++(c)--++(d)--++(e)--++(d)--++(c)--++(b)--++(a)--++ (b)--++(c)--++(b)--++(a)--++(f)--++(a)--++(b)--++(c)--++(d)--++(c)--++(b)--++(a)--++(b)--++(c)--++(b)--++(a)--++(f)--++(e)--++(f)--++(a)--++(f)--++(e)--++(d)--++(e)--++ (f)--++(a)--++(b)--++(a)--++(f)--++(e)--++(f)--++(a)--++(f)--++(e)--++(d)--++(c)--++(d)--++(e)--++(f)--++(e)--++(d)--++(e)--++ (f)--++(a)--++(b)--++(a)--++(f)--++(a)--++(b)--++(c)--++(b)--++(a)--++(f)--++(e)--++(f)--++(a)--++(f)--++(e)--++(d)--++(e)--++(f)--++(a)--++(b)--++(a)--++(f)--++(e)--++ (f)--++(a)--++(f)--++(e)--++(d)--++(c)--++(d)--++(e)--++(d)--++(c)--++(b)--++(c)--++(d)--++(e)--++(f)--++(e)--++(d)--++(e)--++(f)--++(a)--++(f)--++(e)--++ (d)--++(c)--++(d)--++(e)--++(f)--++(e)--++(d)--++(e)--++(f)--++(a)--++(b)--++(a)--++(f)--++ (a)--++(b)--++(c)--++(b)--++(a)--++(f)--++(e)--++(f)--++(a)--++(b)--++(a)--++ (f)--++(a)--++(b)--++(c)--++(d)--++(c)--++(b)--++(a)--++(b)--++(c)--++(b)--++(a)--++(f)--++(e)--++(f)--++(a)--++(f)--++(e)--++(d)--++(e)--++(f)--++(a)--++(b)--++(a)--++(f);
% De oorsprong:
\node[fill=red,circle, inner sep=1.2pt] at (0,0)  {};
\end{tikzpicture}
\caption{ \small Folding in three according to $S=DD$, and opening at $120^o$ degree angles. Shown are the $\tha$-loop and the $5^{\rm th}$ order curve.\label{fig:sier}}
\end{figure}
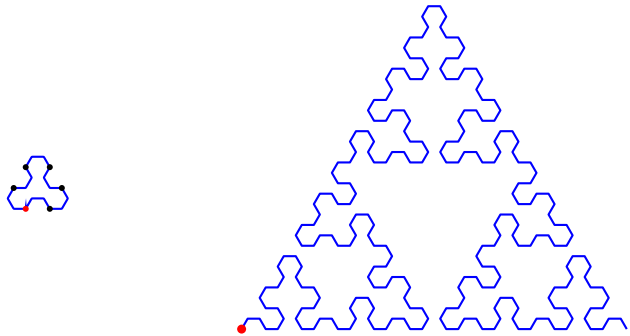

Curves on the hexagonal lattice are obtained when we unfold at $120^o$ degree angles. We need a 6 letter alphabet, say $\{a,b,c,d,e,f\}$, which corresponds to the $6^{\rm th}$ roots of unity. Theorem \ref{th:avoid} remains true, when we define the $\tha$-loop in the usual way. See Figure \ref{fig:sier} for how the 3-fold $S=DD$ generates a self-avoiding curve on the hexagonal lattice. Actually the limiting fractal is an old friend: it is easily shown that the curves $L_\tha^{-n}K[\tha^n(a)]$ are exponentially close in the Hausdorff metric to the usual approximants of the Sierpinski triangle, also known as the Sierpinski gasket. The notion of planefilling is problematic in the hexagonal case, since the grid points  all have degree 3.

The concept of self-avoiding and planefilling curves generated by morphisms on the 4 letter alphabet can be generalized in an other direction. In fact, a slight generalization of Theorems \ref{th:planefill} and \ref{th:folstarfol} will hold for morphisms that preserve $\DD$-words (see \cite{FoldsIII}). From that paper we have the following example of a non-uniform morphism that generates a self-avoiding and planefilling curve:
$$ \tha(a)=abad, \quad  \tha(b)=ab,\quad  \tha(c)=cb,\quad  \tha(d)=cdcd.$$

\subsection*{Acknowledgements}
I am grateful for the enthusiasm and encouragement of Donald Knuth.
Moreover, his criticism and valuable suggestions have greatly improved this paper.
I am also grateful to my young friend Derong Kong for his critical reading of preliminary versions of this paper; his comments have been very useful.

\end{document}